\documentclass[11pt,reqno]{amsproc}
\allowdisplaybreaks
\usepackage{color,soul}
\usepackage[table,xcdraw,dvipsnames]{xcolor}
\usepackage{float}
\usepackage{fullpage}
\usepackage{rotating}
\usepackage{stmaryrd}
\usepackage{subfigure}
\usepackage{listings}
\usepackage{enumerate}
\usepackage[small]{caption}
\usepackage{morefloats}
\usepackage{stmaryrd}
\usepackage{mathtools}
\usepackage{cancel}
\usepackage{siunitx}
\usepackage{algorithm}
\usepackage{algorithmicx}
\usepackage{algpseudocode}

\numberwithin{equation}{section}
\usepackage{graphicx}
\DeclareGraphicsExtensions{.eps,.png}
\usepackage{epstopdf}
\usepackage{array,multirow}
\usepackage{rotating}
\usepackage{booktabs,tabularx}
\usepackage{libertine}
\usepackage[font=small,labelfont=bf,tableposition=top]{caption}
\usepackage[utf8x]{inputenc}
\usepackage[english]{babel}
\usepackage{multicol}
\usepackage[flushleft]{threeparttable}
\makeatother
\usepackage[semicolon,square,authoryear]{natbib}
\usepackage[debug=false, colorlinks=true, pdfstartview=FitV, 
linkcolor=blue, citecolor=blue, urlcolor=blue]{hyperref}

\usepackage{makecell, boldline}
\usepackage{url}
\usepackage{soul}
\usepackage{morefloats}

\definecolor{amethyst}{rgb}{0.6, 0.4, 0.8}
\definecolor{burgundy}{rgb}{0.5, 0.0, 0.13}
\newcommand\bfx{{\bf x}}

\title{\textbf{Bound-preserving discontinuous Galerkin methods for compressible two-phase flows in porous media}}
\author{
    \textbf{M.~S.~Joshaghani}$^\uparrow$$^\circ$\footnote{$^\uparrow$ Department of Computational Applied Mathematics
and Operations Research, Rice University, Houston, TX}\footnote{$^\circ$ Research conducted while at Rice University}
    and
    \textbf{B.~Riviere}$^\uparrow$ $^\ast$\footnote{$^\ast$ Corresponding author: riviere@rice.edu}
}
\newsavebox{\measurebox}
\begin{document}
\date{\today}
\setcounter{figure}{0}
\begin{abstract}
This paper presents a numerical study of immiscible, compressible two-phase flows in porous media, that takes into 
account heterogeneity, gravity, anisotropy and injection/production wells. We formulate a fully implicit stable 
discontinuous Galerkin solver for this system that is accurate, that respects maximum principle for the approximation of saturation,
 and that is locally mass conservative.  To completely eliminate the overshoot and undershoot phenomena,
we construct a flux limiter that produces bound-preserving elementwise average of the saturation. The addition
of a slope limiter allows to recover a pointwise bound-preserving discrete saturation.  Numerical results
show that both maximum principle  and monotonicity of the solution are satisfied. The proposed flux limiter
does not impact the local mass error and the number of nonlinear solver iterations.
\end{abstract}
\maketitle

\section{Introduction}

Compressible multiphase flows in porous media occur in many applications such
as subsurface carbon sequestration. Compressibility is modeled
by the dependence of the fluid mass densities and rock porosity on the phase pressure.
This work formulates a stable discontinuous Galerkin method for compressible two-phase flows
in heterogeneous porous media; the main contribution being that the scheme satisfies a maximum principle
for the numerical phase saturation.

The literature on numerical methods for two-phase flows in porous media is large, particularly for the case
of incompressible phases \citep{chen2006computational,HoteitFiroozabadi2008,Bastian2014,HouChenSunChen2016,doyle2020multinumerics}. It is known that suitable  methods for porous media flows 
should satisfy a local mass conservation property.  Both finite
volume methods and discontinuous Galerkin methods are good candidates. There are other desirable
properties such as maximum principle and monotonicity.  On the one hand, finite volume methods produce piecewise constant approximation
of the saturation that satisfies physical bounds \citep{Michel03,Droniou2014,ghilani2019positive}. On the other hand, finite volume methods are numerically diffusive
and require Voronoi-type grids for unstructured meshes, which can be challenging to construct for anisotropic 
heterogeneous media \citep{aavatsmark2002introduction,de2007node,contreras2021non}.  Discontinuous Galerkin (DG) methods overcome the shortcomings of finite volume
methods because they belong to the class of variational problems like finite element methods. DG methods
can be of arbitrary order, are adapted to any unstructured meshes,
 and in the case of convection-dominated problems, they produce sharp fronts with negligible
numerical diffusion. However in the neighborhood of the saturation front, overshoot and undershoot phenomena may occur
as the maximum principle for the DG solution is not guaranteed \citep{KlieberRiviere2006,epshteyn2007fully,ErnMozolevskiSchuh2010,Bastian2014,JameiGhafouri2016}. These overshoot and undershoot
phenomena remain bounded throughout the simulation and it  is possible to
reduce the amount of overshoot/undershoot by mesh refinement, or by projecting phase velocities into H(div) conforming spaces,
by varying the penalty parameters, or by using slope limiters \citep{Hoteit2004,Krivodonova2007,kuzmin2010vertex,kuzmin2013,kuzmin2012flux}.  However, a complete elimination of the overshoot/undershoot
has been challenging to achieve.

Recently in \citep{JoshaghaniRiviereSekachev}, we proposed a DG method combined with a flux limiter for
solving the immiscible incompressible two-phase flows in porous media. The DG saturations are shown
to satisfy a maximum principle, in the sense that solutions do not exhibit any overshoot and undershoot phenomena. 
This current work is an extension of \citep{JoshaghaniRiviereSekachev}
to the case of compressible phases.  This is a more complicated problem
because of the dependency of the coefficients (phase densities and porous medium porosity) with respect to the pressure.
We have observed that 
the amount of overshoot and undershoot in DG solutions is larger for compressible flows than for incompressible flows.
The numerical method is fully implicit and the nonlinear equations are solved by Newton's method.
The novel contribution is the construction of a new flux limiter that takes into account the dependence of
the densities and porosity on the unknown. The proposed flux limiter is related to flux-corrected transport algorithms
for the solution of conservation laws \citep{frank2019bound,kuzmin2012flux}.
We show that the
resulting method respects the maximum principle and is locally mass conservative for several problems
taking into account gravity and heterogeneity. To our knowledge, this work is the first to present
a DG-based scheme for compressible two-phase flows, that does not violate the maximum principle. 
An outline of the paper is as follows: after a brief introduction of the model equations in Section~\ref{Sec:S1_GEs},
the proposed numerical method is formulated in Section~\ref{sec:method}. Numerical results and conclusions follow.

\section{Governing equations}
\label{Sec:S1_GEs}

The mathematical model describing the flow of a wetting phase (with saturation $s$ and pressure $p$) 
and a non-wetting phase  in a domain $\Omega\subset\mathbb{R}^2$ over the time interval $[0,T]$ is:
\begin{align}
    \label{Eqn:BoM_pres}
	&\frac{\partial}{\partial t} \left(\phi(p)  \rho_\ell(p) (1-s)\right)
	-\nabla \cdot \left(\rho_\ell(p)  \lambda_\ell(s) K (\nabla p - \rho_\ell(p)   \mathbf{g})\right) 
	= \rho_\ell(p) q_\ell,\\
    \label{Eqn:BoM_sat}
	&
	\frac{\partial}{\partial t} \left(\phi(p) \rho_w(p) s\right)
	-\nabla \cdot \left(\rho_w(p) \lambda_w(s) K (\nabla p - \rho_w(p)   \mathbf{g})\right) = \rho_w(p) q_w.
\end{align}
Since the capillary pressure is neglected, the mathematical model is a system of nonlinear hyperbolic equations.
Compressibility of the phases and the medium is modeled by the following dependence of densities and porosity on the pressure:
\begin{align*}
\phi(p) = \phi^0(1+c_r p), 
\quad \rho_\ell(p) = \rho_\ell^0(1+c_\ell p), 
\quad \rho_w(p) = \rho_w^0(1+c_w p),
\end{align*}
where the rock and fluid compressibilities, $c_r, c_\ell, c_w$ and the reference porosity and densities, $\phi_0, \rho_w^0, \rho_\ell^0$ are given constants.  The absolute permeability, $K$,  of the medium is either a positive scalar or a symmetric positive definite matrix $K$
that may vary in space. The phase mobilities are $\lambda_\ell$ and $\lambda_w$ for the non-wetting phase and
wetting phase respectively; they are given functions of saturation and also depend on the phase viscosities $\mu_\ell, \mu_w$. 
In this work, the commonly used
Brooks-Corey model is considered:
\begin{align}
    \lambda_{w}(s_e)=\frac{s_e^{2}}{\mu_w}, 
    \quad
    \lambda_{\ell}(s_e) = \frac{(1-s_e)^2}{\mu_\ell},
\end{align}
where effective saturation is defined as:
\begin{align}
    s_e = \frac{s-s_{rw}}{1-s_{rw}-s_{r\ell}}.
\end{align}
The residual saturation for wetting phase and non-wetting phase are denoted by $s_{rw}$ and $s_{r\ell}$ respectively. 
The functions $q_\ell$ and $q_w$ are given source/sink functions.
The boundary is partitioned into $\partial\Omega =\Gamma^{\mathrm{D},p}\cup\Gamma^{\mathrm{N},p} = \Gamma^{\mathrm{D},s}\cup\Gamma^{\mathrm{out}}\cup\Gamma^{\mathrm{N},s}$.
We prescribe Dirichlet and flux boundary conditions on $\Gamma^{\mathrm{D},p}\cup\Gamma^{\mathrm{D},s}$ and
$\Gamma^{\mathrm{N},p}\cup\Gamma^{\mathrm{N},s}$ respectively, as follows: 
\begin{alignat*}{2}
	& p = g^p 
	&& \quad  \mathrm{on} \; \Gamma^{\mathrm{D},p}, \\
	& s = g^s 
	&& \quad \mathrm{on} \; \Gamma^{\mathrm{D},s}, \\
    & \rho_{\ell}(p) \lambda_{\ell}(p)K (\nabla p-\rho_{\ell}(p)  \mathbf{g})\cdot {\mathbf{n}} = j^p
	&& \quad\mathrm{on} \; \Gamma^{\mathrm{N},p}, \\
    & \rho_w(p) \lambda_w(p) K(\nabla p-\rho_w(p) \mathbf{g})\cdot {\mathbf{n}} = j^s
	&&\quad \mathrm{on} \; \Gamma^{\mathrm{N},s}. 
\end{alignat*}
We will consider the case of flows driven by boundary conditions and the case of flows driven by wells (source/sink functions).
For the latter, only homogeneous Neumann boundary conditions are imposed on the boundary.  The source/sink functions 
depend on the saturation as follows
\begin{align}
\label{eq:sourcesink}
    q_{\alpha}(s)=f_{\alpha}(s_{\mathrm{in}})\bar{q}-f_{\alpha}(s)\underline{q}, \quad \alpha=w,\ell,
\end{align}
where $s_\mathrm{in}$ is the injected saturation value, $\bar{q}$ and $\underline{q}$ are the injection and production
well flow rates respectively, and $f_\alpha$ is the fractional flow for phase $\alpha$.  
The fractional flows are related to the mobilities by:
\begin{align}
    f_w(s) = \frac{\lambda_w(s)}{\lambda_w(s)+\lambda_{\ell}(s)} \quad \mbox{and} \quad
    f_\ell(s) = 1 - f_w(s).
\end{align}

For the case of flows driven by boundary conditions, we assume
that the Dirichlet boundary for the saturation is strictly included in the Dirichlet boundary for the pressure and
the outflow boundary is the complement $\Gamma^\mathrm{out} = \Gamma^{p,\mathrm{D}}\setminus\Gamma^{s,\mathrm{D}}$.  No boundary conditions
are assumed for the saturation on the outflow boundary. The source/sink functions are set to zero.

Finally, the initial pressure and saturation are denoted by $p_0$ and $s_0$.

\section{Numerical Method}
\label{sec:method}

We discretize \eqref{Eqn:BoM_pres}-\eqref{Eqn:BoM_sat} by a fully implicit interior penalty discontinuous Galerkin method. We first set some notation.
The domain $\Omega$ is decomposed into a non-degenerate partition $\mathcal{E}_h=\{E\}_E$
consisting of $N_h$ triangular elements of maximum diameter $h$.   
Let $\Gamma_h$ denote the set of all edges and $\Gamma_h^{i}$ denote the set of interior edges.
For any $e \in \Gamma_h^{i}$, fix a unit normal vector $\mathbf{n}_e$ and denote by $E^+$ and $E^-$ 
the elements that share the edge $e$ such that the vector
$\mathbf{n}_e$ is  directed from $E^+$ to $E^-$.
We define the jump and average of a scalar function $\xi$ on $e$ as follows:
\begin{align}
   \llbracket
          \xi 
   \rrbracket = \xi\vert_{E^+}-\xi\vert_{E^-},\quad
   \{\!\!
   \{
    \xi       
   \}\!\!
   \} = \frac{1}{2}\left( \xi\vert_{E^+}+\xi\vert_{E^-}\right). 
\end{align}
By convention, if $e$ belongs to the boundary $\partial \Omega$, then the jump and average of $\xi$ on $e$ coincide 
with the trace of $\xi$ on $e$ and the normal vector $\mathbf{n}_e$ coincides with the outward normal $\mathbf{n}$.
Let $\mathbb{P}_1(E)$ be the space of linear polynomials  on an element $E$. The discontinuous finite
element space of order one is:
\begin{align}
    \mathcal{D}(\mathcal{E}_h) = \left\{\xi\in L^2(\Omega): \xi\vert_E\in\mathbb{P}_1(E), \,  \forall E \in \mathcal{E}_h\right\}.
\end{align}
The time interval $T$ is divided into $N_\tau$ equal subintervals of length $\tau$. Let $P_{n}$ and $S_{n}$ denote the numerical 
solutions at time $t_n$. The proposed discontinuous Galerkin scheme for equations~\eqref{Eqn:BoM_pres}--\eqref{Eqn:BoM_sat} reads: 
Given $(P_{n},S_{n})\in \mathcal{D}(\mathcal{E}_{h}) \times \mathcal{D}(\mathcal{E}_{h})$, find
$(P_{n+1},S_{n+1})\in \mathcal{D}(\mathcal{E}_{h}) \times \mathcal{D}(\mathcal{E}_{h})$ such that for all $\xi\in\mathcal{D}(\mathcal{E}_h)$: \\
\begin{align}
    &
    \frac{1}{\tau} \sum_{E\in\mathcal{E}_h}\int_E 
	\phi(P_{n+1}) \rho_{\ell} (P_{n+1}) (1-S_{n+1})  \xi
  +
  \sum_{E\in\mathcal{E}_h}\int_E 
  \rho_{\ell}(P_{n+1}) \lambda_{\ell}(S_{n+1})K \big(\nabla P _{n+1} -\rho_{\ell}(P_{n+1})  \mathbf{g}\big) \cdot \nabla \xi \nonumber\\
  &-
  \sum_{e\in\Gamma^{i}_h} \int_e
  (\lambda_{\ell}(S_{n+1}) )^{\uparrow \mathbf{v}_{\ell}^{n}} \, 
  \{\!\!\{
          \rho_{\ell}(P_{n})K\big(\nabla P_{n+1} -\rho_{\ell}(P_{n+1})  \mathbf{g}\big)  \cdot   {\mathbf{n}_e}
  \}\!\!\} \llbracket \xi \rrbracket \nonumber\\
  &-
  \sum_{e\in\Gamma^{\mathrm{D},p}} \int_e
          \lambda_{\ell}(S_{n+1}) K
          {\rho_{\ell}(P_{n+1})\big(\nabla P_{n+1} -\rho_{\ell}(P_{n+1})  \mathbf{g}\big) }   \cdot     {\mathbf{n}_e}
    \; \xi  
   +
      \sum_{e\in\Gamma_h \backslash \Gamma^{\mathrm{N},p }}  \frac{\sigma}{h} \int_e
  \llbracket P_{n+1} \rrbracket \llbracket \xi \rrbracket 
\nonumber\\
  & =  \sum_{E\in\mathcal{E}_h}\int_E 
 \rho_{\ell}(P_n) q_{\ell}(S_n) \xi  
  +\frac{1}{\tau} \sum_{E\in\mathcal{E}_h}\int_E 
 \phi(P_{n}) \rho_{\ell} (P_{n}) (1-S_{n})   \xi
  +\sum_{e\in \Gamma^{\mathrm{D},p}}  \frac{\sigma_p}{h} \int_e
  {g^p}  \xi 
 +\sum_{e\in\Gamma^{\mathrm{N},p}}   \int_e  
 j^{p} \xi, 
    \label{Eqn_presWeak}
\end{align}
and 
\begin{align}
	&  
\frac{1}{\tau}    \sum_{E\in\mathcal{E}_h}\int_E 
	\phi(P_{n+1}) \rho_{w} (P_{n+1}) S_{n+1}  \xi
    +
    \sum_{E\in\mathcal{E}_h}\int_E 
    \rho_{w}(P_{n+1}) \lambda_{w}(S_{n+1})K \big(\nabla P _{n+1} -\rho_{w}(P_{n+1})  \mathbf{g}\big) \cdot \nabla \xi \nonumber\\
  &-
  \sum_{e\in\Gamma_h^i} \int_e
  (\lambda_{w}(S_{n+1}) )^{\uparrow \mathbf{v}_w^n}\,
  \{\!\!\{
          \rho_{w}(P_{n+1})K \big(\nabla P_{n+1} -\rho_{w}(P_{n+1})  \mathbf{g}\big) \cdot  \mathbf{n}_e
  \}\!\!\} \llbracket \xi \rrbracket \nonumber\\
  &-
\sum_{e\in \Gamma^{\mathrm{D},s}} \int_e
\lambda_{w}(g^s)K
\rho_{w}(g^p)\big(\nabla P_{n+1} -\rho_{w}(P_{n+1})  \mathbf{g}\big)  \cdot \mathbf{n}_e \;
    \xi   \nonumber\\
  &-
  \sum_{e\in \Gamma^{\mathrm{out}}} \int_e
  \lambda_{w}(S_{n+1})K
\rho_{w}(g^p)\big(\nabla P_{n+1} -\rho_{w}(P_{n+1}) \mathbf{g}\big)  \cdot  \mathbf{n}_e \;
    \xi   
   + \sum_{e\in\Gamma_h \backslash \Gamma^{\mathrm{N},s }}  \frac{\sigma}{h} \int_e
  \llbracket S_{n+1} \rrbracket \llbracket \xi \rrbracket \nonumber\\
  &  =\sum_{E\in\mathcal{E}_h}\int_E 
  \rho_{w}(P_n) q_{w}(S_n) \xi
	+\frac{1}{\tau} \sum_{E\in\Gamma_h}\int_E 
    \phi(P_n)\rho(P_n) S_n \xi
    +\sum_{e\in \Gamma^{\mathrm{D},s}}  \frac{\sigma_s}{h} \int_e
    g^s  \xi 
    + \sum_{e \in \Gamma^{\mathrm{N},s}} \int_e
    j^s \xi.
    \label{Eqn_satWeak}
\end{align}

The penalty parameter $\sigma$ is constant on the interior edges and its value is chosen $10$ times larger on the Dirichlet boundaries.
The quantities $(\cdot)^{\uparrow\mathbf{v}_\ell^n}$ and $(\cdot)^{\uparrow\mathbf{v}_w^n}$ denote
the upwind values with respect to the vector functions $\mathbf{v}_\ell^n$ and $\mathbf{v}_w^n$ that are scaled quantities of the phase velocities. They depend on the pressure and saturation evaluated at the previous time $t_n$: 
\[
    \label{Eqn:Velocities}
\mathbf{v}_\ell^n  =  -\rho_{\ell}(P_n)K\big(\nabla P_{n} -  \rho_{\ell}(P_n)\mathbf{g}\big), \quad
\mathbf{v}_w^n  =  -\rho_{w}(P_n)K\big(\nabla P_{n} -  \rho_{w}(P_n)\mathbf{g}\big).
\]
The definition of the upwind operator with respect to a generic discontinuous vector field $\mathbf{v}$ is:
\[
\forall e  = \partial E^+\cap \partial E^-, \quad
\xi^{\uparrow\mathbf{v}}|_e = \left\{
\begin{array}{c}
\xi\vert_{E^+}, \quad\mbox{if} \quad \{\!\!  \{ \mathbf{v} \}\!\!  \}\cdot\mathbf{n}_e >0,\\
\xi\vert_{E^-}, \quad\mbox{if} \quad \{\!\!  \{ \mathbf{v} \}\!\!  \}\cdot\mathbf{n}_e \leq 0.
\end{array}
\right.
\]
At the initial time, the discrete saturation and pressure are the $L^2$ projection of the initial conditions.
\[
\int_\Omega S_0 v = \int_\Omega s_0 v, \quad 
\int_\Omega P_0 v = \int_\Omega p_0 v, \quad 
\forall v\in\mathcal{D}(\mathcal{E}_h).
\]

To solve the nonlinear system \eqref{Eqn_presWeak}-\eqref{Eqn_presWeak}, we use Newton's method.
Let the superscript $(i)$ denote the current Newton iteration. 
We solve for the updates $\delta P$ and $\delta S$ at each iteration:
\begin{align*}
    P^{(i+1)}_{n+1} = P^{(i)}_{n+1} + \delta P, \quad
    S^{(i+1)}_{n+1} = S^{(i)}_{n+1} + \delta S.
\end{align*}
Once the Newton iterations converge, we apply the flux and slope limiters described in the next section (see Algorithm~\ref{alg}).
The novelty of this work is in the formulation of the flux limiters described in details in Section~\ref{Sub:Flux_limiter}.
For the slope limiter, we employ the vertex-based slope limiter introduced by \citep{kuzmin2010vertex}.

\begin{algorithm}
\caption[Scheme]{DG+FL+SL method}\label{alg}
\begin{algorithmic}
\State Compute initial saturation $S_0$ and pressure $P_0$
\For{$n=0,\dots, (N_\tau-1)$}
\State Apply Newton solver to obtain $(P_{n+1},S_{n+1})$
\State Apply flux limiter: $S_{n+1}^\mathrm{FL} = \mathcal{L}_\mathrm{flux}(S_{n+1})$
\State Apply slope limiter:  $S_{n+1} = \mathcal{L}_\mathrm{slope}(S_{n+1}^\mathrm{FL})$
\EndFor
\end{algorithmic}
\end{algorithm}

\subsection{Flux limiter}
\label{Sub:Flux_limiter}
The flux limiter will enforce that the element-wise average of the saturation satisfies the desired physical bounds. 
It is applied every time step and we assume that the saturation at the previous time step, $t_n$,  satisfies a maximum principle:
\begin{align}
    s_\ast \le S_n(\mathbf{x}) \le s^\ast, \quad \forall \mathbf{x} \in \Omega.
\label{eq:physbounds}
\end{align}
for some constants $0\leq s_\ast\leq s^\ast\leq 1$; these constants depend on the residual saturations,
namely  $s_\ast = s_{rw}$ and $s^\ast = 1-s_{r\ell}$. 
The flux limiting is applied to each element $E$ given the element-wise average of the saturation at the previous and current time steps and
given a flux function defined on each face $e\subset\partial E$. We denote the element-wise
average of the saturation at time $t_n$ and $t_{n+1}$, by $\overline{S_n}$ and $\overline{S_{n+1}}$ defined by:
\begin{equation}
\label{eq:defavg}
\overline{S_{i}}|_E = \overline{S_{i,E}}, \quad
\overline{S_{i,E}}= \frac{1}{|E|}\int_{E} S_{i}, 
\quad \forall E\in\mathcal{E}_h, \quad i=n,n+1.
\end{equation}
Next, for a fixed element $E$, let $\mathbf{n}_E$ be the unit normal vector outward to $E$. We 
define the flux function $\mathcal{H}_{n+1}|_E = \mathcal{H}_{n+1,E}$ as follows:
   \begin{align*}
           \forall e = \partial E\cap \partial E', \quad \mathcal{H}_{n+1,E}(e) =&
          - \int_e (\lambda_{w}(S_{n+1}) )^{\uparrow\mathbf{v}_w^n}
          \{\!\!
          \{
              \rho_w(P_{n+1}) K \big(\nabla P_{n+1} -\rho_{w}(P_{n+1})  \mathbf{g}\big) \cdot  {\mathbf{n}}_{E}       
          \}\!\!
          \}
            \nonumber\\
          & + \frac{\sigma}{h} \int_e (S_{n+1}|_E-S_{n+1}|_{E'}) \\
           \forall e \in\partial E \cap \Gamma^{\mathrm{D},s}, \quad \mathcal{H}_{n+1,E}(e)& =
           - \int_e \lambda_{w}(g^s) K \rho_{w}(g^p) \big(\nabla P_{n+1} -\rho_{w}(P_{n+1})  \mathbf{g}\big) \cdot {\mathbf{n}}_{E} 
           + \frac{\sigma}{h} \int_e (S_{n+1} - g^s),\\
           \forall e \in\partial E \cap \Gamma^{\mathrm{N},s}, \quad \mathcal{H}_{n+1,E}(e) &= \int_e j^s,\\
           \forall e \in\partial E \cap \Gamma^{\mathrm{out}}, 
           \quad \mathcal{H}_{n+1,E}(e) &= \int_e \lambda_w(S_{n+1}) K \rho_{w}(g^p)(\nabla P_{n+1}-\rho_w(P_{n+1}) \mathbf{g})\cdot\mathbf{n}_E.
\end{align*}
 For an interior face $e$ of the element $E$, the quantity $\mathcal{H}_{n+1,E}(e)$ measures the net mass flux 
across $e$ into the neighboring element $E'$ that also shares the face $e$. We note that: 
\[\mathcal{H}_{n+1,E}(e) = - \mathcal{H}_{n+1,E'}(e).
\]
The flux limiter updates the saturation in each mesh element such that its new element-wise average satisfies
the maximum principle \eqref{eq:physbounds}.
\begin{equation}
\mathcal{L}_\mathrm{flux}(S_{n+1})(\bfx) 
= S_{n+1}^{\mathrm{FL}}(\bfx) = 
S_{n+1}(\bfx)-\overline{S_{n+1,E}}+\bar{S}^{\mathrm{FL}}_{n+1}|_E, \quad \forall \bfx \in E.
\label{eq:updatecellaverage}
\end{equation}
The new cell-average of the saturation is obtained by an iterative process, that takes for input the cell average at the previous time
step and the flux function:
\[
\bar{S}^{\mathrm{FL}}_{n+1} = \mathcal{L}_\mathrm{avg}(\overline{S_{n}},\mathcal{H}_{n+1}).
\]

Next, we describe the algorithm for the operator $\mathcal{L}_\mathrm{avg}$.
For a fixed element $E$, we  denote by $\mathcal{N}_E$ the set of elements that include $E$ and all neighboring elements $E'$ that share a face $e$ with $E$. 
The algorithm constructs a sequence of flux functions and element-wise averages for $E$ and its neighbors $E'$.   While the construction
of the element-wise averages are local to $E$ and its neighbors $E'$, the stopping criterion is global to ensure bound-preserving solutions.
We first initialize the sequences with the input arguments:
\[
\bar{S}_{\tilde{E}}^{(0)} =\overline{S_{n,\tilde{E}}},\quad \mathcal{H}_{\tilde{E}}^{(0)} = \mathcal{H}_{n+1,\tilde{E}},
\quad \forall \tilde{E} \in \mathcal{N}_E.
\]
Next, for $k\geq 1$, we have the following steps:
\begin{itemize}
\item[Step~1.] Compute inflow and outflow fluxes:
        \begin{align}
        P_{\tilde{E}}^+ = \tau \sum_{e \in \partial \tilde{E}}\max(0,-\mathcal{H}_{\tilde{E}}^{(k-1)}(e)),
        \quad
        P_{\tilde{E}}^- = \tau \sum_{e \in \partial \tilde{E}}\min(0,-\mathcal{H}_{\tilde{E}}^{(k-1)}(e)), \quad \forall \tilde{E} \in \mathcal{N}_E.
        \end{align}
\item[Step~2.] Compute admissible upper and lower bounds for all $\tilde{E}\in \mathcal{N}_E$:
        \begin{align}\label{eq:defQplus}
            Q_{\tilde{E}}^+ =& |\tilde{E}|\left({\color{black}
            \overline{\rho_w(P_{n+1})\phi(P_{n+1})}
            }s^{*}
            -{\color{black}
            \overline{\rho_w(P_{n})\phi(P_{n})}
        }\bar{S}_{\tilde{E}}^{(k-1)} \right)\\
&- |\tilde{E}| \gamma_{1k}\tau \left(
            {\color{black}
                \overline{\rho_w(P_n) f_w(s_{\mathrm{in}})}\bar{q}_{\tilde{E}}
+\overline{\rho_w(P_n)} f_w(\bar{S}_{\tilde{E}}^{(k-1)}) \underline{q}_{\tilde{E}}
            } \right),\nonumber\\
            Q_{\tilde{E}}^- =& |\tilde{E}|\left({\color{black}
            \overline{\rho_w(P_{n+1})\phi(P_{n+1})}
            }s_{*}
            -{\color{black}
            \overline{\rho_w(P_{n})\phi(P_{n})}
        }\bar{S}_{\tilde{E}}^{(k-1)}\right)\\
            & -\vert \tilde{E}\vert \gamma_{1k}\tau \left(
            {\color{black}
                \overline{\rho_w(P_n) f_w(s_{\mathrm{in}})}\bar{q}_{\tilde{E}}
+\overline{\rho_w(P_n)} f_w(\bar{S}_{\tilde{E}}^{(k-1)})
                \underline{q}_{\tilde{E}}
            } \right).\nonumber\\
        \end{align}
The quantity        $Q_{\tilde{E}}^{+}$ measures the amount of mass that can be stored in element $\tilde{E}$ without creating a mean-value overshoot. 
        Similarly, $Q_{\tilde{E}}^{-}$ is a measure for the amount of mass that should be removed from element $\tilde{E}$ without creating a mean-value undershoot.
The scalar factor $\gamma_{1k}$ is equal to $1$ if $k=1$ and $0$ otherwise.
The injection and production well rates, restricted to any element $\tilde{E}$, are denoted by $\bar{q}_{\tilde{E}}$ and $\underline{q}_{\tilde{E}}$ respectively.
They are assumed to be piecewise constant fields; otherwise we take the element-wise average of the flow rates.

\item[Step~3.]   Compute limiting factors $\alpha_E^{(k-1)}(e)$ for all faces $e \subset \partial E$.
If $e$ is an interior face such that $e=\partial E\cap\partial E'$:
\[
\alpha_E^{(k-1)}(e) = \left\{
\begin{array}{cc}
\min\left( \min(1,Q^+_E/P^+_E), \min(1,Q^-_{E'}/P^-_{E'}) \right) & \mbox{if} \quad \mathcal{H}_E^{(k-1)}(e) < 0, \\
\min\left( \min(1,Q^-_E/P^-_E), \min(1,Q^+_{E'}/P^+_{E'}) \right) & \mbox{if} \quad \mathcal{H}_E^{(k-1)}(e) > 0.
\end{array}
\right.
\]
If $e$ is a boundary face:
\[
\alpha_E^{(k-1)}(e) = \left\{
\begin{array}{cc}
\min(1,Q^+_E/P^+_E) & \mbox{if} \quad \mathcal{H}_E^{(k-1)}(e) < 0, \\
\min(1,Q^-_E/P^-_E) & \mbox{if} \quad \mathcal{H}_E^{(k-1)}(e) > 0.
\end{array}
\right.
\]
\item[Step~4.]  Update $\bar{S}_{E}^{(k)}$ and $\mathcal{H}_E^{(k)}$ as follows:
                  \begin{align}
                      \bar{S}_{E}^{(k)}
                      =&
                      {\color{black}
                      \frac{\overline{\rho_w(P_{n})\phi(P_{n})} }
                     {\overline{\rho_w(P_{n+1})\phi(P_{n+1})} }
                         }
                      \bar{S}_{E}^{(k-1)}
                      -
                      \frac{\tau}{
                          {\color{black}
                              \overline{\rho_w(P_{n+1})\phi(P_{n+1})}
                      }|E|
                      }
                      \sum_{e\subset\partial E}\alpha_E(e)^{(k-1)}
                      \mathcal{H}_E^{(k-1)}(e)
                      \nonumber\\
                      &
                      {\color{black}
                          {\color{black}+ \; \gamma_{1k}} 
                          \frac{\tau}{\color{black}
                              \overline{
                          \rho_w(P_{n+1}) \phi(P_{n+1})}}
                          \Big(\overline{\rho_w(P_n) f_w(s_{\mathrm{in}})} \bar{q}_E
                          -\overline{\rho_w(P_n)} f_w(\bar{S}_{E}^{(k-1)})\underline{q}_E\Big)
                      },\label{eq:updateFL} \\
      \mathcal{H}_E^{(k)}(e) &= (1-\alpha_E^{(k-1)}(e))\, \mathcal{H}_E^{(k-1)}(e),
      \quad \forall e \subset \partial E.
          \end{align}
\item[Step~5.] Define a global stopping criterion\\
If  $\left(\max_{E\in\mathcal{E}_h}  \vert \mathcal{H}_E^{(k)}\vert
<\epsilon_1\right)$ or $\left(\max_{E\in\mathcal{E}_h}  \vert \mathcal{H}_E^{(k)}-\mathcal{H}_E^{(k-1)}\vert
<\epsilon_2\right)$ for $k\geq 2$\\
\hspace*{0.5cm} return $\bar{S}^{\mathrm{FL}}_{n+1}|_E = \bar{S}_E^{(k)}$.\\
Else\\
\hspace*{0.5cm} set $k\leftarrow k+1$ and go to Step~1.
\end{itemize}

\section{Numerical Results}%
\label{sec:numerical_results}
In this section, we study the effect of limiters by comparing numerical solutions
obtained without limiters (unlimited DG), and with flux and slope limiters (limited DG or DG+FL+SL).
We utilize the vertex-based slope limiter introduced in \citep{kuzmin2010vertex}.
For all problems addressed in this section, we assume the following parameters unless otherwise mentioned:
\begin{align*}
    &\rho_{w}^0 = 1000~\si{\kilogram\per\meter\cubed}, \quad
    \rho_{\ell}^0 = 850~\si{\kilogram\per\meter\cubed}, \quad
    \mu_{w} = 5\times10^{-4}~\si[inter-unit-product = \cdot]{\pascal\second}, \quad
    \mu_{\ell} = 2\times10^{-3}~\si[inter-unit-product = \cdot]{\pascal\second},  \quad\phi^0 = 0.15,
    \\
    &c_r = 9\times10^{-10}, \quad c_w ={10}^{-10}, \quad c_{\ell}=10^{-6} \quad
    s_{rw} = s_{r\ell} =  0.15, \quad s_0 = 0.15, \quad p_0 = 10^6~\mbox{Pa}, \quad \sigma = 100.
\end{align*}
\begin{figure}
    \subfigure[Pressure-driven flow problem \label{Fig:Sch_patch}]{
        \includegraphics[clip,width=0.48\linewidth,trim=0 0cm 0cm
        0]{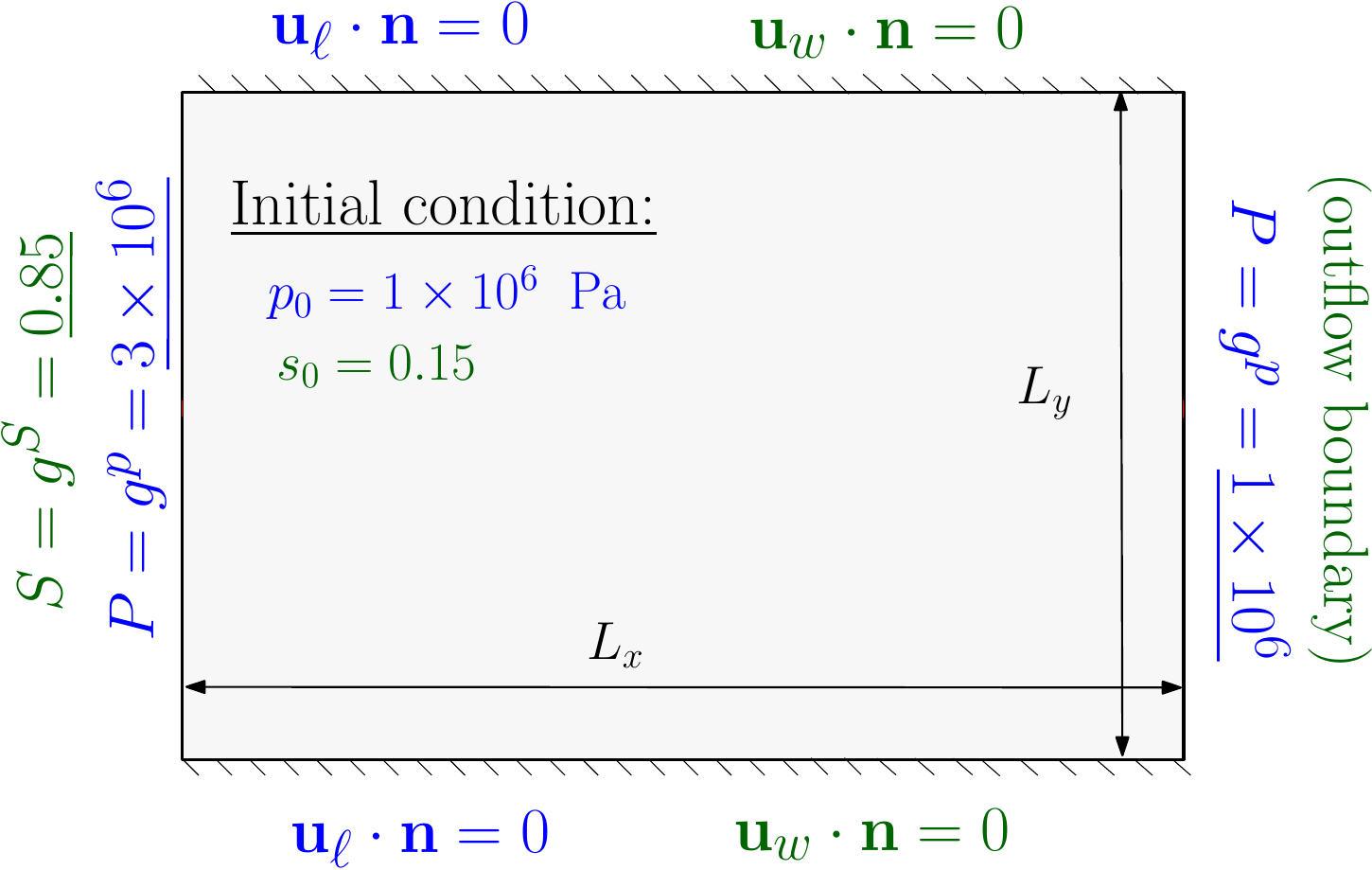}}
        \hspace{.6cm}
        \subfigure[Quarter five-spot flow problem \label{Fig:Sch_Q5}]{
        \includegraphics[clip,width=0.38\linewidth,trim=0 0cm 0cm 0]{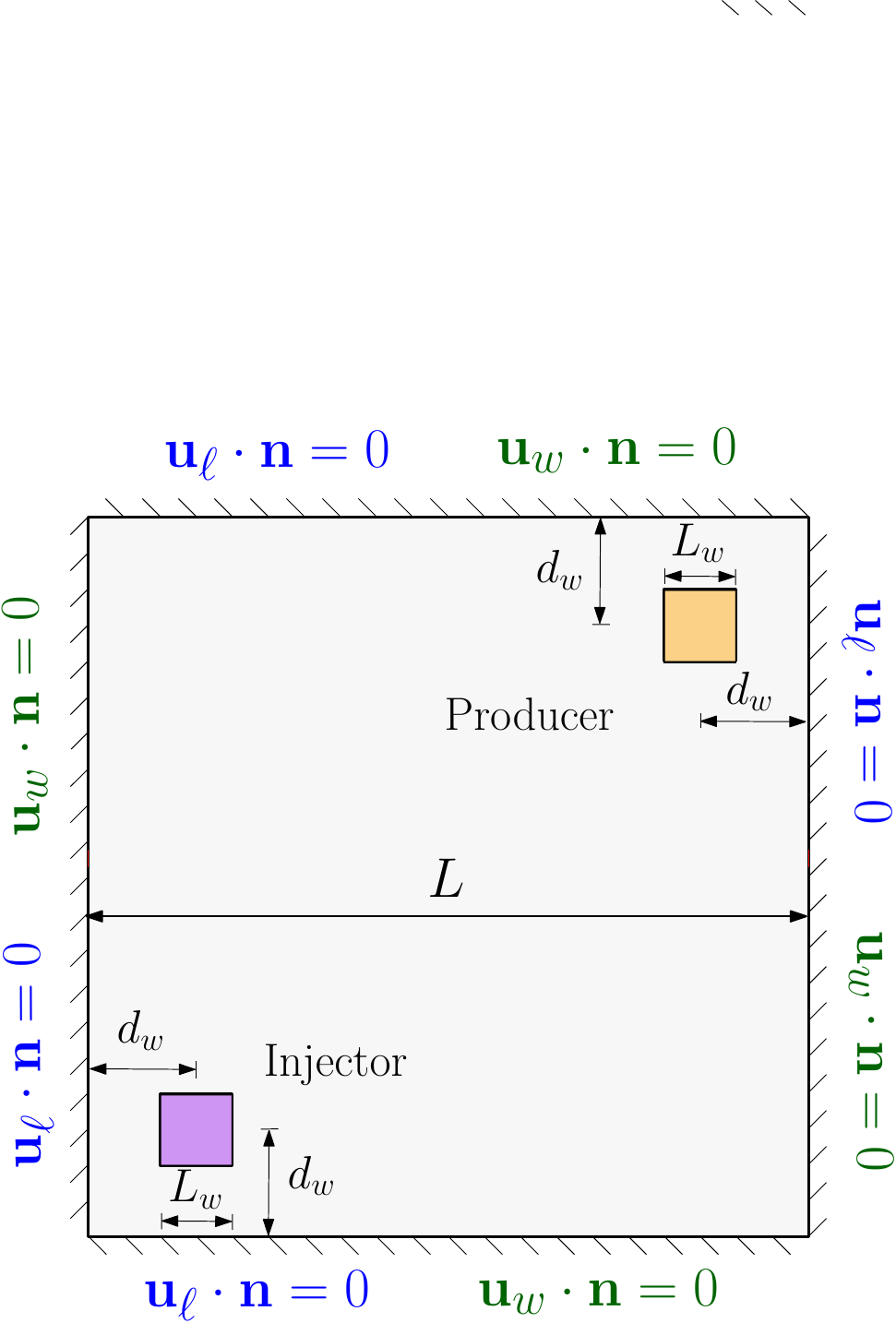}} 
        \caption{Pictorial description of computational domain and boundary conditions of pressure-driven flow problem and quarter five-spot problem.
        \label{Fig:main_schematic}}
\end{figure}

\subsection{Analytical problem and $h-$convergence study}%
\label{sub:convergence_study}
We first perform an $h-$convergence study on 2D structured triangular meshes of size $h$. The computational domain
is the unit square and the exact solutions are:
\begin{subequations}
    \begin{align}
        \label{Eqn:exact_S}
        &s(x,y,t) =  0.4+0.4xy+\cos{t+x}, \\
        \label{Eqn:exact_P}
        &p(x,y,t) =  2+x^2y-y^2+x^2\sin{y+t}-\frac{1}{3}\cos{t}+\frac{1}{3}\cos{t+1}-\frac{11}{6}.
    \end{align}
\end{subequations}
We replace the right-hand side of equations \eqref{Eqn:BoM_pres}--\eqref{Eqn:BoM_sat} by body forces 
obtained via the manufactured solutions. Dirichlet boundary conditions
are prescribed on $\partial\Omega$ on both saturation and pressure fields and the other parameters are taken 
as:
\begin{align*}
    &\phi=1,\quad
    K=1~\si{\meter\squared},\quad
    c_r=c_w=c_{\ell}=10^{-10},\quad\\
    &\rho_{w}^0=\rho_{\ell}^0=1~\si{\kilogram\per\meter\cubed},\quad
    \mathbf{g}=\mathbf{0},\quad
    \mu_{w}=\mu_{\ell}=1~\si[inter-unit-product =\cdot]{\pascal\second}.
\end{align*}
Table \ref{tab:MMF_Sat} and \ref{tab:MMF_pressure} show the errors in $L^2$ norm and rates evaluated at $T=1$ \si{\second} for saturation and pressure solutions. At each refinement level, the time step $\tau$ is set to $h^2$; and at every time instance $t$, the admissible bounds $s_*$ and $s^*$ are updated to the maximum and minimum of the exact saturation solution \eqref{Eqn:exact_S}. We compare the rates for three different cases of
unlimited DG, limited DG (DG+FL+SL) and also the case of DG with flux limiters only (DG+FL).
For both unknowns, DG and DG+FL yield expected optimal rate of $2$ in the $L^2$ norm whereas the application of slope limiters lead to suboptimal rates. We should highlight that the proposed flux limiter is rate-preserving and is independent of the slope limiter. Devising a rate preserving slope limiters still remains an open challenge.

\begin{table}
\centering
\caption{\label{tab:MMF_Sat}Errors in $L^2$ norm and convergence rates for saturation solutions at $T=1~\si{\second}$.}
\vspace{-0.25cm}
\resizebox{1\textwidth}{!}{%
\begin{tabular}{cc@{\hskip 0.5in}cc@{\hskip 0.5in}cc@{\hskip 0.5in}cc}
\Xhline{2\arrayrulewidth}
\\[-0.9em]
\multicolumn{1}{c}{} &
  \multicolumn{1}{c}{} &
  \multicolumn{2}{l@{\hskip 0.5in}}{DG} &
  \multicolumn{2}{l@{\hskip 0.5in}}{DG+FL} &
  \multicolumn{2}{l}{DG+FL+SL} \\ 
  \cmidrule(l{-0,1in}r{0.45 in}){3-4} \cmidrule(l{-0.1in}r{0.45in}){5-6} \cmidrule(l{-0.1in}r){7-8}
  \multicolumn{1}{c}{\multirow{-2}{*}{$h$ (\si{\meter})}} & \multicolumn{1}{c@{\hskip 0.5in}}{\multirow{-2}{*}{dofs}} & $\| S_n-s(T)|\vert_{L^2(\Omega)}$ & rate & $\| S_n-s(T)|\vert_{L^2(\Omega)}$ & rate & $\| S_n-s(T)|\vert_{L^2(\Omega)}$ & rate 
\\[-0.9em]
  \\ \hline
\\[-0.9em]
1/2 &
  48 &
  $5.08\times10^{-2}$ &
  \cellcolor[HTML]{EFEFEF} $-$ &
  $4.92\times10^{-2}$ &
  \cellcolor[HTML]{EFEFEF}$-$ &
  $3.39\times10^{-2}$ &
  \cellcolor[HTML]{EFEFEF}$-$ \\
1/4 &
  192 &
  $1.65\times10^{-2}$ &
  \cellcolor[HTML]{EFEFEF}1.62 &
  $1.67\times10^{-2}$ &
  \cellcolor[HTML]{EFEFEF}1.56 &
  $1.35\times10^{-2}$ &
  \cellcolor[HTML]{EFEFEF}1.33 \\
1/8 &
  768 &
  $4.03\times10^{-3}$ &
  \cellcolor[HTML]{EFEFEF}2.03 &
  $4.05\times10^{-3}$ &
  \cellcolor[HTML]{EFEFEF}2.04 &
  $6.54\times10^{-3}$ &
  \cellcolor[HTML]{EFEFEF}1.04 \\
1/16 &
  3072 &
  $1.02\times10^{-3}$ &
  \cellcolor[HTML]{EFEFEF}1.98 &
  $1.02\times10^{-3}$ &
  \cellcolor[HTML]{EFEFEF}1.99 &
  $3.28\times10^{-3}$ &
  \cellcolor[HTML]{EFEFEF}1.00 \\
1/32 &
  12288 &
  $2.61\times10^{-4}$ &
  \cellcolor[HTML]{EFEFEF}1.96 &
  $2.61\times10^{-4}$ &
  \cellcolor[HTML]{EFEFEF}1.96 &
  $1.39\times10^{-3}$ &
  \cellcolor[HTML]{EFEFEF}1.24 \\
1/64 &
  49152 &
  $6.92\times10^{-5}$ &
  \cellcolor[HTML]{EFEFEF}1.92 &
  $6.92\times10^{-5}$ &
  \cellcolor[HTML]{EFEFEF}1.92 &
  $5.56\times10^{-4}$ &
  \cellcolor[HTML]{EFEFEF}1.32 \\ 
\\[-0.9em]
\Xhline{2\arrayrulewidth}
\end{tabular}%
}
\end{table}

\begin{table}
\centering
\caption{\label{tab:MMF_pressure}Errors in $L^2$ norm and convergence rates for pressure solutions at $T=1~\si{\second}$.}
\vspace{-0.25cm}
\resizebox{1\textwidth}{!}{%
\begin{tabular}{cc@{\hskip 0.5in}cc@{\hskip 0.5in}cc@{\hskip 0.5in}cc}
\Xhline{2\arrayrulewidth}
\\[-0.9em]
\multicolumn{1}{c}{} &
  \multicolumn{1}{c}{} &
  \multicolumn{2}{l@{\hskip 0.5in}}{DG} &
  \multicolumn{2}{l@{\hskip 0.5in}}{DG+FL} &
  \multicolumn{2}{l}{DG+FL+SL} \\ 
  \cmidrule(l{-0,1in}r{0.45 in}){3-4} \cmidrule(l{-0.1in}r{0.45in}){5-6} \cmidrule(l{-0.1in}r){7-8}
  \multicolumn{1}{c}{\multirow{-2}{*}{$h$ (\si{\meter})}} & \multicolumn{1}{c@{\hskip 0.5in}}{\multirow{-2}{*}{dofs}} & $\| P_n-p(T)|\vert_{L^2(\Omega)}$ & rate & $\| P_n-p(T)|\vert_{L^2(\Omega)}$ & rate & $\| P_n-p(T)|\vert_{L^2(\Omega)}$ & rate 
\\[-0.9em]
  \\ \hline
\\[-0.9em]
1/2 &
  48 &
  $2.68\times10^{-3}$ &
  \cellcolor[HTML]{EFEFEF}$-$ &
  $2.68\times10^{-3}$ &
  \cellcolor[HTML]{EFEFEF}$-$ &
  $1.74\times10^{-3}$ &
  \cellcolor[HTML]{EFEFEF}$-$ \\
1/4 &
  192 &
  $2.25\times10^{-3}$ &
  \cellcolor[HTML]{EFEFEF}0.25 &
  $2.25\times10^{-3}$ &
  \cellcolor[HTML]{EFEFEF}0.25 &
  $1.45\times10^{-3}$ &
  \cellcolor[HTML]{EFEFEF}0.26 \\
1/8 &
  768 &
  $6.74\times10^{-4}$ &
  \cellcolor[HTML]{EFEFEF}1.74 &
  $6.77\times10^{-4}$ &
  \cellcolor[HTML]{EFEFEF}1.73 &
  $5.14\times10^{-4}$ &
  \cellcolor[HTML]{EFEFEF}1.50 \\
1/16 &
  3072 &
  $1.80\times10^{-4}$ &
  \cellcolor[HTML]{EFEFEF}1.91 &
  $1.80\times10^{-4}$ &
  \cellcolor[HTML]{EFEFEF}1.91 &
  $1.32\times10^{-4}$ &
  \cellcolor[HTML]{EFEFEF}1.96 \\
1/32 &
  12288 &
  $4.55\times10^{-5}$ &
  \cellcolor[HTML]{EFEFEF}1.98 &
  $4.55\times10^{-5}$ &
  \cellcolor[HTML]{EFEFEF}1.98 &
  $3.20\times10^{-5}$ &
  \cellcolor[HTML]{EFEFEF}2.04 \\
1/64 &
  49152 &
  $1.14\times10^{-5}$ &
  \cellcolor[HTML]{EFEFEF}1.99 &
  $1.14\times10^{-5}$ &
  \cellcolor[HTML]{EFEFEF}2.00 &
  $8.08\times10^{-6}$ &
  \cellcolor[HTML]{EFEFEF}1.99 \\[-0.9em] \\ 
  \Xhline{2\arrayrulewidth}
\end{tabular}%
}
\end{table}

\subsection{Pressure-driven flow}%
\label{sub:pressure_driven_flow}
In this section, we perform various pressure-driven flow problems with homogeneous and heterogeneous permeabilities to study the
efficacy and robustness of limiters on capturing accurate and bound-preserving solutions.
For all problems, we take a rectangular computational domain $\Omega = [0,L_x]\times[0,L_y]$~\si{\meter\squared}. A wetting phase is injected along the left boundary and displaces  the non-wetting phase  out of the domain
through the right boundary. As depicted in Figure~\ref{Fig:Sch_patch}, Dirichlet boundary conditions are set to: 
$g^p=3\times10^{6}$~\si{\pascal} and  $g^s=0.85$ on $\{0\}\times [0,L_y]$~\si{\meter}; and  $g^p=10^6$~\si{\pascal} 
on $\{L_x\}\times [0,L_y]$~\si{\meter}. Outflow boundary condition is prescribed on the right boundary for saturation and the top/bottom boundaries are set as no-flow ($j^s=j^p=0$). We note that due to the residual saturations, the exact saturation satisfies the
maximum principle:
\[
0.15 \leq s \leq 0.85.
\]
We will highlight below the behavior of the phase saturation regarding these physical bounds.
\subsubsection{Example~1: Homogeneous domain}%
\begin{figure}
    \subfigure[DG (without limiter) at $t=5.1$ days \label{Fig:PD_homogen_sat_1a}]{
        \includegraphics[clip,scale=0.11,trim=0 0cm 0cm 0]{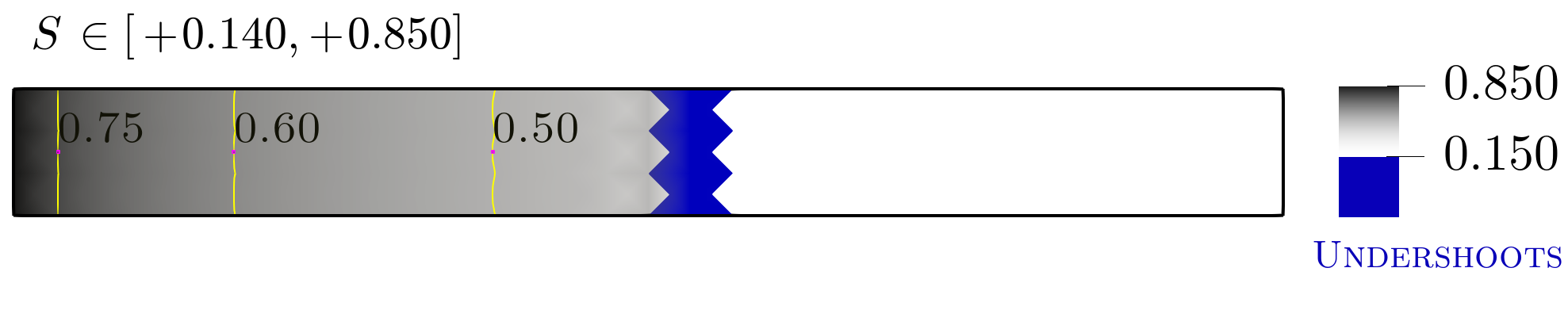}}
        \hspace{.25cm}
    \subfigure[DG (without limiter) at $t=10$ days \label{Fig:PD_homogen_sat_1b}]{
        \includegraphics[clip,scale=0.107,trim=0 0cm 0cm 0]{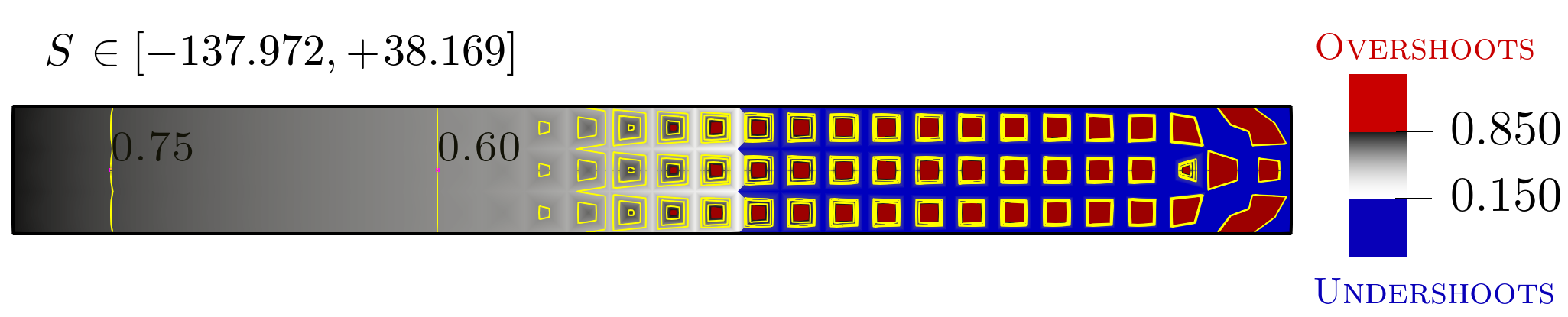}} \\
    \subfigure[DG + SL at $t=5.1$ days \label{Fig:PD_homogen_sat_1c}]{
        \includegraphics[clip,scale=0.108,trim=0 0cm 0cm 0]{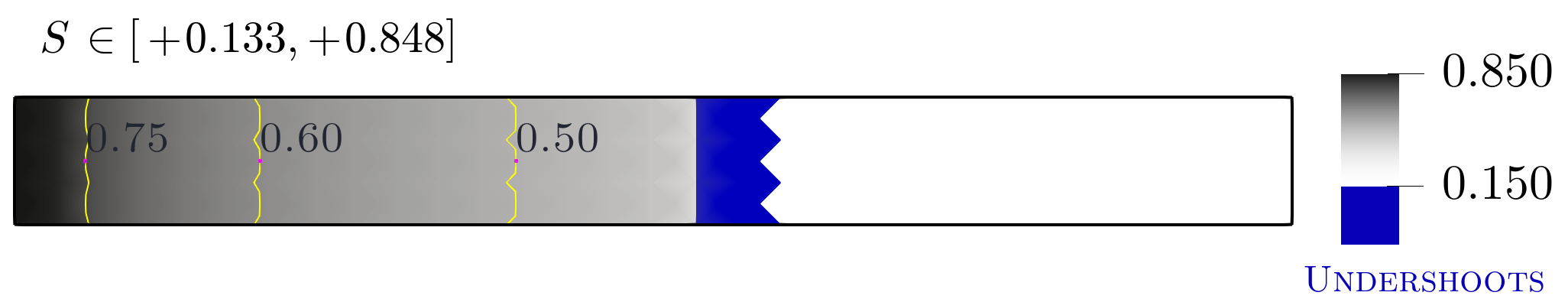}}
        \hspace{.25cm}
    \subfigure[DG + SL at $t=10$ days \label{Fig:PD_homogen_sat_1d}]{
        \includegraphics[clip,scale=0.108,trim=0 0cm 0cm 0]{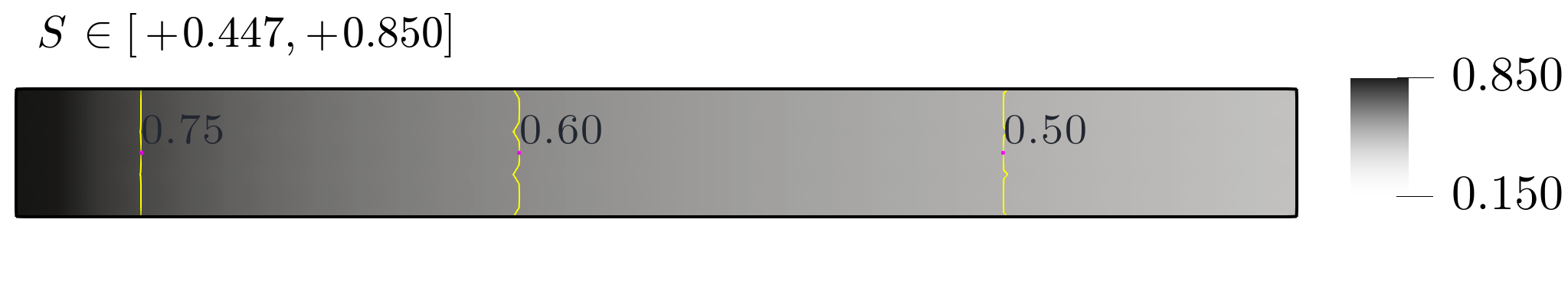}} \\
    \subfigure[DG + FL + SL at $t=5.1$ days \label{Fig:PD_homogen_sat_1e}]{
        \includegraphics[clip,scale=0.108,trim=0 0cm 0cm 0]{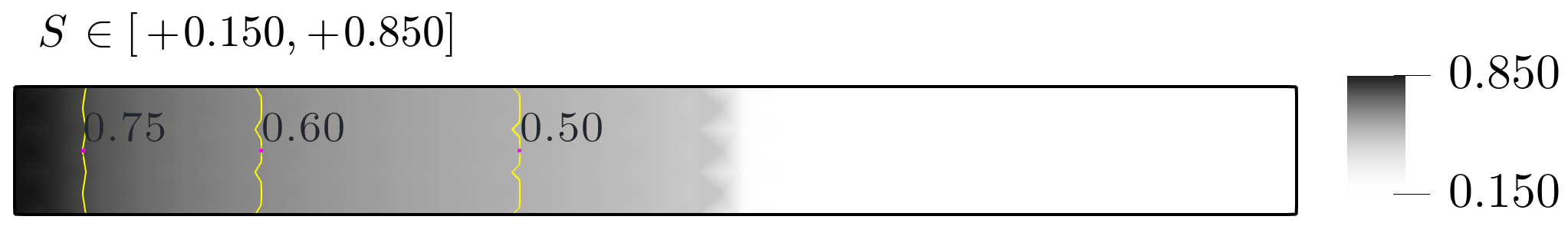}}
        \hspace{.25cm}
    \subfigure[DG + FL + SL at $t=10$ days \label{Fig:PD_homogen_sat_1f}]{
        \includegraphics[clip,scale=0.108,trim=0 0cm 0cm 0]{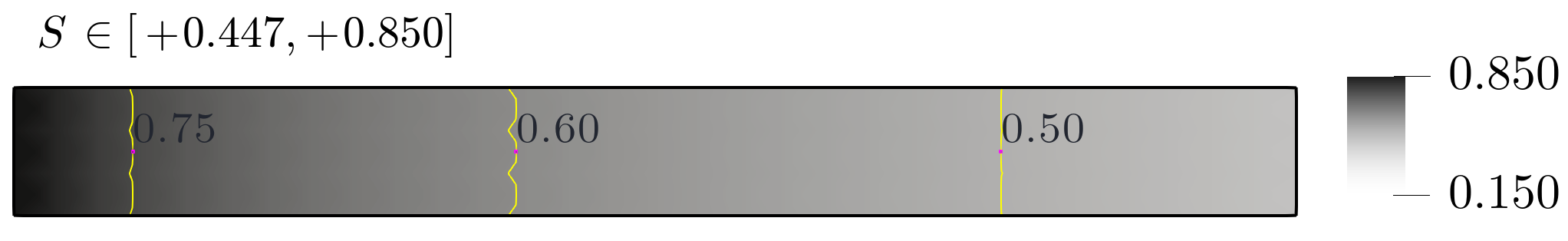}} \\
        \caption{\textsf{Pressure-driven flow in homogeneous domain:~}The evolution of saturation solutions obtained with
            DG (top), DG+SL (middle), and DG+FL+SL (bottom) at two different time steps. 
            The color mapping for physical range of saturation is presented
            in grayscale, whereas the values below and above bounds are colored blue and red, respectively.
            Results suggest that DG+FL+SL, unlike the two other schemes, provides maximum-principle-satisfying solutions at
            all time steps.
        \label{Fig:PD_homogen_sat_main}}
\end{figure}

\label{ssub:homogeneous_doamin}
We consider a domain of length $L_x=100$~\si{\meter} and $L_y=30$~\si{\meter} with constant permeability of $K=10^{-12}$~\si{\meter\squared} partitioned into structured crossed triangular meshes of size $h=10$~\si{\meter}. The time step is chosen as $0.05$~days and the
final time is $T=10$~days. Gravity is neglected. Figure~\ref{Fig:PD_homogen_sat_main} compares
the saturation profiles obtained from unlimited DG, DG with only vertex-based slope limiter \citep{kuzmin2010vertex} (i.e., DG+SL), and the proposed limited DG (i.e., DG+FL+SL) at two different time steps. It is seen that the saturation front, 
under all three approximations, propagates  with the same speed. However,
limited DG, unlike its unlimited counterparts, produces a numerical saturation that 
remains physically bounded and neither undershoots (blue-colored elements) nor overshoots (red-colored elements) are detected throughout the simulation.  Table~\ref{Tab:Homogen_compare} shows the minimum and maximum values of the saturation over all time steps. While the slope limiter removes the overshoot for this simulation, there is still significant undershoot. 
The percentage of these overshoot and undershoot with respect to the physical range $[0.15, 0.85]$ are also displayed.

\begin{table}
\centering
\caption{This table shows the efficacy of the limiting schemes when applied to the pressure-driven flow problem with homogeneous domain. \label{Tab:Homogen_compare}}
\vspace{-0.25cm}
    \begin{tabular}{l@{\hskip 0.5in}cc@{\hskip 0.5in}cc@{\hskip 0.5in}cc}
\Xhline{2\arrayrulewidth}\\[-0.9em]
               & \multicolumn{2}{l}{DG}                                     & \multicolumn{2}{l}{DG+SL} & \multicolumn{2}{l}{DG+FL+SL} \\ \cmidrule(l{-0.1in}r{0.42 in}){2-3} \cmidrule(l{-0.1in}r{0.42 in}){4-5}\cmidrule(l{-0,1in}r{0.0 in}){6-7}
               & value                        & \%                  & value     & \%    & value     & \%  \\[-0.9em]  \\ \hline\\[-0.9em]
minimum saturation  & -137.97 &  19731 & 0.033 & 16.7 & 0.15 & 0 \\
maximum saturation  & 38.16 &   5330 & 0.85      & 0             & 0.85      & 0          \\[-0.9em] \\
\Xhline{2\arrayrulewidth}
\end{tabular}%
\end{table}
\subsubsection{Example~2: Heterogeneous domain}%
\begin{figure}[htpb]
    \centering
    \includegraphics[width=0.8\linewidth]{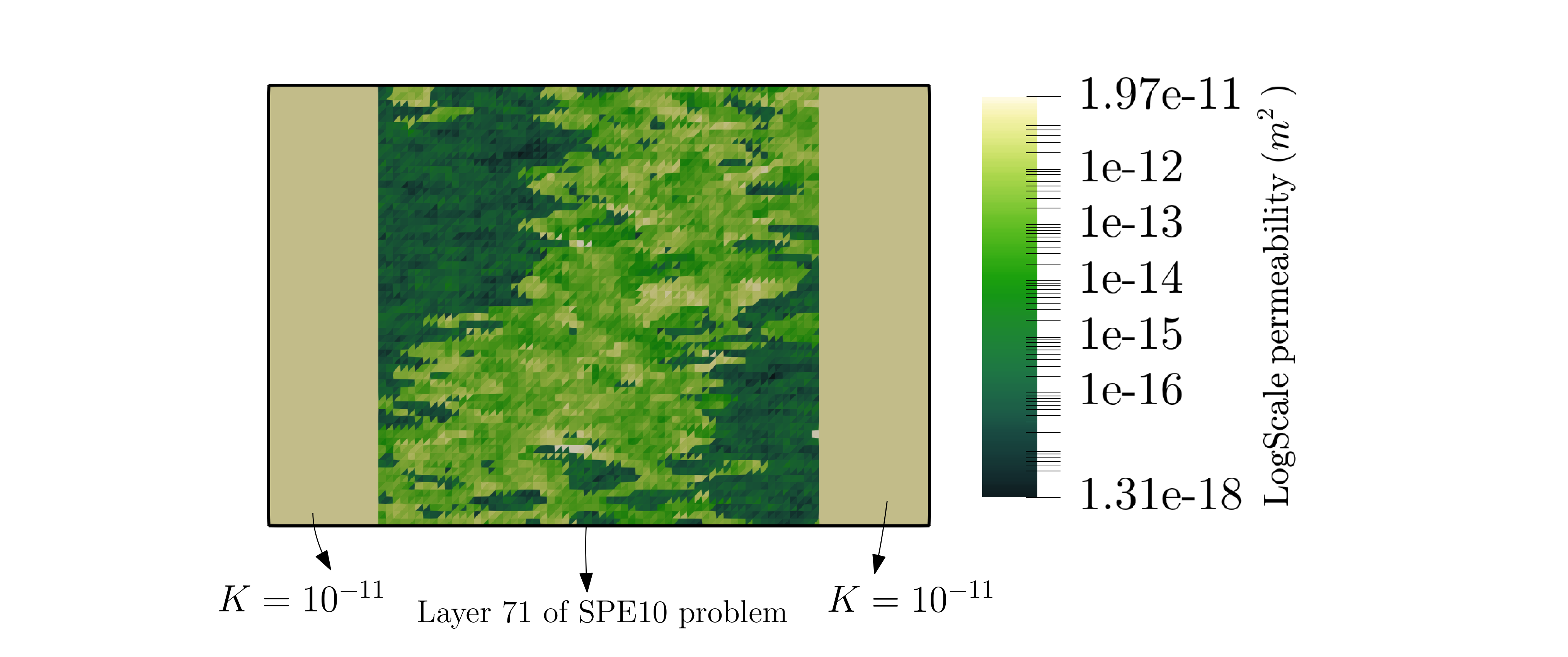}
    \caption{\textsf{Pressure-driven flow in heterogeneous domain:~}The permeability field in the center is adopted from layer 71 of SPE10 benchmark problem. Values are shown in logarithmic scale \label{Fig:pressure_driven_het}}
\end{figure}

\begin{figure}
    \subfigure[Unlimited DG at $t=25$ days \label{Fig:PD_heterogen_sat_a}]{
        \includegraphics[clip,scale=0.11,trim=0 0cm 0cm 1]{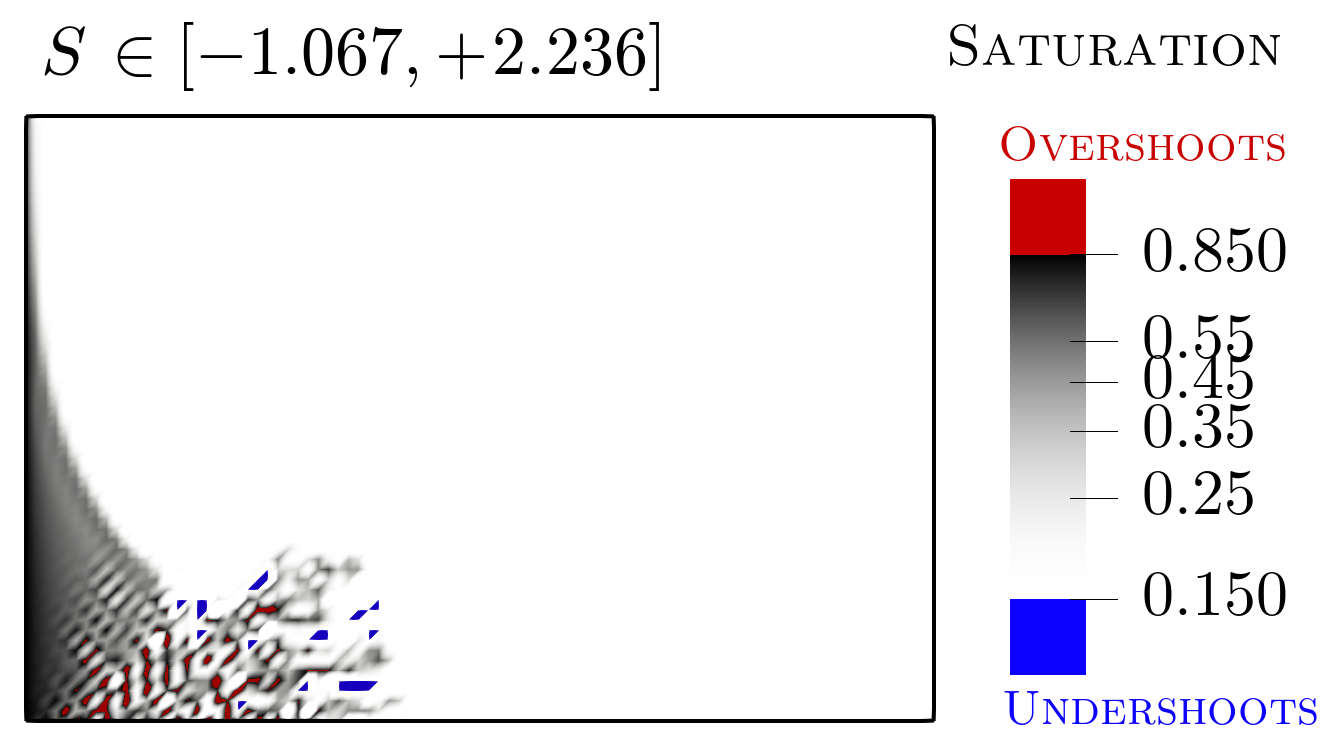}}
        \hspace{.0cm}
        \subfigure[Unlimited DG at $t=50$ days \label{Fig:PD_heterogen_sat_b}]{
        \includegraphics[clip,scale=0.11,trim=0 0cm 0cm 0]{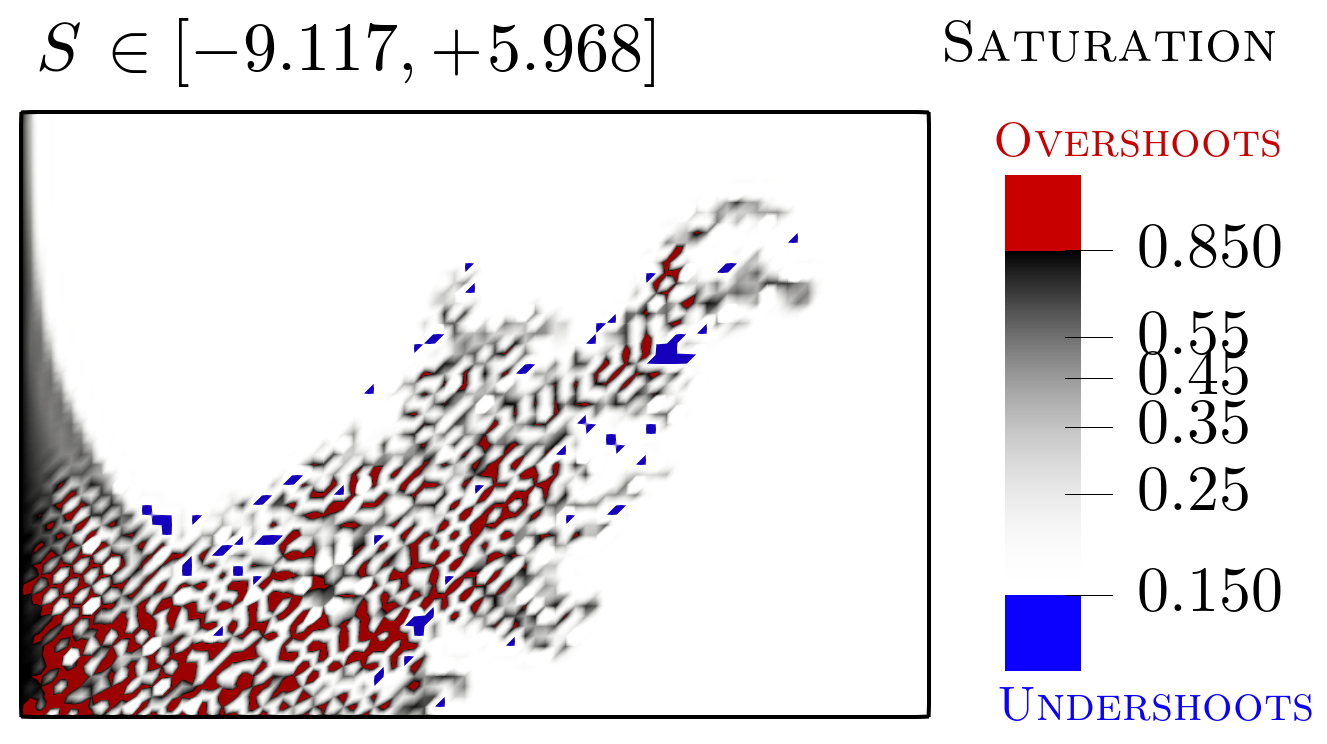}} 
        \hspace{.0cm}
        \subfigure[Unlimited DG at $t=68$ days \label{Fig:PD_heterogen_sat_c}]{
        \includegraphics[clip,scale=0.11,trim=0 0cm 0cm 0]{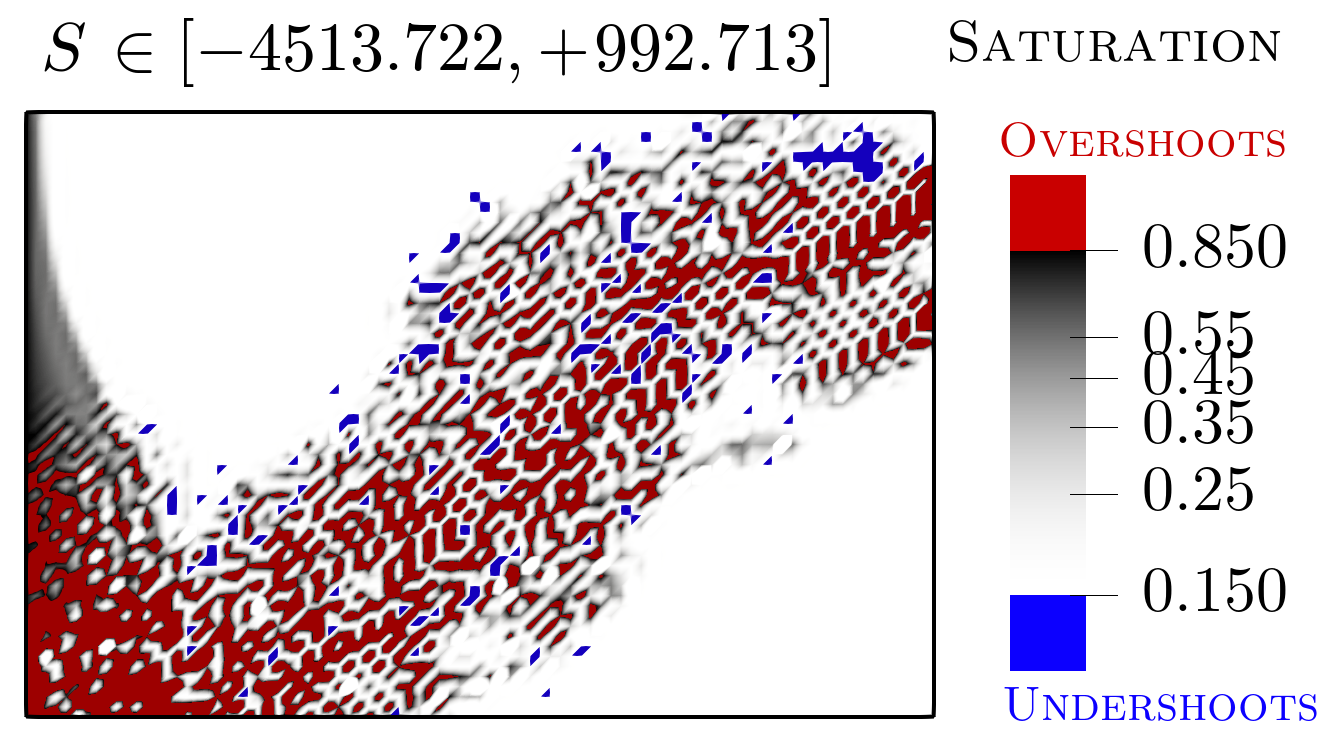}} \\
    \subfigure[Limited DG at $t=25$ days \label{Fig:PD_heterogen_sat_d}]{
        \includegraphics[clip,scale=0.11,trim=0 0cm 0cm 1]{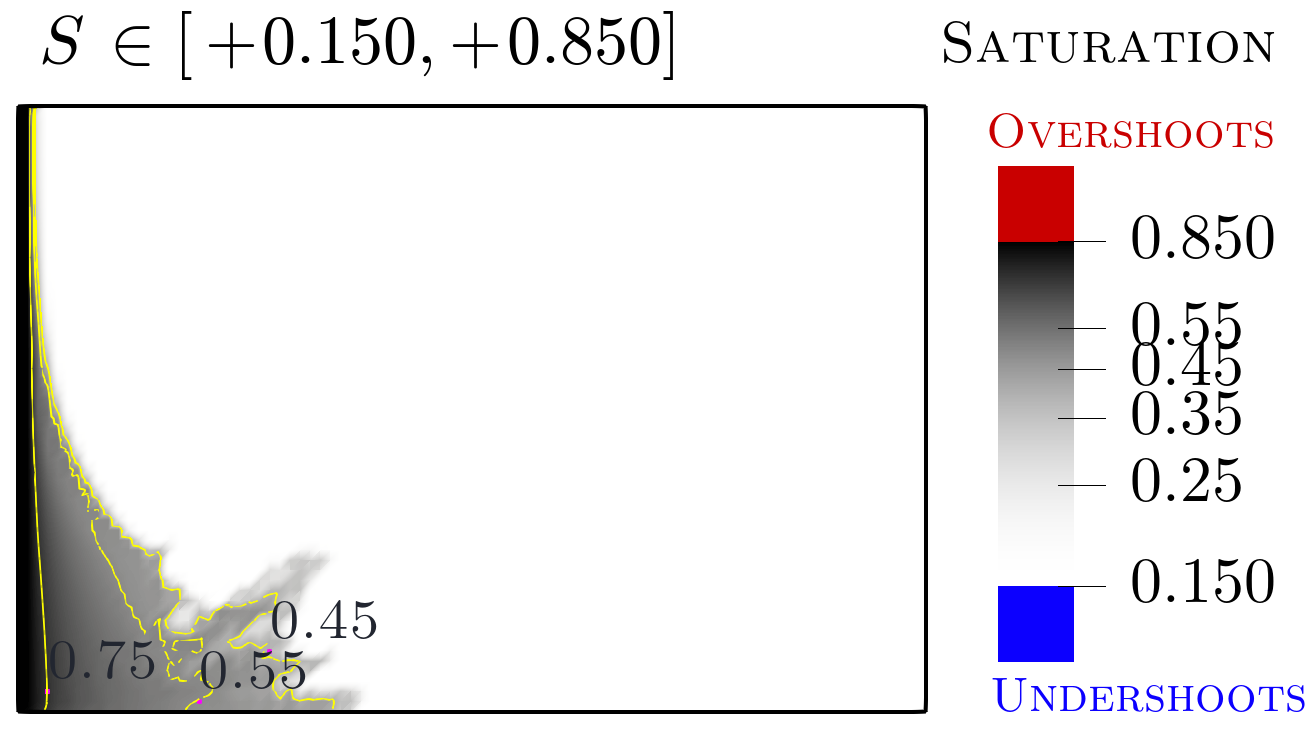}}
        \hspace{.0cm}
        \subfigure[Limited DG at $t=50$ days \label{Fig:PD_heterogen_sat_e}]{
        \includegraphics[clip,scale=0.11,trim=0 0cm 0cm 0]{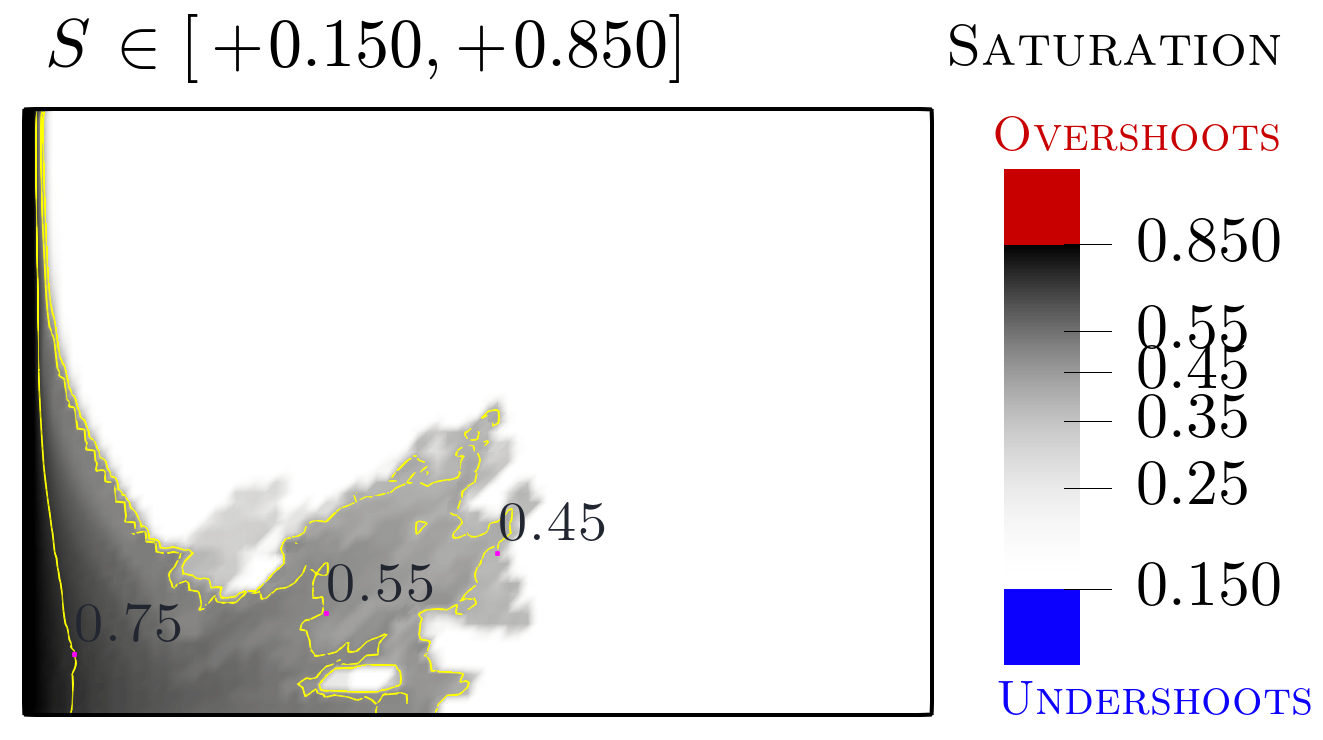}} 
        \hspace{.0cm}
        \subfigure[Limited DG at $t=68$ days \label{Fig:PD_heterogen_sat_f}]{
        \includegraphics[clip,scale=0.11,trim=0 0cm 0cm 0]{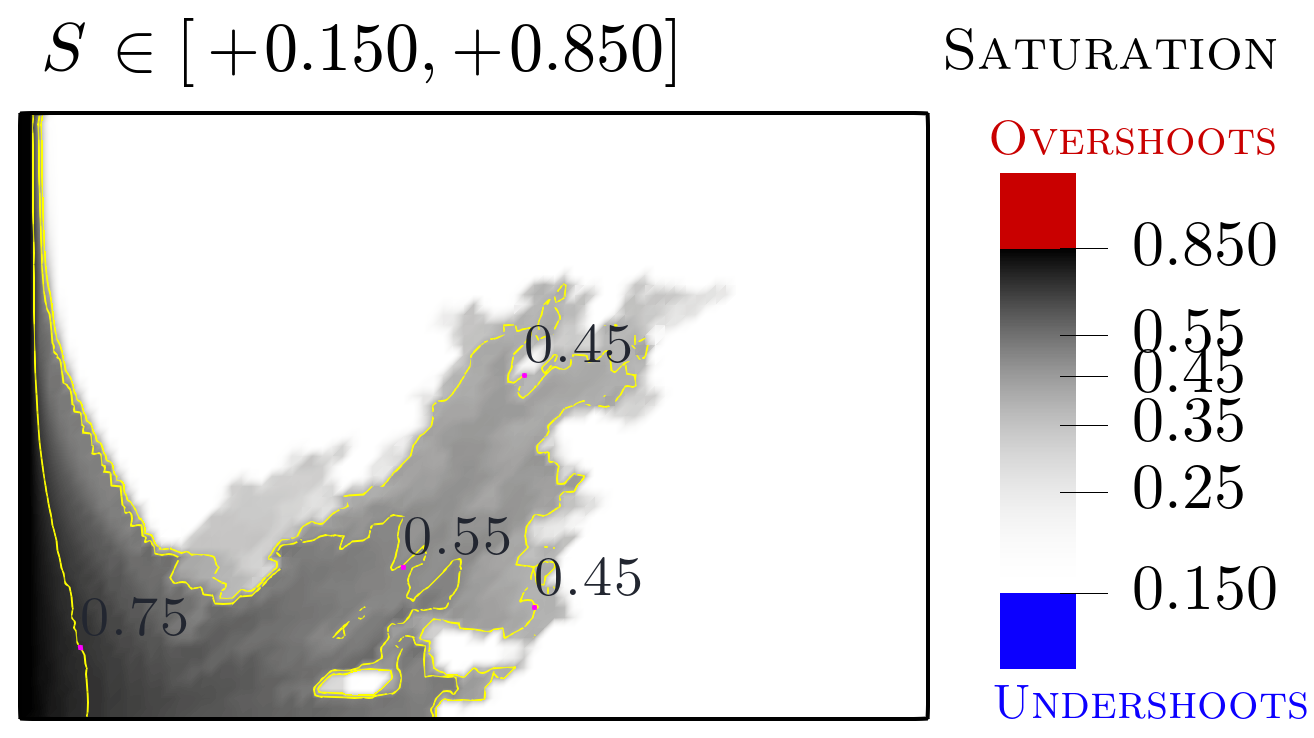}} \\
        \caption{\textsf{Pressure-driven flow in heterogeneous domain:~}The evolution of saturation profile using DG scheme without limiter (top row) and with the proposed limiters (bottom row). The color mapping for physical range of saturation is presented
        in grayscale, whereas the values below and above bounds are colored blue and red, respectively. DG approximation 
    give rise to noticeable violations but limited DG is capable of producing maximum-principle-satisfying results.
        \label{Fig:PD_heterogen_sat_main}}
\end{figure}

\begin{figure}
    \subfigure[Unlimited DG at $t=68$ days \label{Fig:PD_heterogen_vel_a}]{
        \includegraphics[clip,scale=0.12,trim=0 0cm 0cm 0]{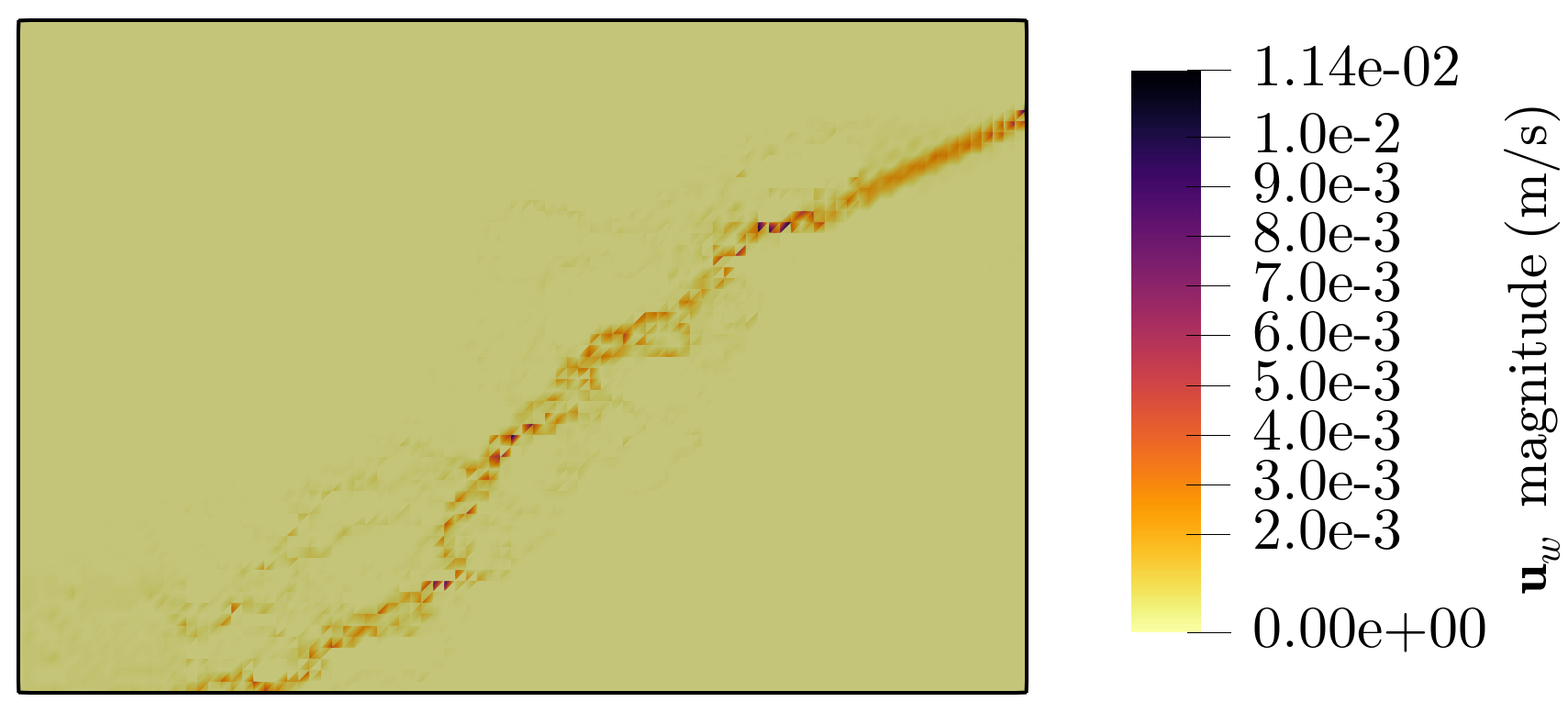}}
        \hspace{.25cm}
    \subfigure[Limited DG at $68$ days \label{Fig:PD_heterogen_vel_b}]{
        \includegraphics[clip,scale=0.12,trim=0 0cm 0cm 0]{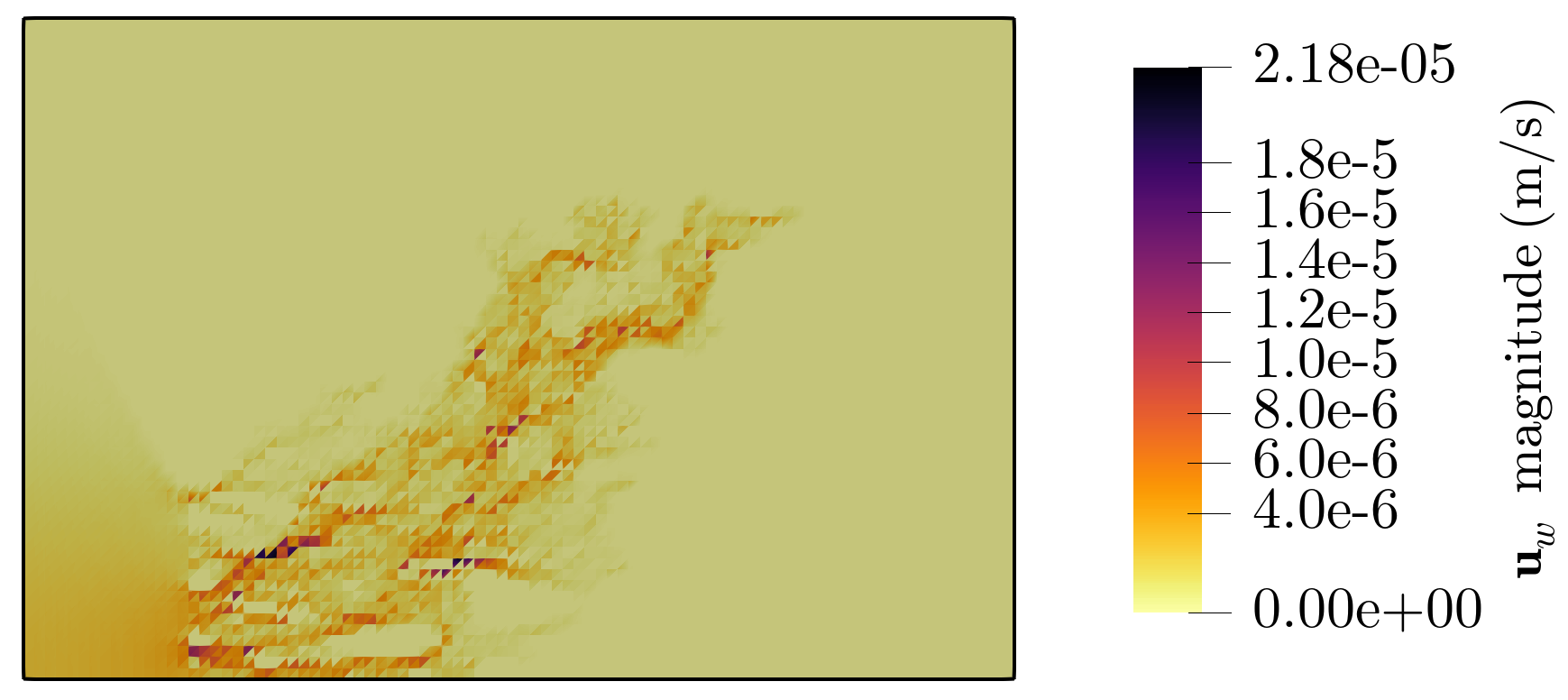}} \\
        \caption{\textsf{Pressure-driven flow in heterogeneous domain:} The magnitude of velocity field
            at final time $t=68$ days using DG scheme without limiter (left) and with limiters (right). 
        \label{Fig:PD_heterogen_vel_main}}
\end{figure}

\begin{figure}
    \subfigure[Unlimited DG at $t=68$ days \label{Fig:PD_heterogen_pres_a}]{
        \includegraphics[clip,scale=0.15,trim=0 0cm 0cm 0]{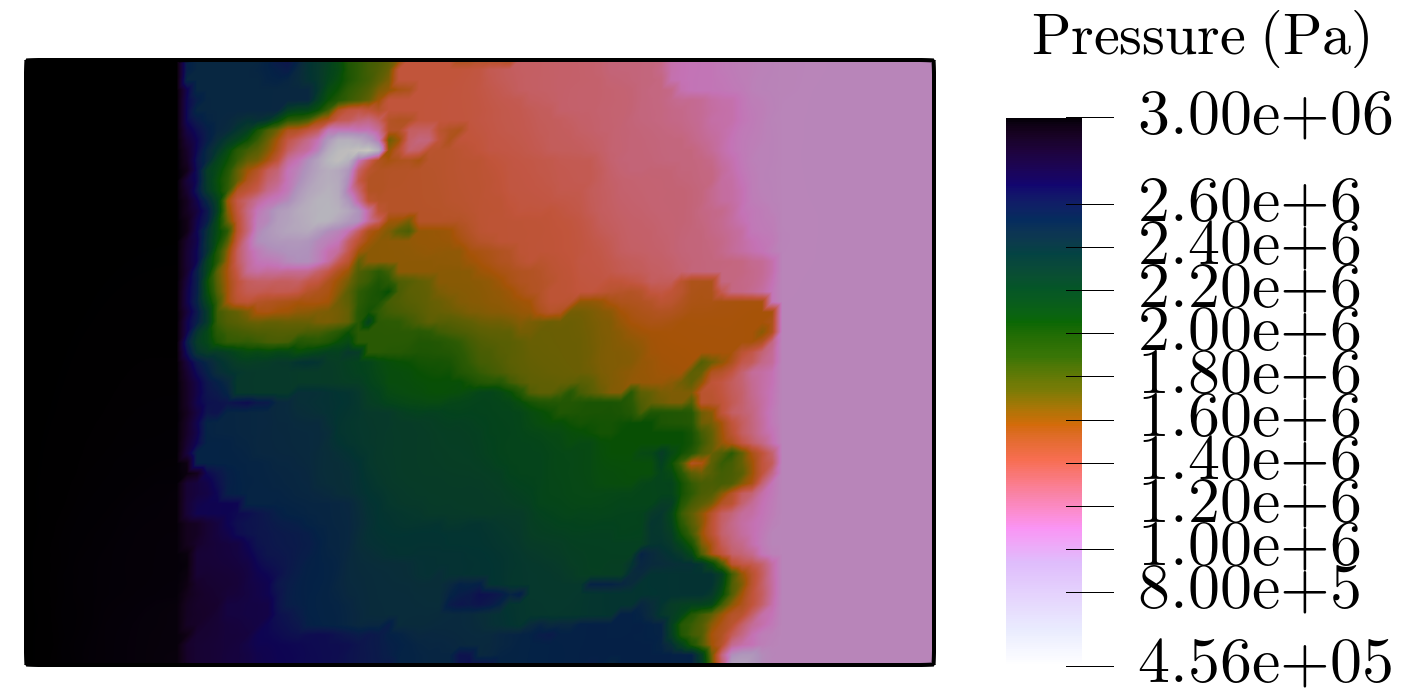}}
        \hspace{.25cm}
    \subfigure[Limited DG at $68$ days \label{Fig:PD_heterogen_pres_b}]{
        \includegraphics[clip,scale=0.15,trim=0 0cm 0cm 0]{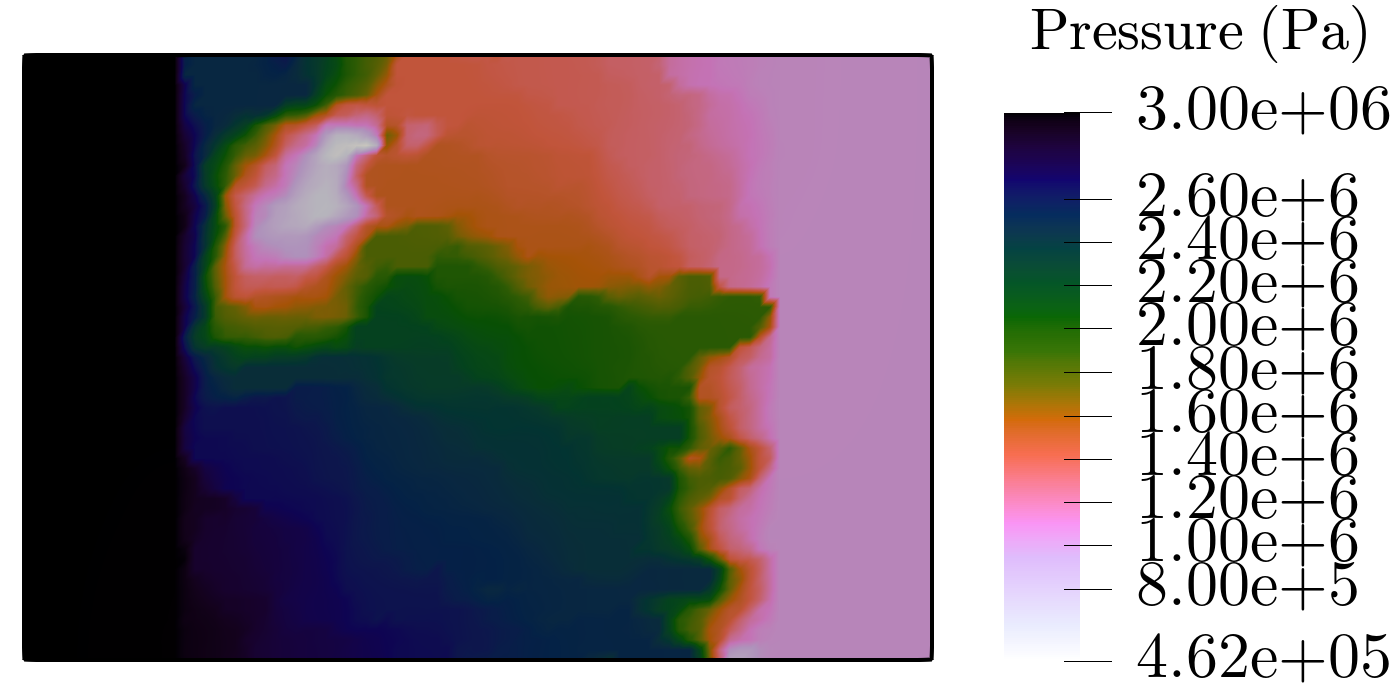}} \\
        \caption{\textsf{Pressure-driven flow in heterogeneous porous media:} The pressure contours at
time $t=68$ days using unlimited DG (left) and limited DG (right). 
        \label{Fig:PD_heterogen_pres_main}}
\end{figure}

\label{ssub:heterogeneous_domain}
We repeat the experiment in Section \ref{ssub:homogeneous_doamin} with an heterogeneous domain of size $L_x=150$ \si{\meter} and 
$L_y=100$ \si{\meter}. The permeability field is composed of a highly discontinuous central block sandwiched by two 
buffer zones of $K=10^{-11}$ \si{\meter\squared} (see Figure~\ref{Fig:pressure_driven_het} for a description of the 
permeability field). The data for the central block is taken from the horizontal permeability slice  number $71$ of the 
SPE10 benchmark model \citep{SPE10reference} that is scaled to a $100\times100$ \si{\meter\squared} grid. 
The domain is discretized with structured triangular mesh of size $h=5/3$ \si{\meter}, time step is $\tau=1/12$ days, and 
total simulation time is $T=68$ days. Figure~\ref{Fig:PD_heterogen_sat_main} shows the saturation profile under limited and 
unlimited DG at three different time instances $t=25$, $50$, and $68$ days. The wetting phase floods the domain from
the left buffer zone toward the right buffer zone while avoiding the low permeable regions. As expected, even 
for highly heterogeneous domains, 
the proposed 
limited DG completely eliminates the violation of maximum principle that appeared as overshoots and undershoots in 
unlimited DG approximations. 
Figure~\ref{Fig:PD_heterogen_vel_main} depicts the magnitude of the wetting phase velocity, $\mathbf{u}_w$, computed at $t=68$ days. The velocity
is computed at time $t_n$ in each mesh element by
\begin{equation}\label{eq:defuw}
\mathbf{u}_w^n = -\lambda_w(S_n) K (\nabla P_n -\rho_\ell(P_n) \mathbf{g}).
\end{equation}
Velocity obtained under DG with no limiters  does not accurately follow the path of saturation propagation and exhibits overestimation
of the magnitude of the velocity. On the other hand, the limited DG scheme eliminates these shortcomings
            and results in 
            distinguishable flow paths that match those of the saturation contours.
The pressure contours are shown in Figure~\ref{Fig:PD_heterogen_pres_main} for both unlimited DG and limited DG. Both methods
produce the same pressure range, but there are visible differences in the pressure field in the heterogeneous region.

\subsubsection{Example~3: Non-homogeneous domain with gravitational force}%
\begin{figure}
    \subfigure[No gravity; $t=30$ days \label{Fig:PD_gravity_sat_a}]{
        \includegraphics[clip,scale=0.16,trim=0 0cm 0cm 0]{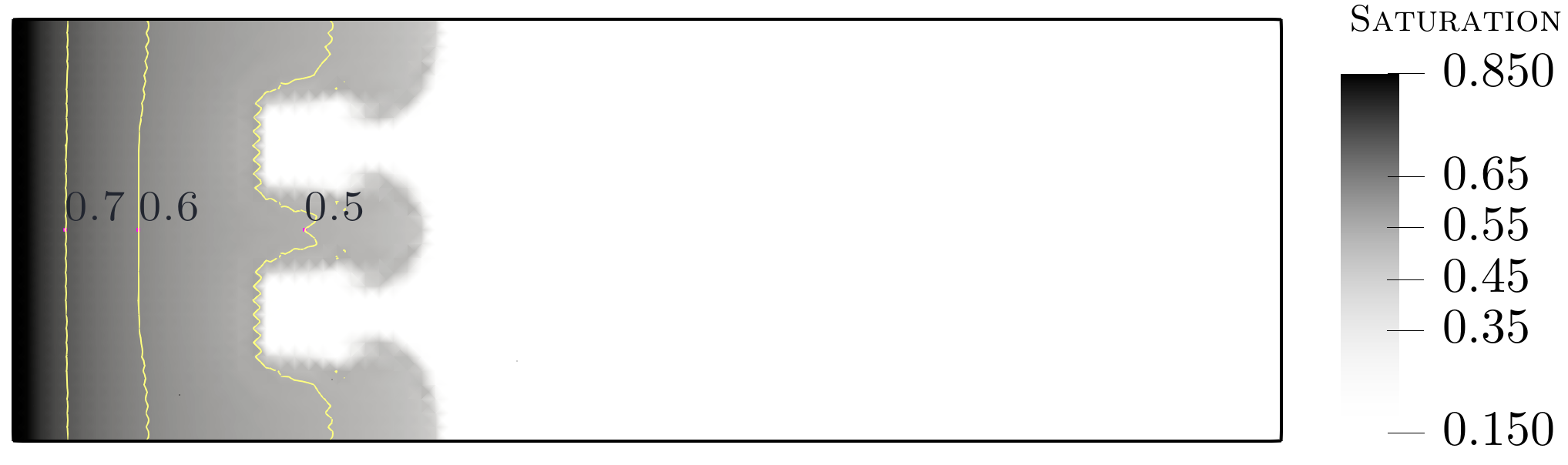}}
        \hspace{.25cm}
    \subfigure[With gravity; $t=30$ days \label{Fig:PD_gravity_sat_b}]{
        \includegraphics[clip,scale=0.172,trim=0 0cm 0cm 0]{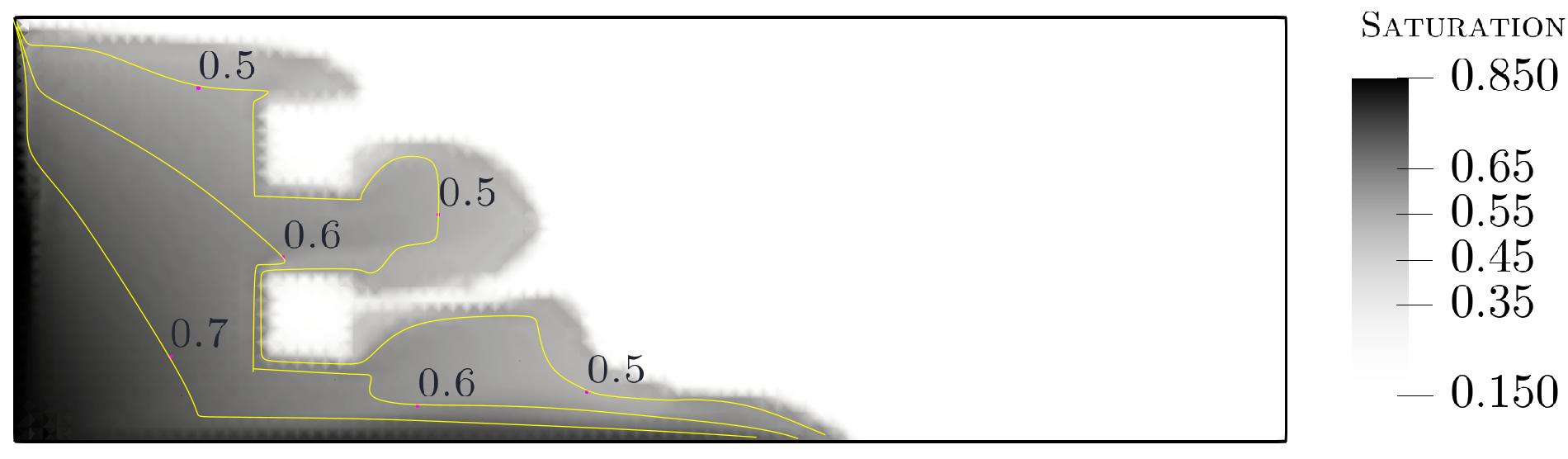}} \\
        \caption{\textsf{Pressure-driven flow in non-homogeneous domain with gravitational force:} This figure shows
            the wetting phase saturation contours at $t=30$ days for two different cases where the gravitational force
            is absent (top) and is present (bottom). Both flux and slope limiter are utilized in generating results and
            hence no violation of maximum principle is observed.
        \label{Fig:PD_Gravity_sat_main}}
\end{figure}

\begin{figure}
    \subfigure[Without gravitational forces; $t=30$ days \label{Fig:PD_gravity_pres_a}]{
        \includegraphics[clip,scale=0.18,trim=0 0cm 0cm 0]{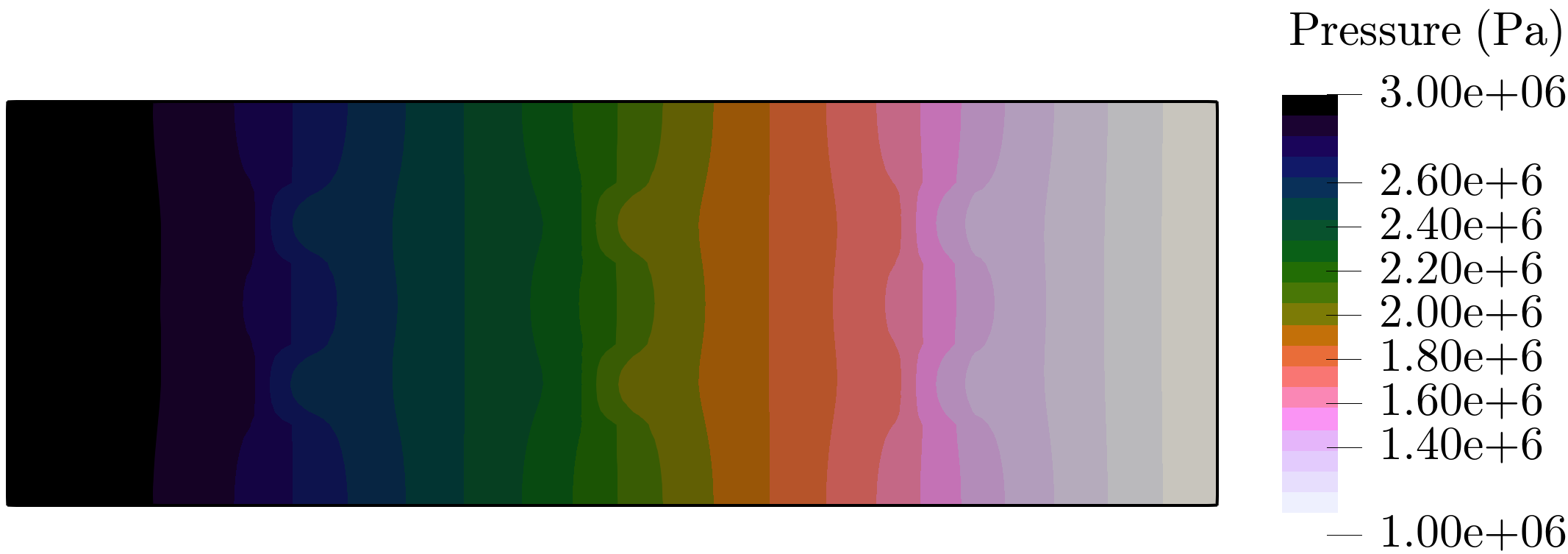}}
        \hspace{.25cm}
        \subfigure[With gravity body force ($\mathrm{Gr}=0.4$); $t=30$ days \label{Fig:PD_gravity_pres_b}]{
        \includegraphics[clip,scale=0.174,trim=0 0cm 0cm 0]{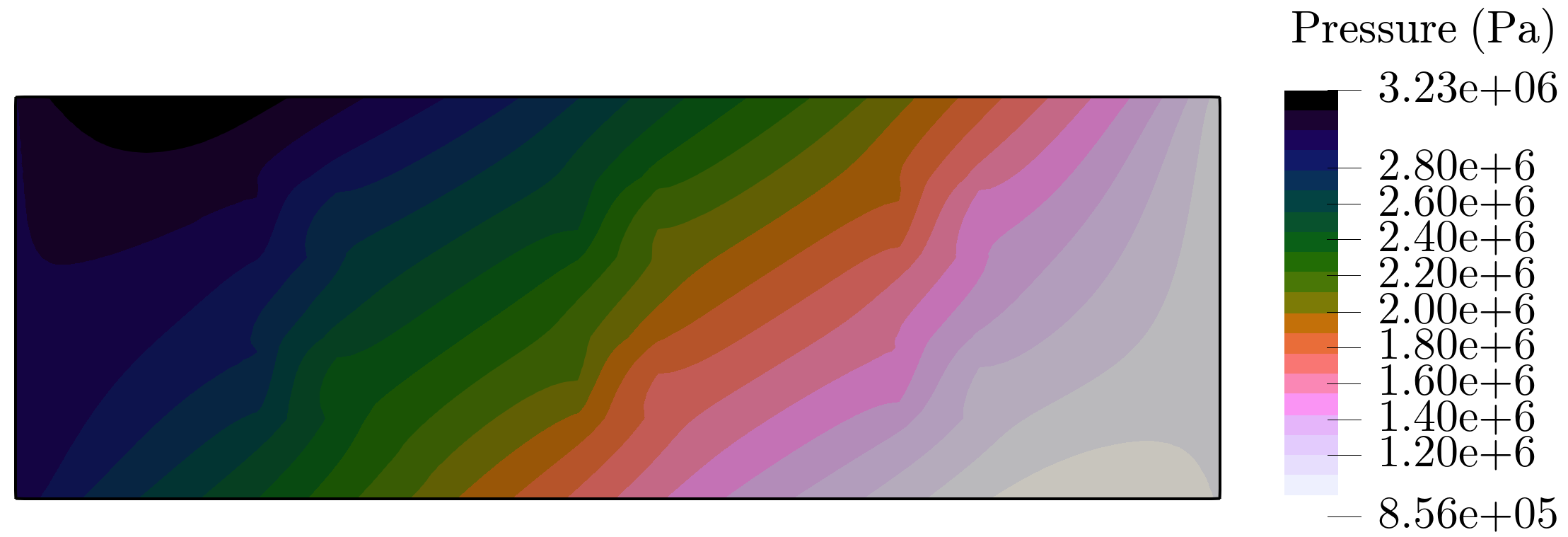}} \\
        \caption{\textsf{Pressure-driven flow in non-homogeneous domain with gravitational force:~}This figure
            shows the wetting phase pressure contours at $t=30$ days for $\mathrm{Gr}=0$ and $\mathrm{Gr}=0.4$.
        \label{Fig:PD_Gravity_pres_main}}
\end{figure}

\begin{figure}
    \subfigure[No gravity; $t=30$ days \label{Fig:PD_gravity_vel_a}]{
        \includegraphics[clip,scale=0.175,trim=0 0cm 0cm 0]{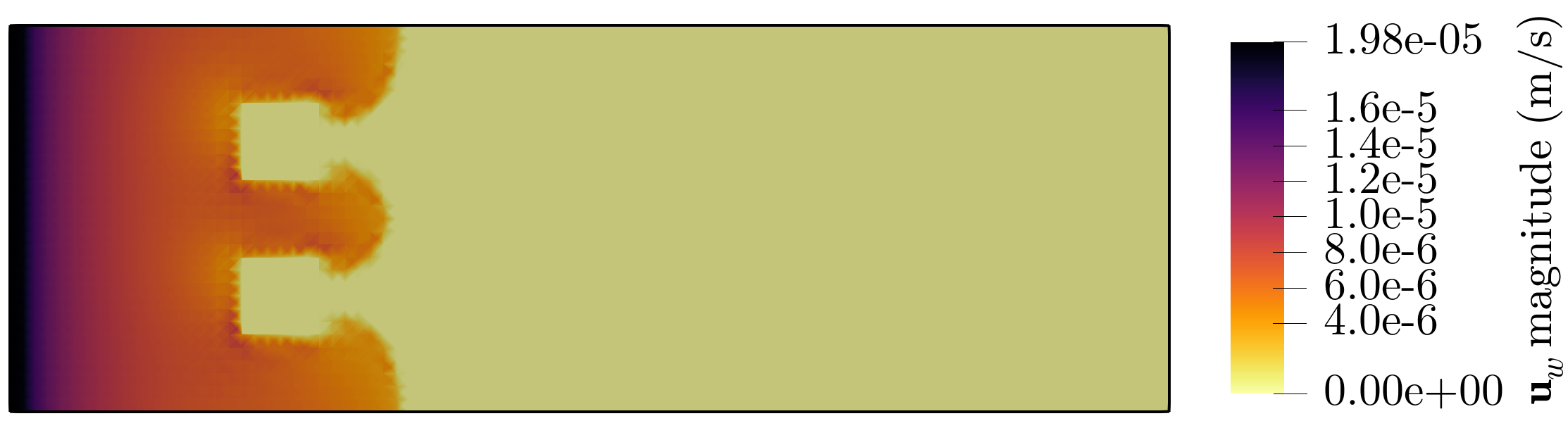}}
        \hspace{.25cm}
    \subfigure[With gravity; $t=30$ days \label{Fig:PD_gravity_vel_b}]{
        \includegraphics[clip,scale=0.171,trim=0 0cm 0cm 0]{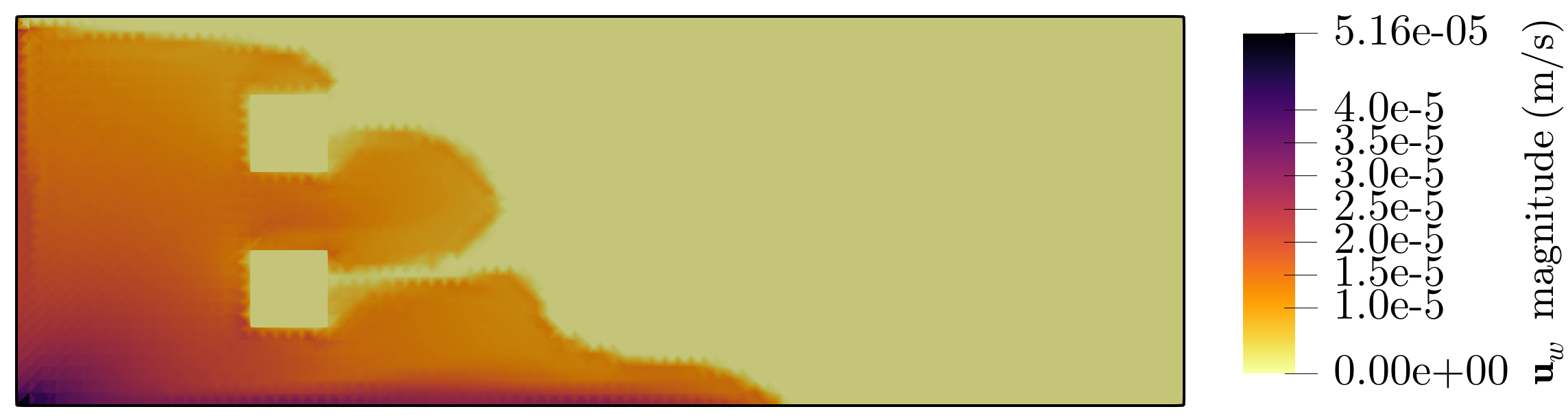}} \\
        \caption{\textsf{Pressure-driven flow in non-homogeneous domain with gravitational force:~}This figure
            shows the magnitude of wetting phase velocity solutions at $t=30$ days for $\mathrm{Gr}=0$ and $\mathrm{Gr}=0.4$.
        \label{Fig:PD_Gravity_vel_main}}
\end{figure}

\label{ssub:inhomogeneous_domain_with_gravity_fileds}
We now examine the performance of our limiting scheme in the presence of gravity field. For this problem, the domain of 
length $L_x=300$~\si{\meter} and $L_y=100$~\si{\meter} is partitioned into a crossed triangular mesh of size $h=10/3$ 
\si{\meter}. The gravity number $\mathrm{Gr}$ depends on the difference between phase densities. 
Permeability is set to $10^{-12}$ \si{\meter\squared}
everywhere except inside six square inclusions of length $20$ \si{\meter} centered at coordinates $(70,30)$, $(150,30)$,
$(230,30)$, $(70,70)$, $(150,70)$, and $(230,70)$ \si{\meter}, where the permeability is $10^3$ times smaller. Time
step is set to $\tau=1/12$ days and the final time is $T=30$ days. The proposed DG scheme with flux and slope limiters
is applied for two scenarios of $\mathrm{Gr}=0$ (i.e., no gravity) and $\mathrm{Gr}=0.4$. Figure~\ref{Fig:PD_Gravity_sat_main}
shows the saturation contours at the time $t=30$ days. In the presence of the gravitational body force, the horizontal symmetry of flow is broken and the wetting phase, which is heavier, starts to deposit at the bottom edge. As flow advances, the 
gravitational tongue at the bottom of domain becomes more distinct. For both problems, the limiting scheme exhibits 
satisfactory results with respect to the maximum principle. Pressure contours and the magnitude of the velocity field (see \eqref{eq:defuw}) are displayed
in Figure~\ref{Fig:PD_Gravity_pres_main} and \ref{Fig:PD_Gravity_vel_main}, respectively. The impact of gravity in both
solutions is noticeable.


\subsection{Quarter five-spot problem}%
\begin{figure}
    \subfigure[Unlimited DG; $t=4$ days\label{Fig:Q5_homogen_sat_a}]{
        \includegraphics[clip,scale=0.175,trim=0 0cm 0cm 0]{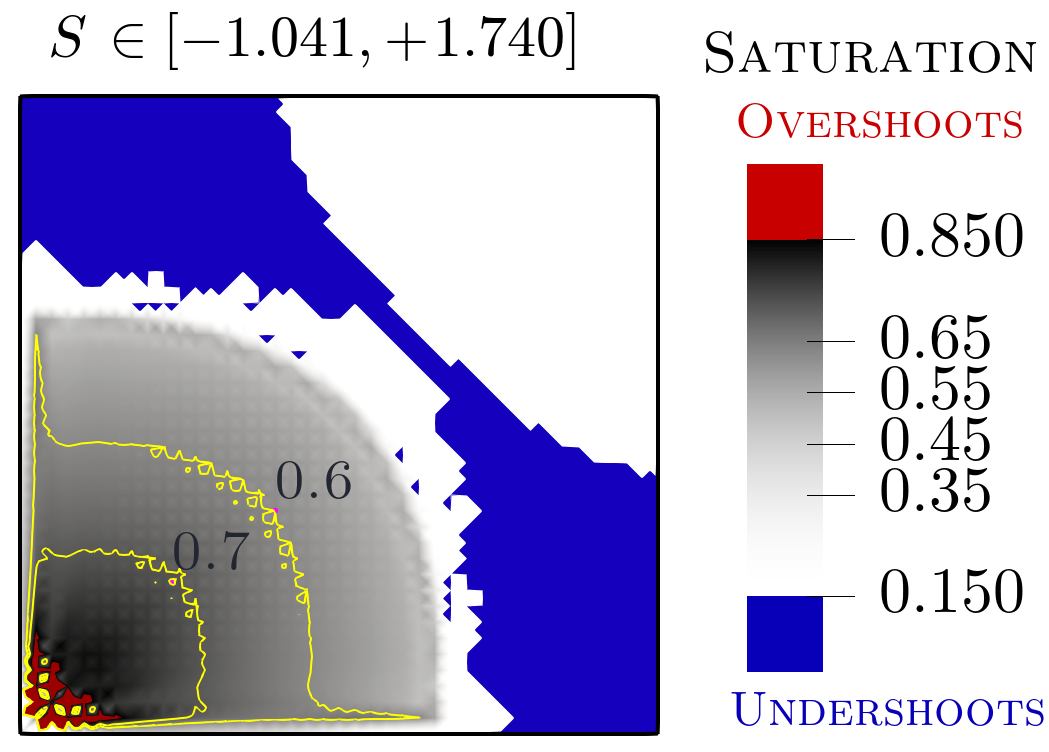}}
        \hspace{.3cm}
    \subfigure[Unlimited DG; $t=11$ days \label{Fig:Q5_homogen_sat_b}]{
        \includegraphics[clip,scale=0.175,trim=0 0cm 0cm 0]{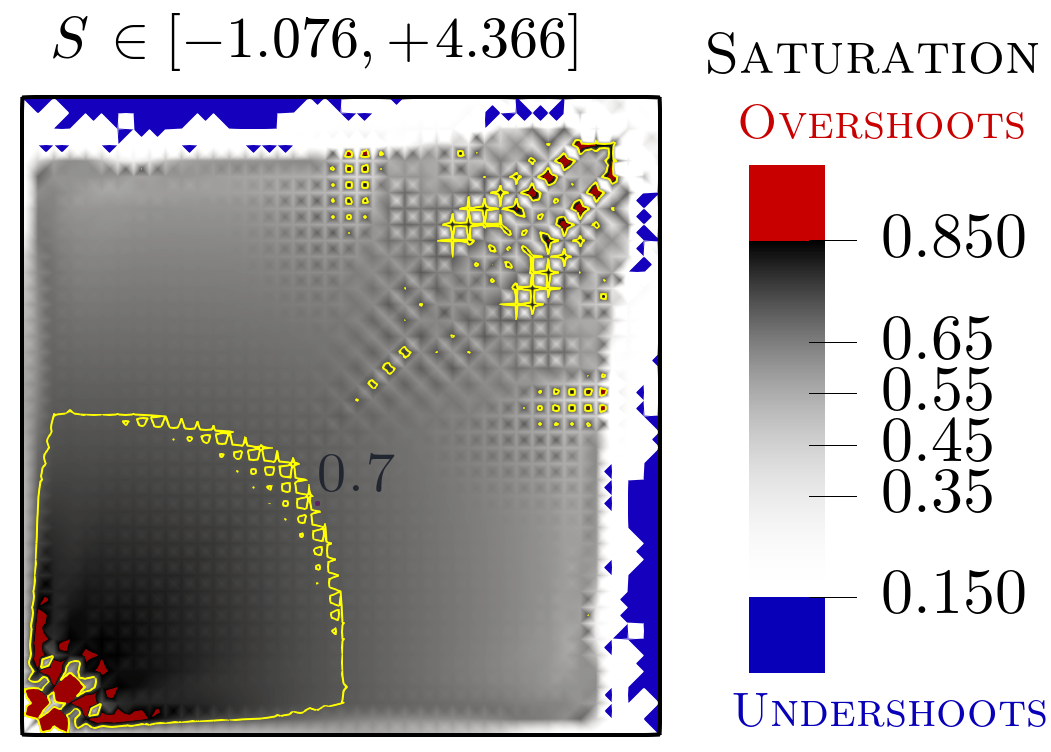}} \\
    \subfigure[DG+SL; $t=4$ days\label{Fig:Q5_homogen_sat_c}]{
        \includegraphics[clip,scale=0.175,trim=0 0cm 0cm 0]{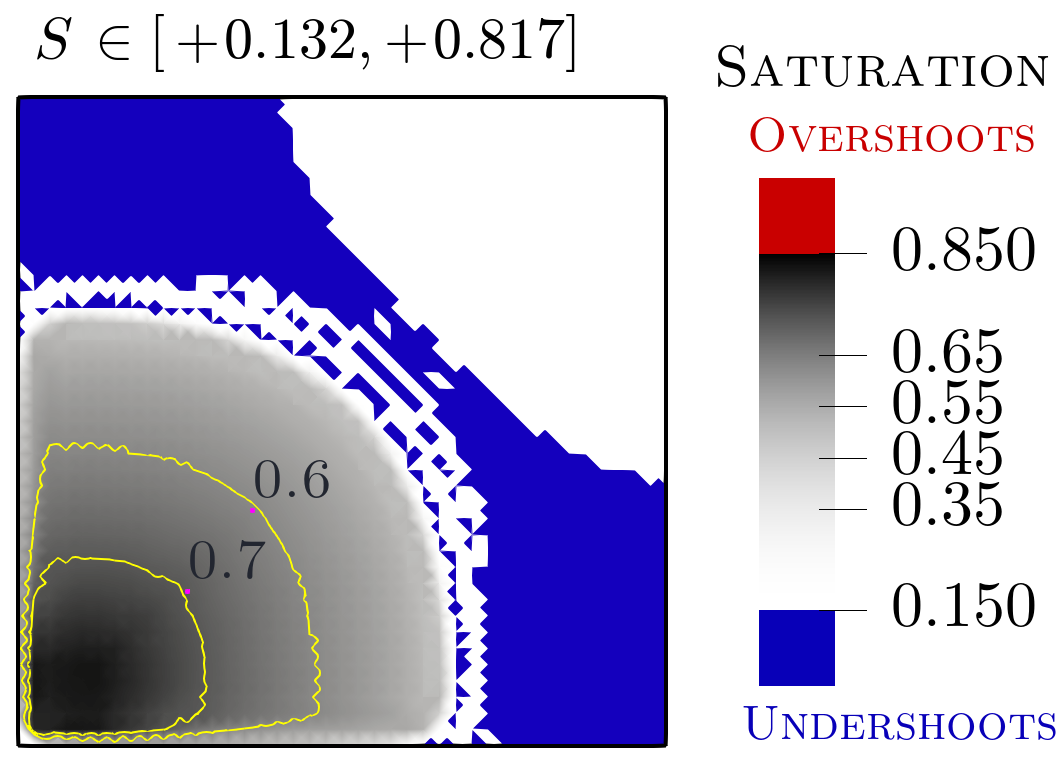}}
        \hspace{.3cm}
    \subfigure[DG+SL; $t=11$ days \label{Fig:Q5_homogen_sat_d}]{
        \includegraphics[clip,scale=0.175,trim=0 0cm 0cm 0]{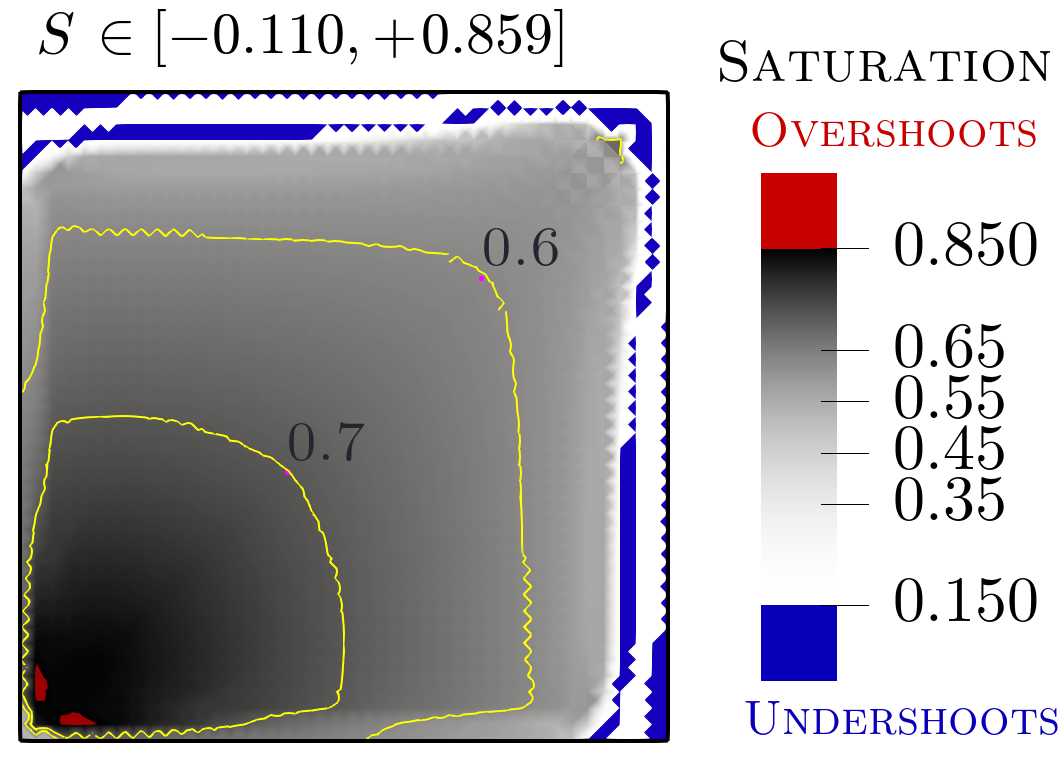}} \\
    \subfigure[DG+FL+SL; $t=4$ days\label{Fig:Q5_homogen_sat_e}]{
        \includegraphics[clip,scale=0.175,trim=0 0cm 0cm 0]{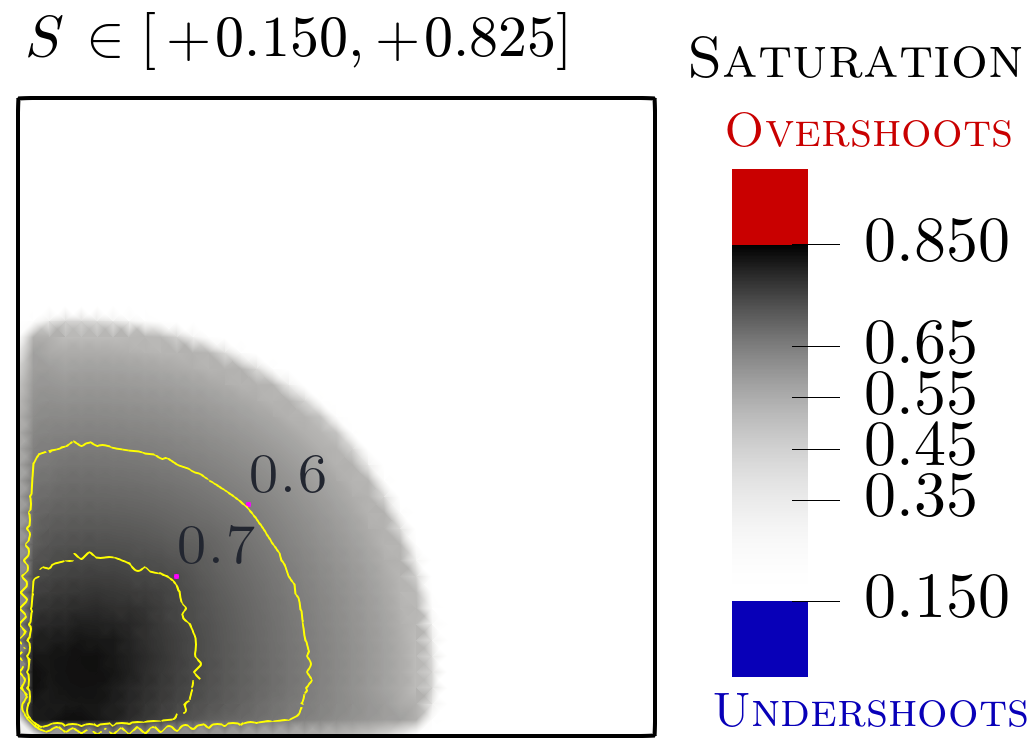}}
        \hspace{.3cm}
    \subfigure[DG+FL+SL; $t=11$ days \label{Fig:Q5_homogen_sat_f}]{
        \includegraphics[clip,scale=0.175,trim=0 0cm 0cm 0]{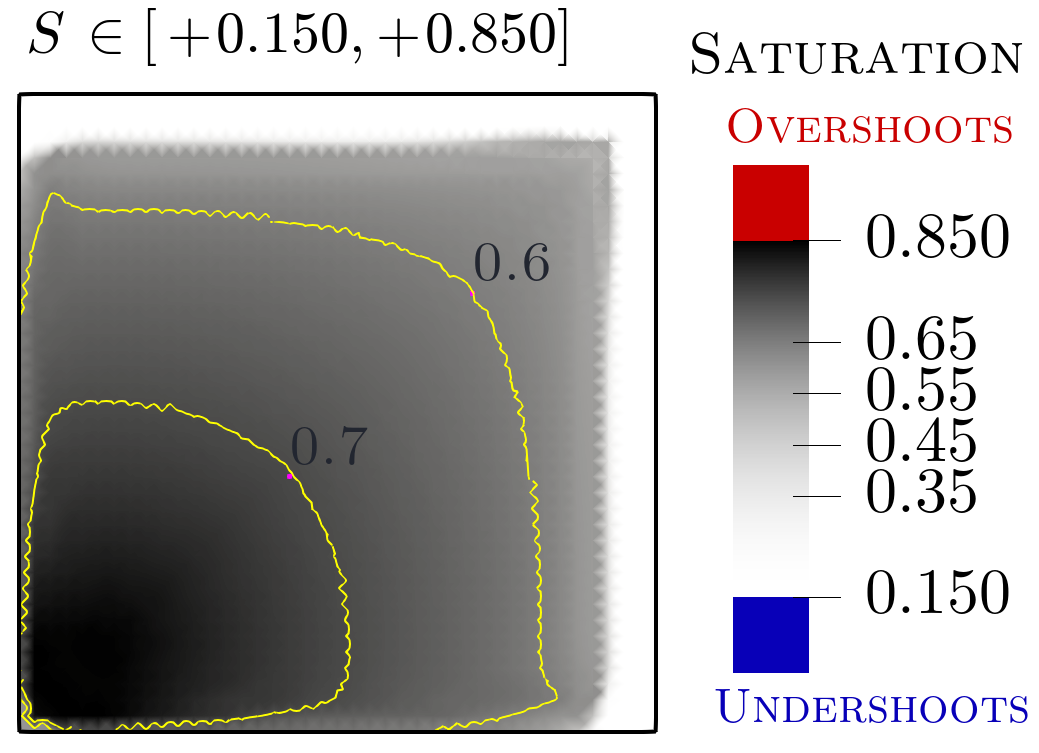}} \\
        \caption{\textsf{Quarter five-spot problem in homogeneous domain:}~This figure shows the saturation contours obtained with DG (top), DG+SL (middle), DG+FL+SL (bottom) at two different time steps. The proposed DG+FL+SL scheme, unlike the two other schemes, respects maximum principle throughout the simulation.
        \label{Fig:Q5_homogen_sat_main}}
\end{figure}

\begin{figure}
    \subfigure[Limited DG; $4$ days \label{Fig:Q5_homogen_pres_a}]{
        \includegraphics[clip,scale=0.175,trim=0 0cm 0cm 0]{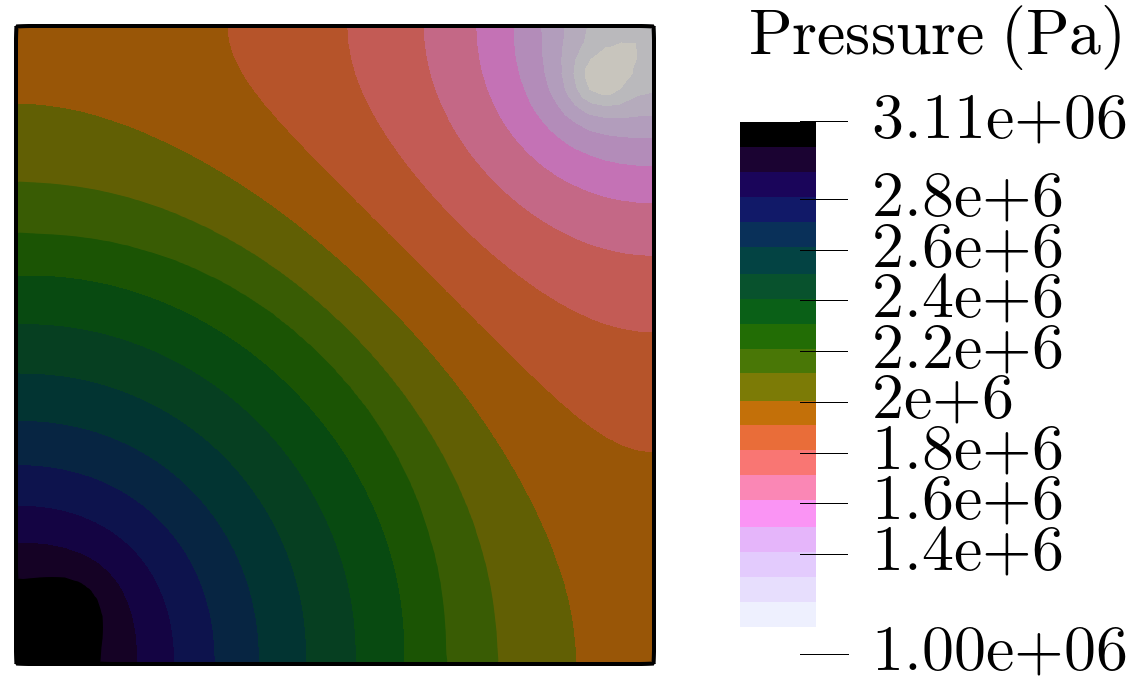}}
        \hspace{.3cm}
    \subfigure[Limited DG; $11$ days \label{Fig:Q5_homogen_pres_b}]{
        \includegraphics[clip,scale=0.175,trim=0 0cm 0cm 0]{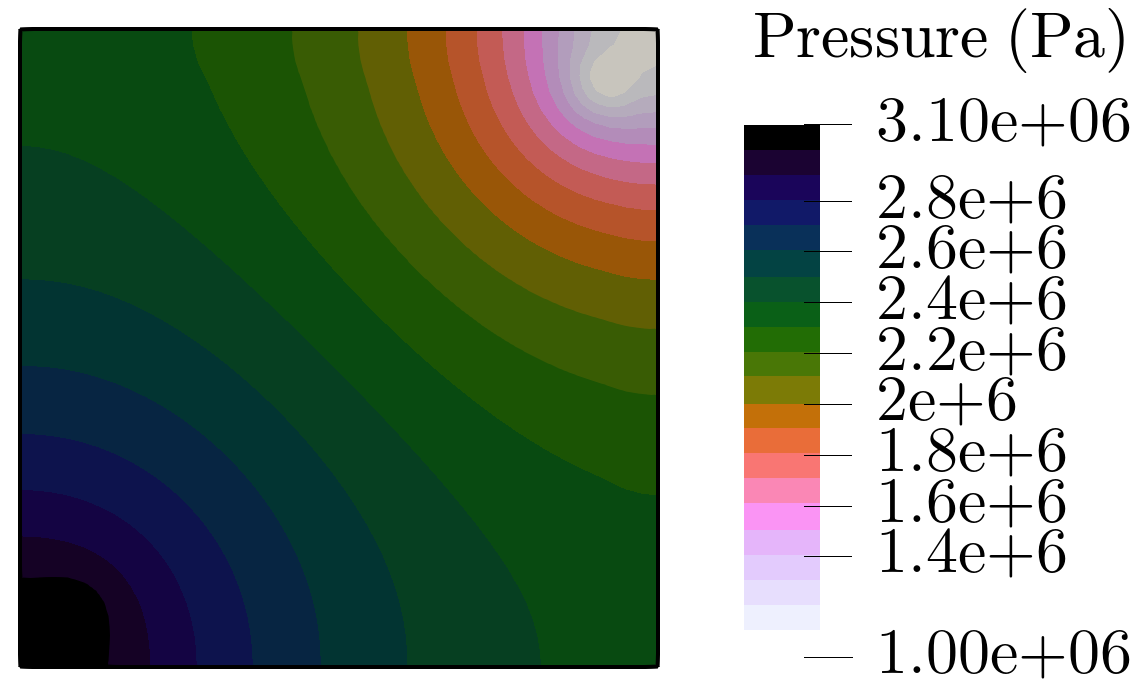}} \\
        \caption{\textsf{Quarter five-spot problem in homogeneous domain:}~This figure shows the wetting phase pressures obtained under
            the proposed DG+FL+SL algorithm.
        \label{Fig:Q5_homogen_pres_main}}
\end{figure}

\begin{figure}
    \subfigure[Limited DG; $4$ days \label{Fig:Q5_homogen_vel_a}]{
        \includegraphics[clip,scale=0.17,trim=0 0cm 0cm 0]{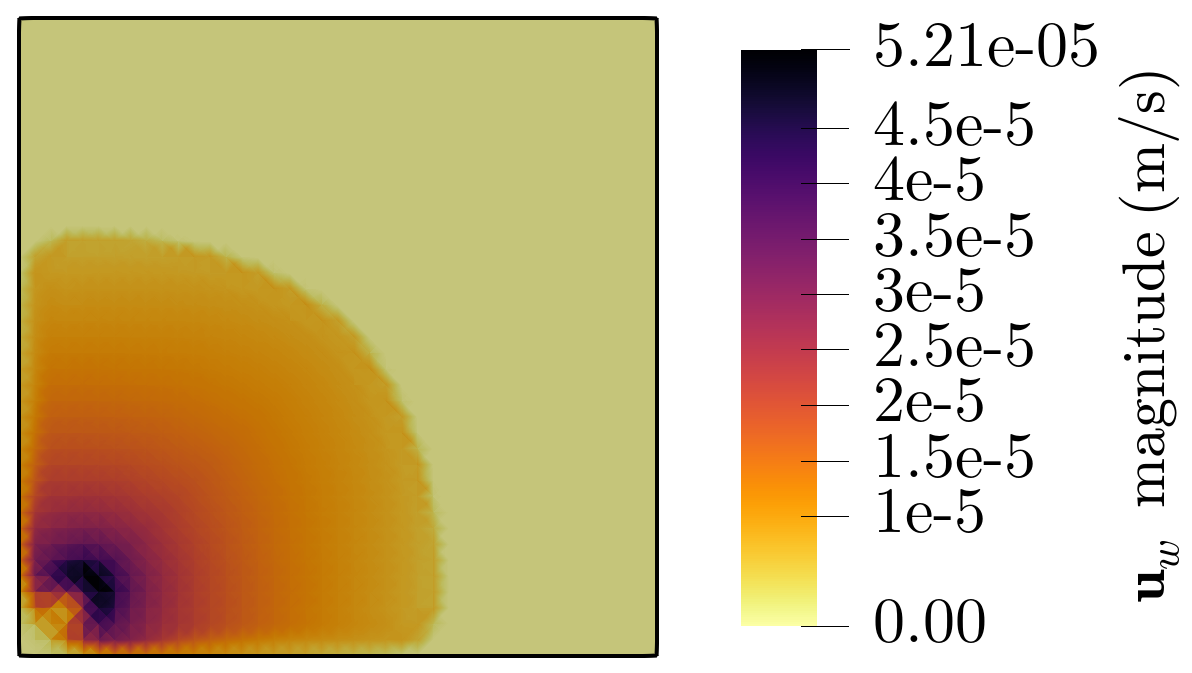}}
        \hspace{.3cm}
    \subfigure[Limited DG; $11$ days \label{Fig:Q5_homogen_vel_b}]{
        \includegraphics[clip,scale=0.17,trim=0 0cm 0cm 0]{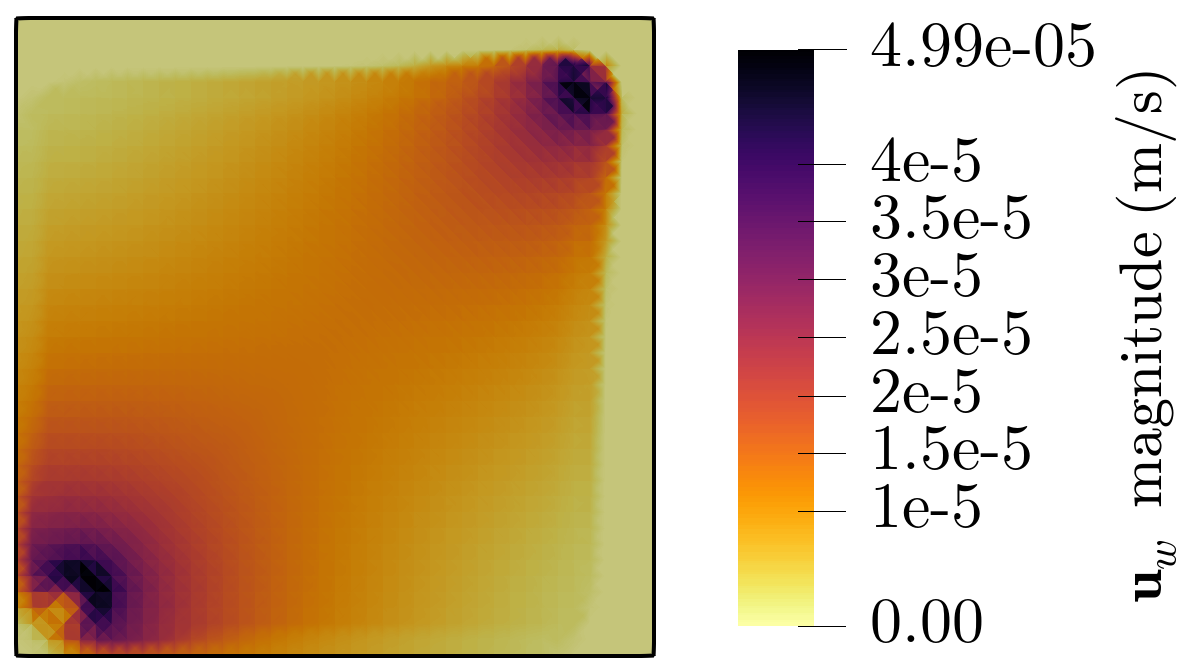}} \\
        \caption{\textsf{Quarter five-spot problem in homogeneous domain:}~This figure shows the wetting phase velocities obtained under
            the proposed limiting algorithm.
        \label{Fig:Q5_homogen_vel_main}}
\end{figure}

\begin{figure}
    \subfigure[Unlimited DG \label{Fig:Q5_homogen_error_a}]{
        \includegraphics[clip,scale=0.135,trim=0 0cm 0cm 0]{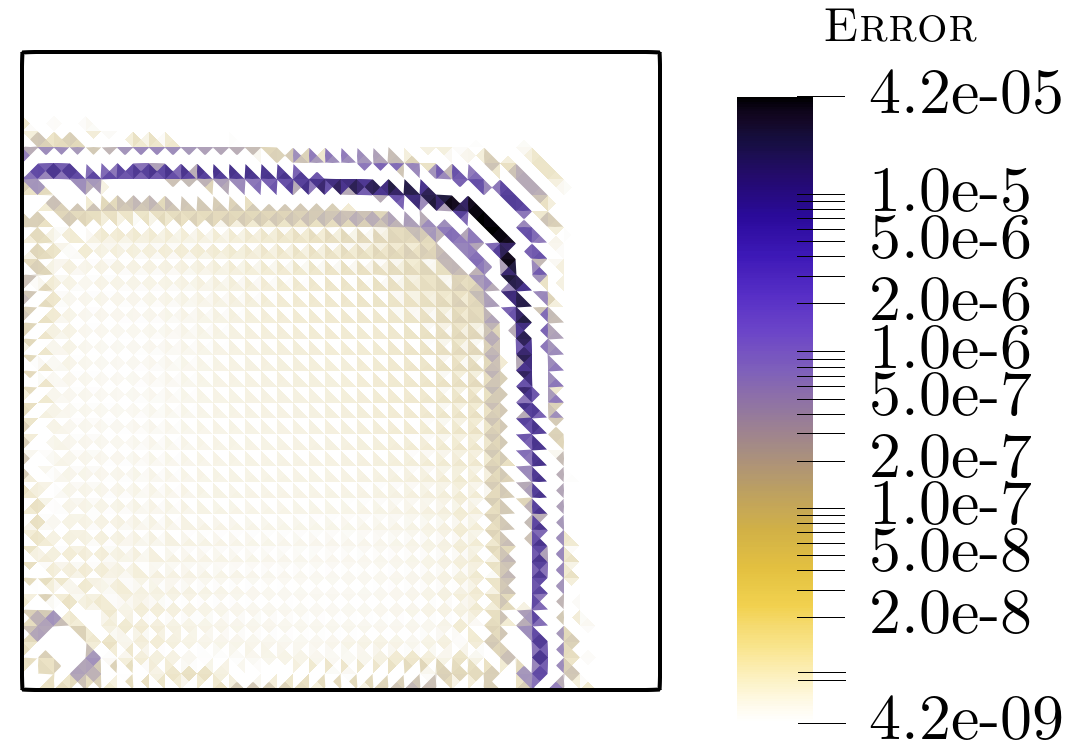}}
        \hspace{.05cm}
    \subfigure[DG+SL \label{Fig:Q5_homogen_error_b}]{
        \includegraphics[clip,scale=0.135,trim=0 0cm 0cm 0]{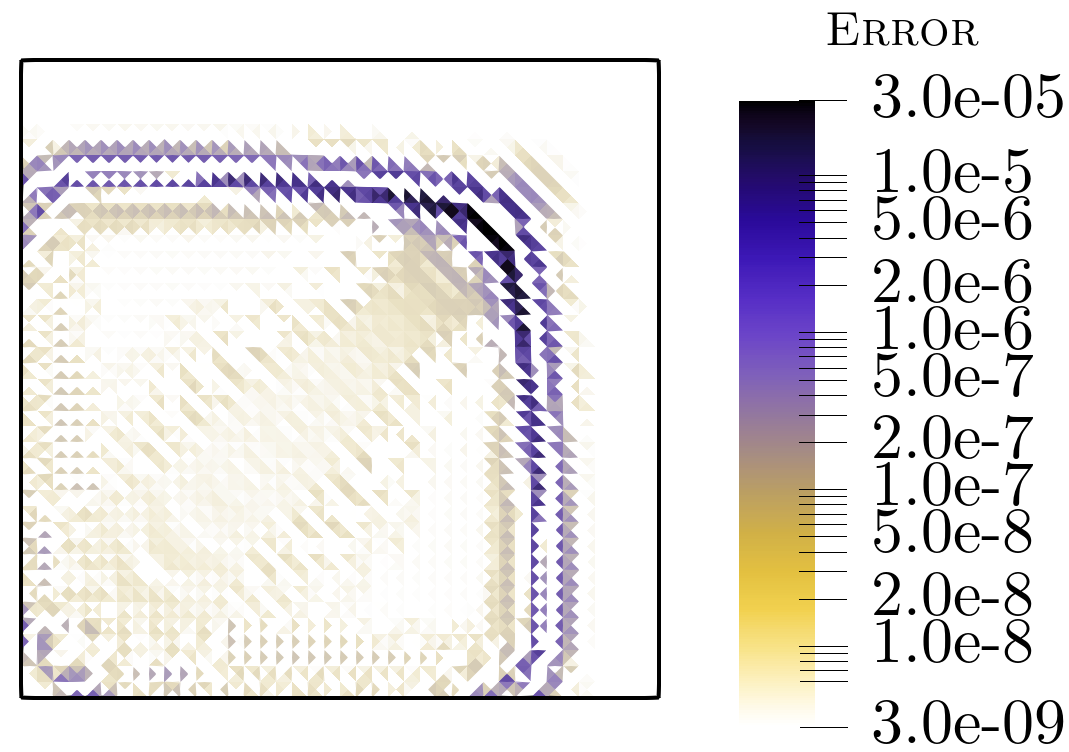}} 
        \hspace{.05cm}
    \subfigure[DG+FL+SL \label{Fig:Q5_homogen_error_c}]{
        \includegraphics[clip,scale=0.135,trim=0 0cm 0cm 0]{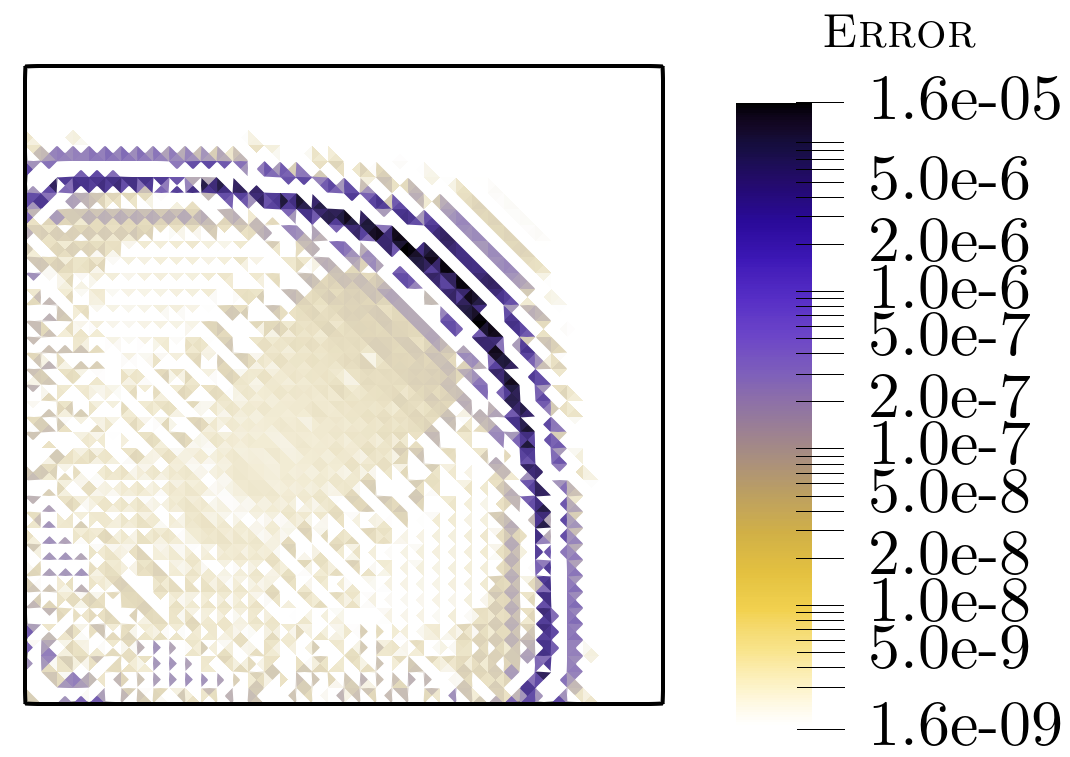}} 
        \caption{\textsf{Quarter five-spot problem in homogeneous domain:}~This figure shows the element-wise mass balance error
            $\mathcal{B}(E)$ at time $t=7$ days. Errors remain small regardless of the scheme used.
        \label{Fig:Q5_homogen_error_main}}
\end{figure}

\label{sub:quarter_five_spot_problem}
In this section, we evaluate the performance of limiters in the presence of wells.
We take a square computational domain of size $L=100$ \si{\meter} with permeability of $10^{-12}$ \si{\meter\squared} 
everywhere. The domain is partitioned into a crossed triangular mesh of size $h=2.5$ \si{\meter}. 
As shown in Figure~\ref{Fig:Sch_Q5}, no flow boundary condition is prescribed on $\partial\Omega$ and the flow 
is driven by injector/producers (source/sink functions) (see \eqref{eq:sourcesink}).
The injection saturation is set to $s_\mathrm{in}=0.85$ and the injection and production flow rate of wells 
are set to:
\begin{align}
    \int_{\Omega}\bar{q}=\int_{\Omega}\underline{q}= 9.8437\times10^{-4},
\end{align}
where $\bar{q}$ is piecewise constant on $[5,12.5]^2$ \si{\meter\squared} and $\bar{q}=0$ elsewhere and $\underline{q}$
is piecewise constant on $[87.5,95]^2$ \si{\meter\squared} and $\underline{q}=0$ elsewhere. 
Time step is set to $\tau=0.05$ days and the simulation advances up to  $T=11$ days. Figure~\ref{Fig:Q5_homogen_sat_main} depicts the wetting phase saturations at two different time of $t=4$ and $t=11$
days for unlimited DG, DG+SL, and DG+FL+SL schemes. Unlimited DG
returns oscillatory non-monotone solutions and violations of the maximum principle are noticeable in the vicinity of injector and after the saturation front. The slope limiter improves the accuracy of the solution by reducing the amount of overshoot near
the injection well, but falls short in generating bound-preserving solutions throughout the simulation. On the other hand, the proposed DG+FL+SL limiting scheme returns monotone solutions and fully eliminates undershoots (i.e., blue-colored cells) and overshoots (red-colored cells). Table~\ref{Tab:Q5_compare} displays the minimum and maximum values of the saturations over the whole simulation, as well as the amounts
of overshoot and undershoot in percentages. 
Figure~\ref{Fig:Q5_homogen_pres_main} and \ref{Fig:Q5_homogen_vel_main} show the wetting phase pressure and the magnitude of the velocity (defined by \eqref{eq:defuw})  under the proposed limiting scheme. As time progresses, higher pressure differences build up near the producer and hence the magnitude of the velocity increases in that region. Finally, we investigate the impact of applying slope and flux limiters on the local mass conservation properties of the DG formulation. Following \eqref{eq:defavg}, we denote by $\overline{\varphi}$
the element-wise average of a function $\varphi$. By choosing a test
function equal to $1/\vert E\vert$ on one element and $0$ elsewhere, we obtain
the local mass balance of an element
$E\in \mathcal{E}_{h}$ at time $t_n$ as follows: 
\begin{align}
    \mathcal{B}(E)= &
    \frac{1}{\tau}\Big(\overline{\rho_w(P_{n+1})\phi(P_{n+1}) S_{n+1}}
  - \overline{\rho_w(P_{n})\phi(P_{n}) S_n} \Big)
  +
  \frac{1}{|E|}
  \sum_{e\subset\partial E} \mathcal{H}_{n+1,E}(e)
    \nonumber\\&-
\overline{\rho_w(P_n) f_w(s_{\mathrm{in}}) \bar{q}_E
-\rho_w(P_n) f_w(S_{n}|_E)\underline{q}_E}.
  \label{eq:updateFLL} 
\end{align}
The values of $\mathcal{B}(E)$ are displayed in Figure~\ref{Fig:Q5_homogen_error_main} for unlimited DG, DG+SL and DG+FL+SL
at time $t=7$ days. 
We observe that the mass error is of the order of $10^{-9}$ everywhere except in a small neighborhood of the
saturation front where the mass error increases to $10^{-5}$. The slope and flux limiters do not change the
magnitude of the local mass error.
\begin{table}
\centering
\caption{This table shows the efficacy of the limiting schemes when applied to the quarter five-spot flow problem with homogeneous domain.\label{Tab:Q5_compare}}
\vspace{-0.25cm}
\resizebox{0.75\textwidth}{!}{%
    \begin{tabular}{l@{\hskip 0.5in}cc@{\hskip 0.5in}cc@{\hskip 0.5in}cc}
\Xhline{2\arrayrulewidth}\\[-0.9em]
               & \multicolumn{2}{l}{DG}                                     & \multicolumn{2}{l}{DG+SL} & \multicolumn{2}{l}{DG+FL+SL} \\ \cmidrule(l{-0.1in}r{0.42 in}){2-3} \cmidrule(l{-0.1in}r{0.42 in}){4-5}\cmidrule(l{-0,1in}r{0.0 in}){6-7}
               & value                        & \%                  & value     & \%    & value     & \%  \\[-0.9em]  \\ \hline\\[-0.9em]
minimum saturation & -1.39 & 220 & -0.1137 & 37.7 & 0.15 & 0 \\
maximum saturation & 10.51 &  1380 & 0.86      &   1.4             & 0.85      & 0     \\[-0.9em] \\
\Xhline{2\arrayrulewidth}
\end{tabular}%
}
\end{table}

\subsection{Effect of rock and phases compressibility factors}%
\begin{figure}
    \subfigure[Various rock compressibilites \label{Fig:PD_Rock_compress}]{
        \includegraphics[clip,scale=0.28,trim=0 0cm 0cm 0]{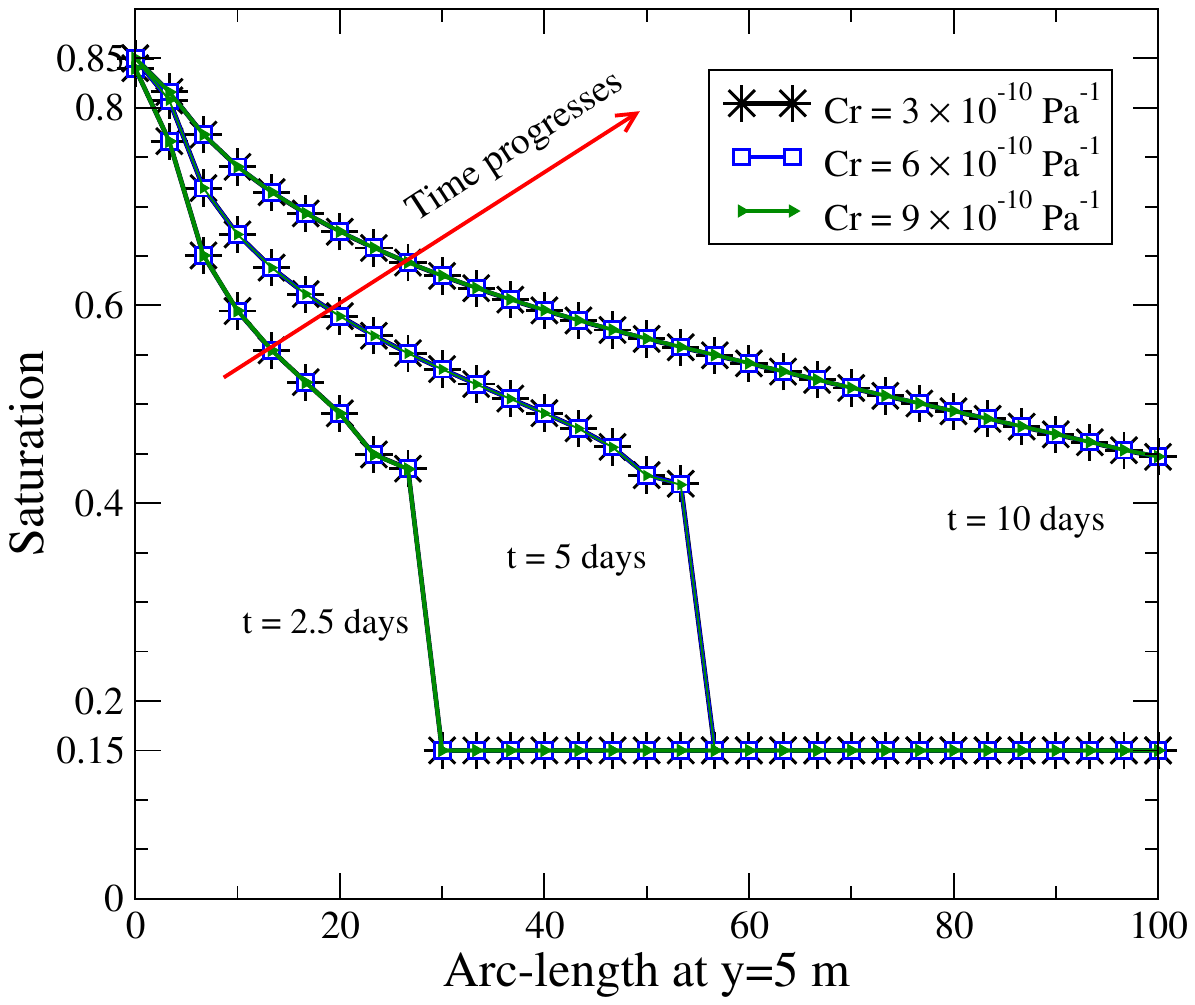}} 
    \subfigure[Compressible vs incompressible fluid phases \label{Fig:PD_compressibility}]{
        \includegraphics[clip,scale=0.28,trim=0 0cm 0cm 0]{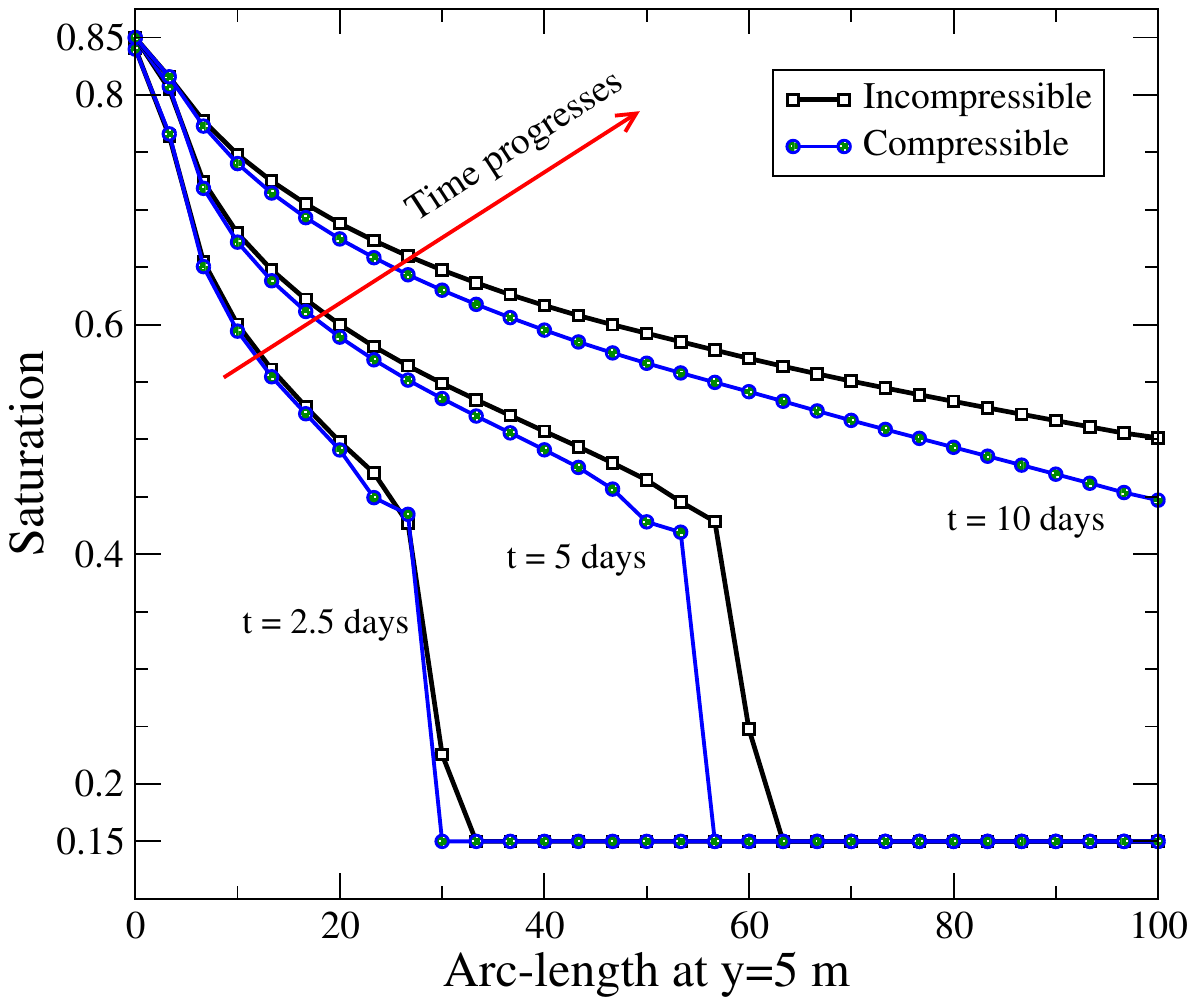}}
        \caption{\textsf{Effect of compressibility factors for a pressure-driven flow problem:~}
            This figure shows the evolution of
            saturation profiles obtained from the limited DG scheme for 
            (a) three different rock compressibility factors with fixed $c_w=10^{-10}$ and
            $c_{\ell}=10^{-6}$ and (b) for different phase compressibility factors with fixed $c_{r}=6\times10^{-10}$.
            Rock compressibility has negligible effect on solutions. It is also evident that fluid compressibility slows
            down the propagation of wetting phase saturation.
        \label{Fig:PD_compare_main}}
\end{figure}

\begin{figure}
    \subfigure[Compressible case; saturation \label{Fig:Q5_compare_a}]{
        \includegraphics[clip,scale=0.15,trim=0 0cm 0cm 0]{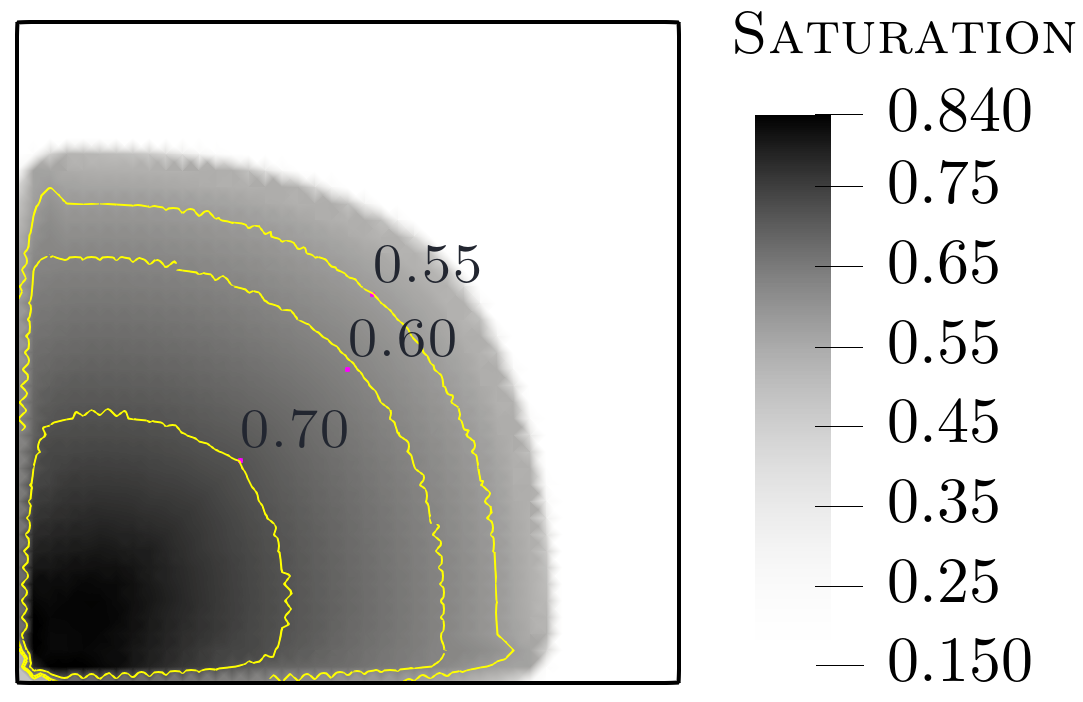}}
        \hspace{.25cm}
    \subfigure[Incompressible case; saturation \label{Fig:Q5_compare_b}]{
        \includegraphics[clip,scale=0.15,trim=0 0cm 0cm 0]{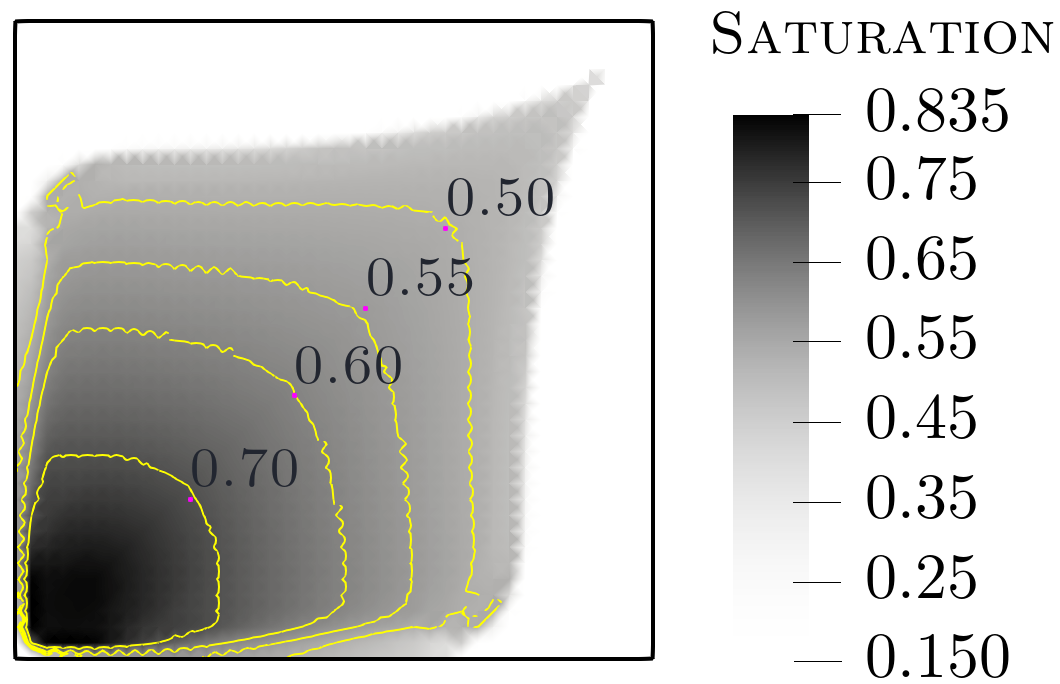}} \\
    \subfigure[Compressible case; pressure \label{Fig:Q5_compare_c}]{
        \includegraphics[clip,scale=0.15,trim=0 0cm 0cm 0]{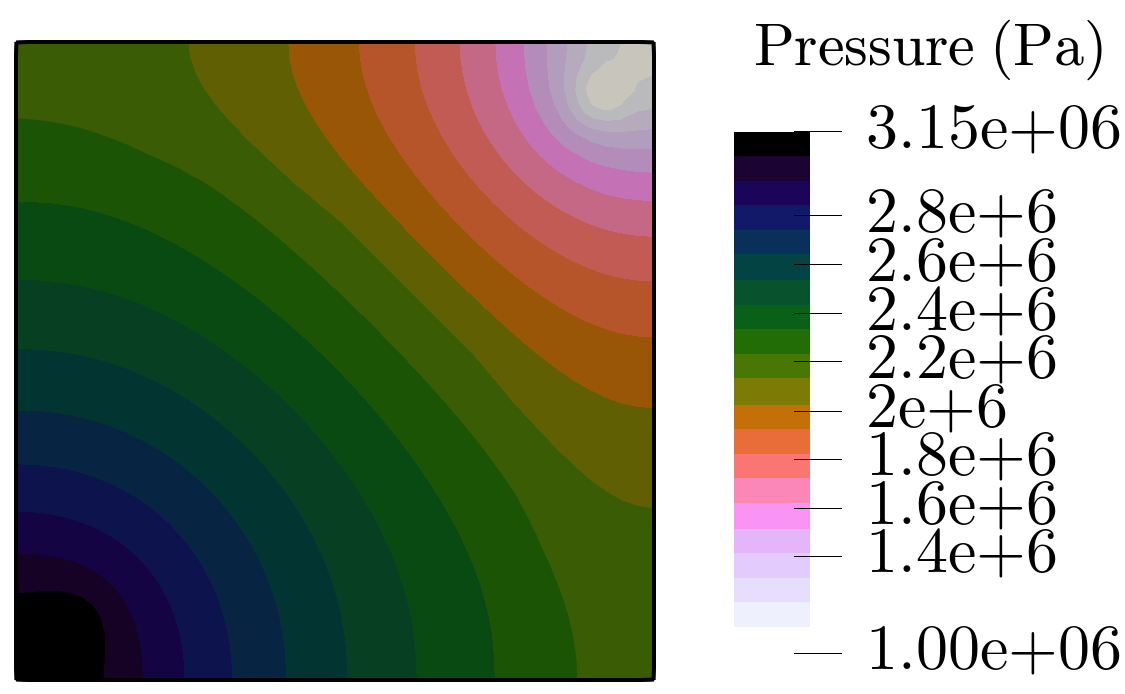}}
        \hspace{.25cm}
    \subfigure[Incompressible case; pressure \label{Fig:Q5_compare_d}]{
        \includegraphics[clip,scale=0.15,trim=0 0cm 0cm 0]{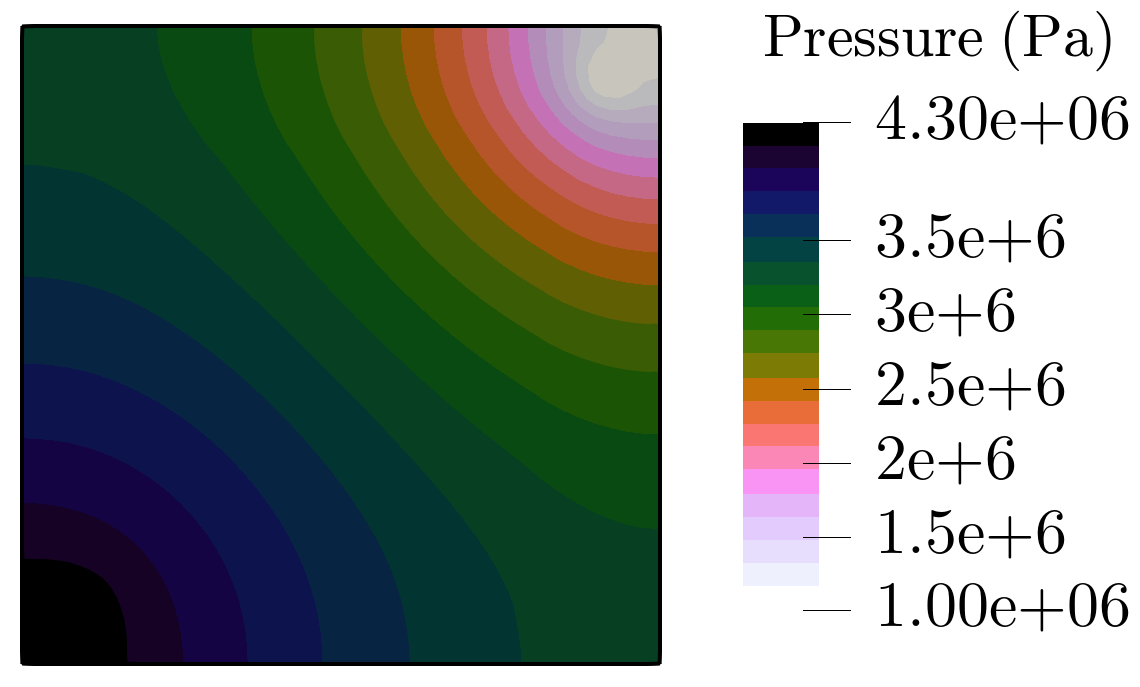}} \\
    \subfigure[Compressible case; velocity \label{Fig:Q5_compare_e}]{
        \includegraphics[clip,scale=0.15,trim=0 0cm 0cm 0]{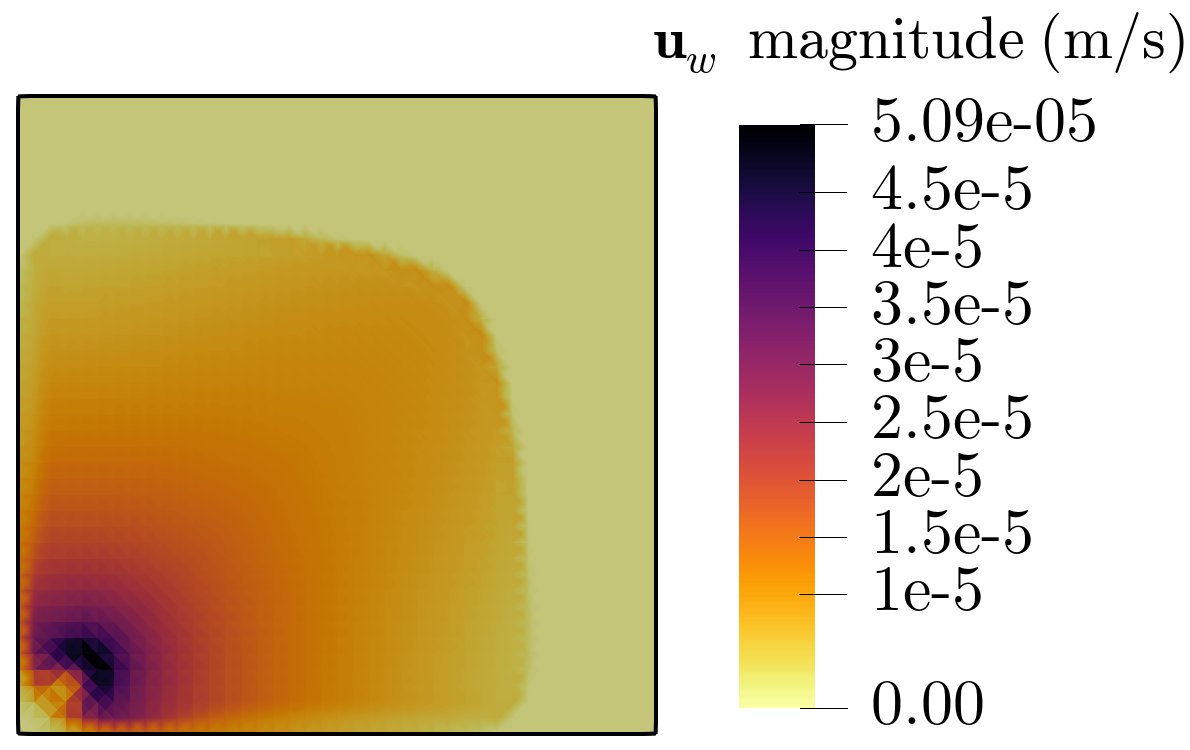}}
        \hspace{.25cm}
    \subfigure[Incompressible case; velocity \label{Fig:Q5_compare_f}]{
        \includegraphics[clip,scale=0.15,trim=0 0cm 0cm 0]{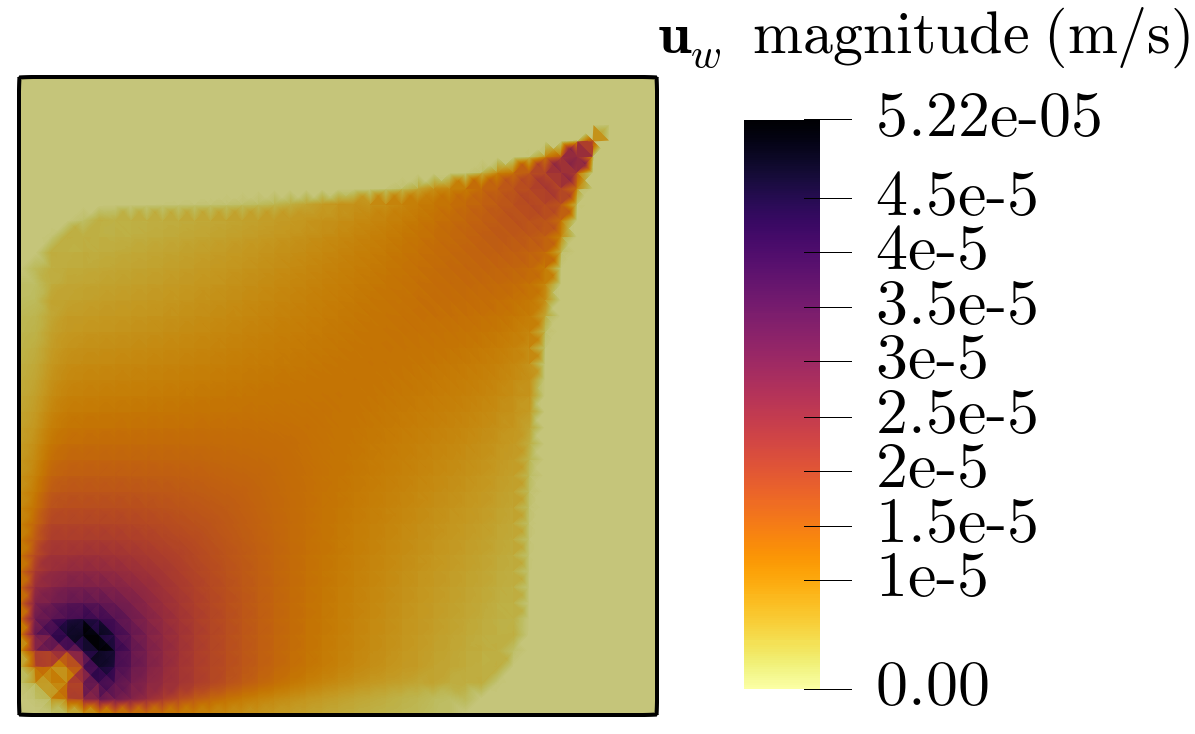}} \\
        \caption{\textsf{Quarter five-spot problem in homogeneous domain:}~
            This figure compares solutions for compressible and incompressible phases obtained from limited DG scheme  
            at $t=6.5$ days. Wetting phase fluid propagates faster for the incompressible case.
        \label{Fig:Q5_compare_main}}
\end{figure}

\label{sub:effect_of_compressibility}
In this section we use the limited DG scheme to study the impact of compressibility factors on pressure-driven flow problem discussed in Section~\ref{ssub:homogeneous_doamin} and the quarter-five spot problem discussed in Section~\ref{sub:quarter_five_spot_problem}. We first set the rock compressibility to take three different physical values of $3\times10^{-10}$, $6\times10^{-10}$, and $9\times10^{-10}$ \si{\per\pascal}
\citep{baker20157} and examine two cases of compressible phases (with $c_w=10^{-10}$ and $c_{\ell}=10^{-6}$) and incompressible phases (i.e., $c_w=c_{\ell}=0$). Other parameters and boundary conditions remain unchanged. Figure~\ref{Fig:PD_Rock_compress} displays the saturation solution for compressible phases along the line $y=5$ \si{\meter}. We observe that the rock compressibility factor yields  negligible changes in solutions. From Figures
\ref{Fig:PD_compressibility} and \ref{Fig:Q5_compare_main} it can be seen that the wetting phase floods the domain faster
in the incompressible case than in  the compressible case.

\subsection{Quarter five-spot problem with a highly anisotropic permeability}%
\begin{figure}[htpb]
    \centering
    \includegraphics[width=0.4\linewidth]{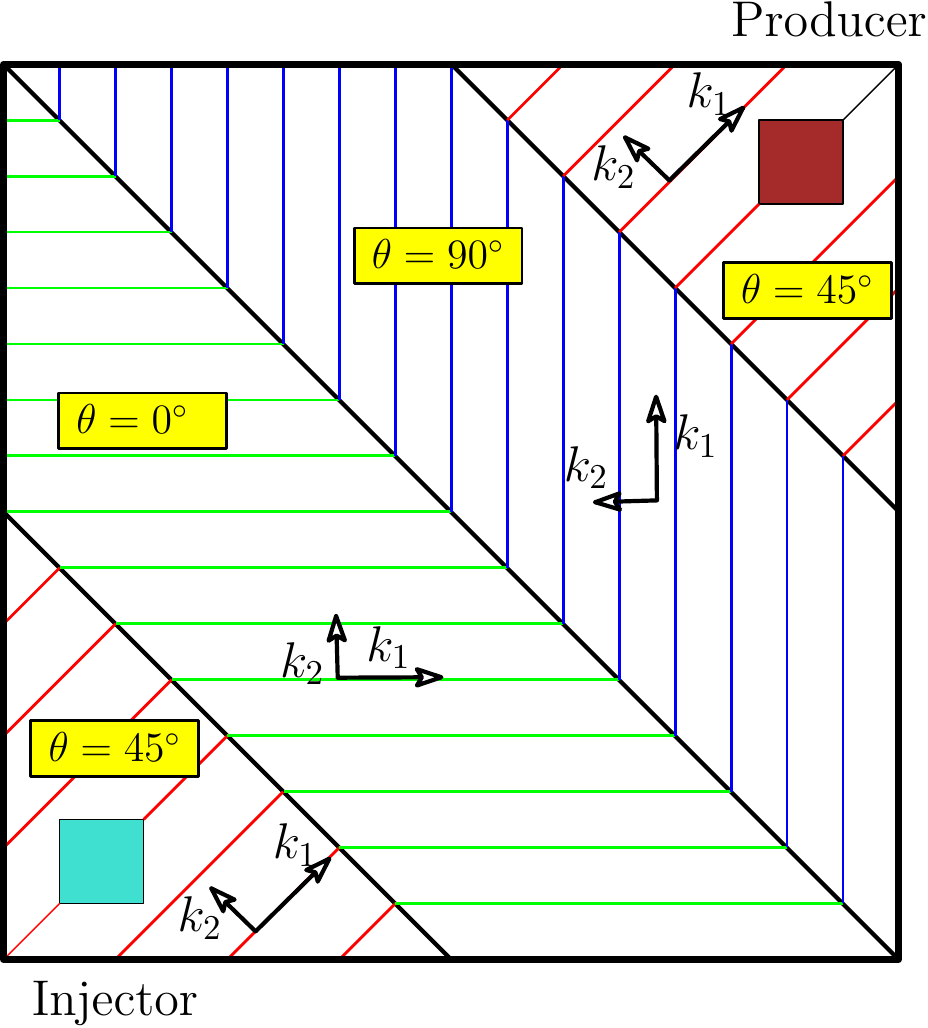}
    \caption{Quarter five-spot problem with a discontinuous full tensor and highly anisotropic permeability. Principal permeabilities are set to $k_1=2.25\times10^{-12}$ and $k_2=2.25\times10^{-14}$.}%
    \label{fig:aniso_schematic}
\end{figure}

\begin{figure}
    \subfigure[$t=1.25$ days; saturation \label{Fig:Q5_anisotropy_sat_a}]{
        \includegraphics[clip,scale=0.13,trim=0 0cm 0cm 0]{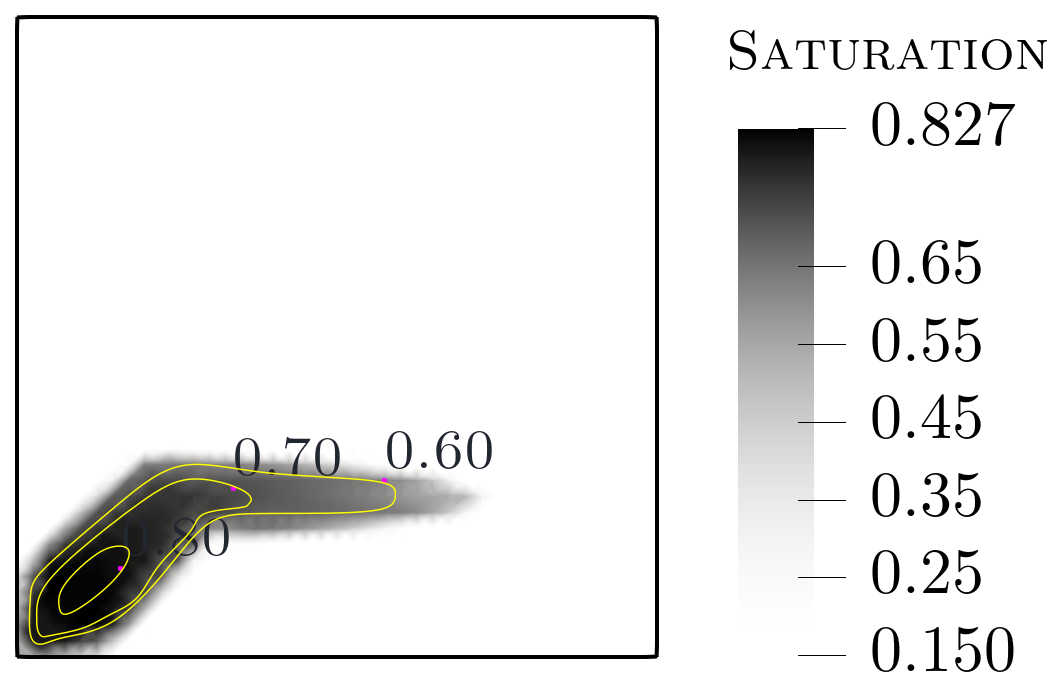}}
        \hspace{.07cm}
    \subfigure[$t=2.5$ days; saturation \label{Fig:Q5_anisotropy_sat_b}]{
        \includegraphics[clip,scale=0.13,trim=0 0cm 0cm 0]{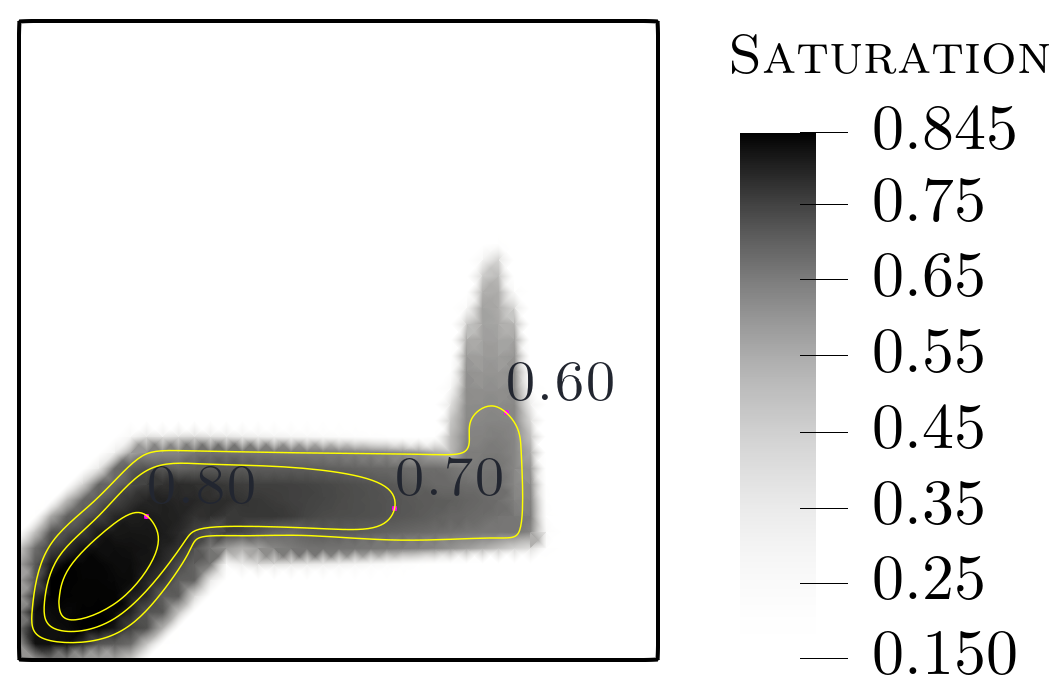}} 
        \hspace{.05cm}
    \subfigure[$t=5$ days; saturation \label{Fig:Q5_anisotropy_sat_c}]{
        \includegraphics[clip,scale=0.13,trim=0 0cm 0cm 0]{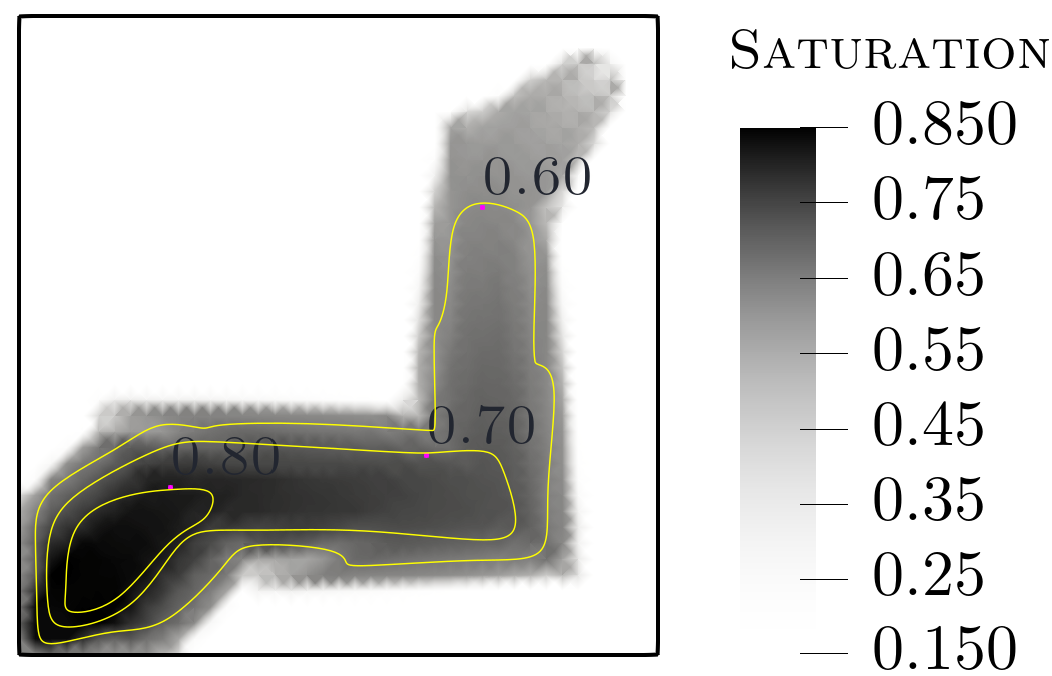}} \\
    \subfigure[$t=1.25$ days; pressure \label{Fig:Q5_anisotropy_pres_a}]{
        \includegraphics[clip,scale=0.125,trim=0 0cm 0cm 0]{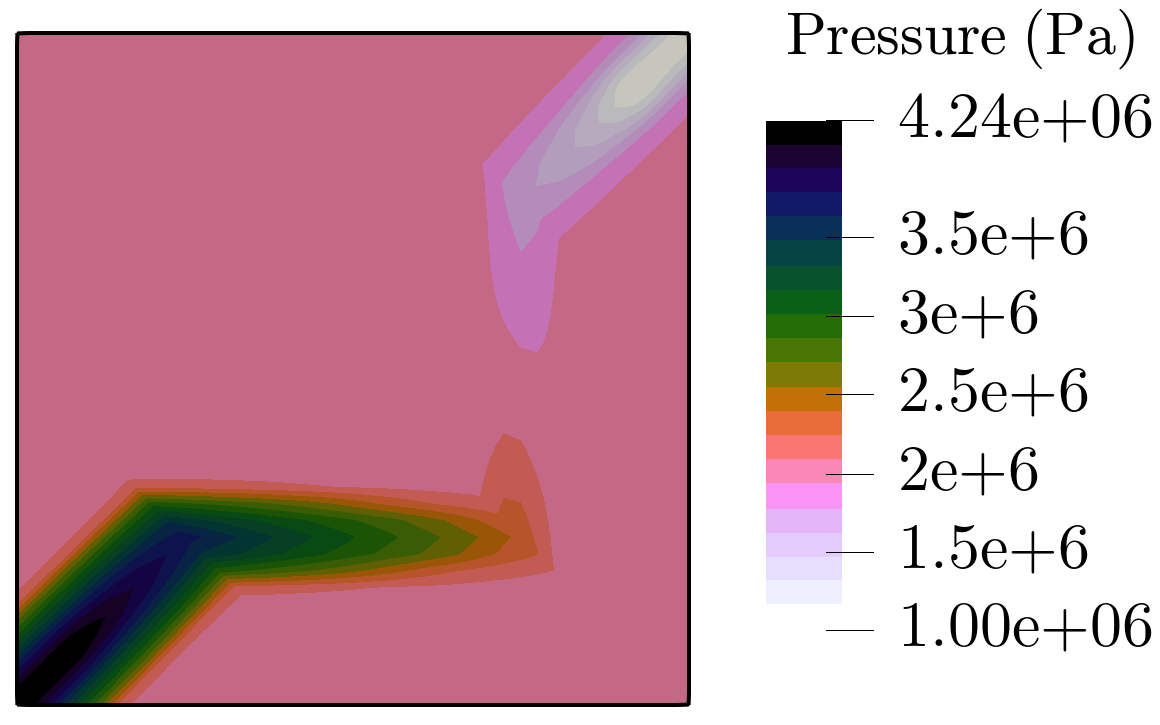}}
        \hspace{.05cm}
    \subfigure[$t=2.5$ days; pressure \label{Fig:Q5_anisotropy_pres_b}]{
        \includegraphics[clip,scale=0.13,trim=0 0cm 0cm 0]{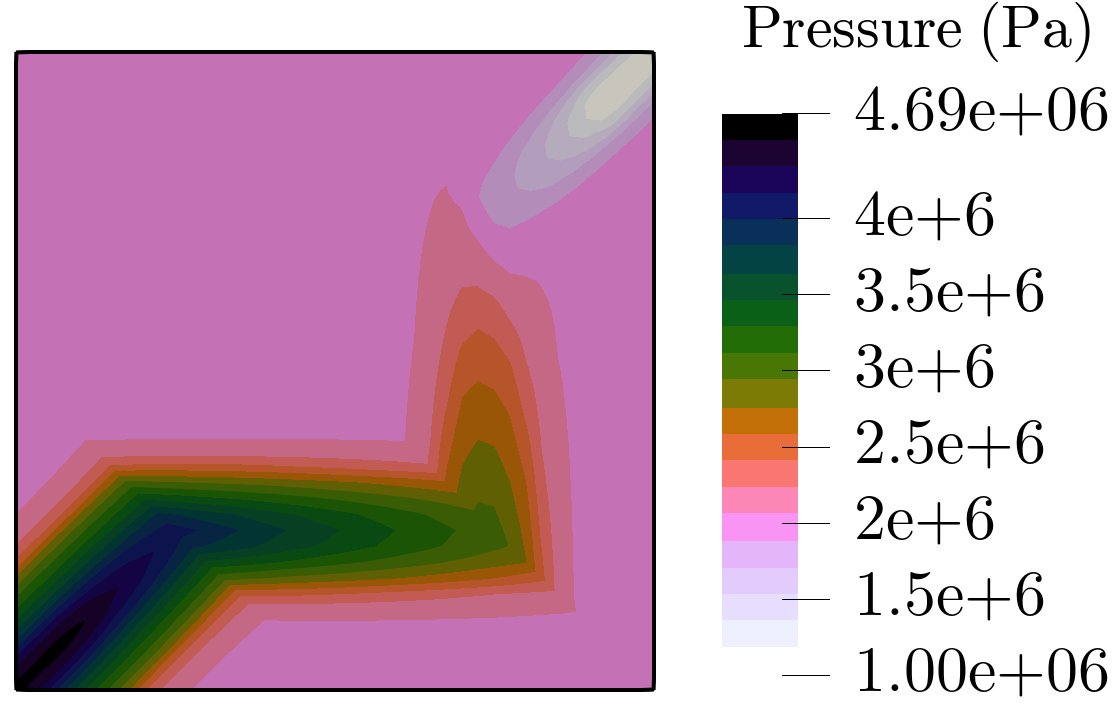}} 
        \hspace{.05cm}
    \subfigure[$t=5$ days; pressure \label{Fig:Q5_anisotropy_pres_c}]{
        \includegraphics[clip,scale=0.13,trim=0 0cm 0cm 0]{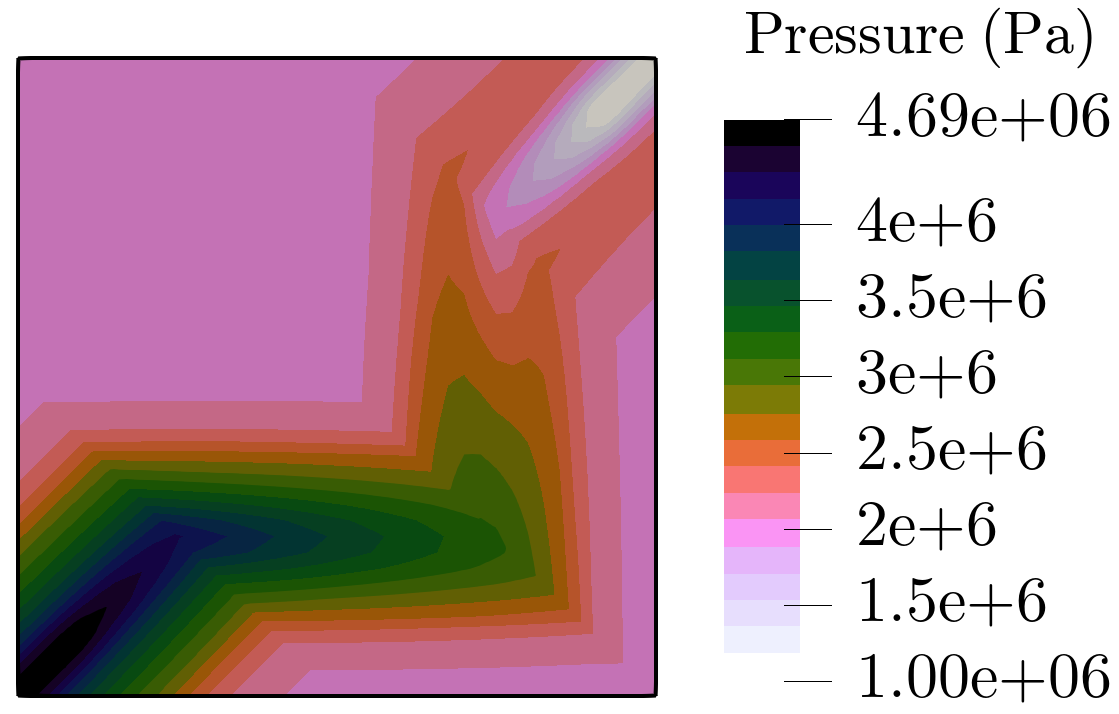}} \\
    \subfigure[$t=1.25$ days; velocity \label{Fig:Q5_anisotropy_vel_a}]{
        \includegraphics[clip,scale=0.12,trim=0 0cm 0cm 0]{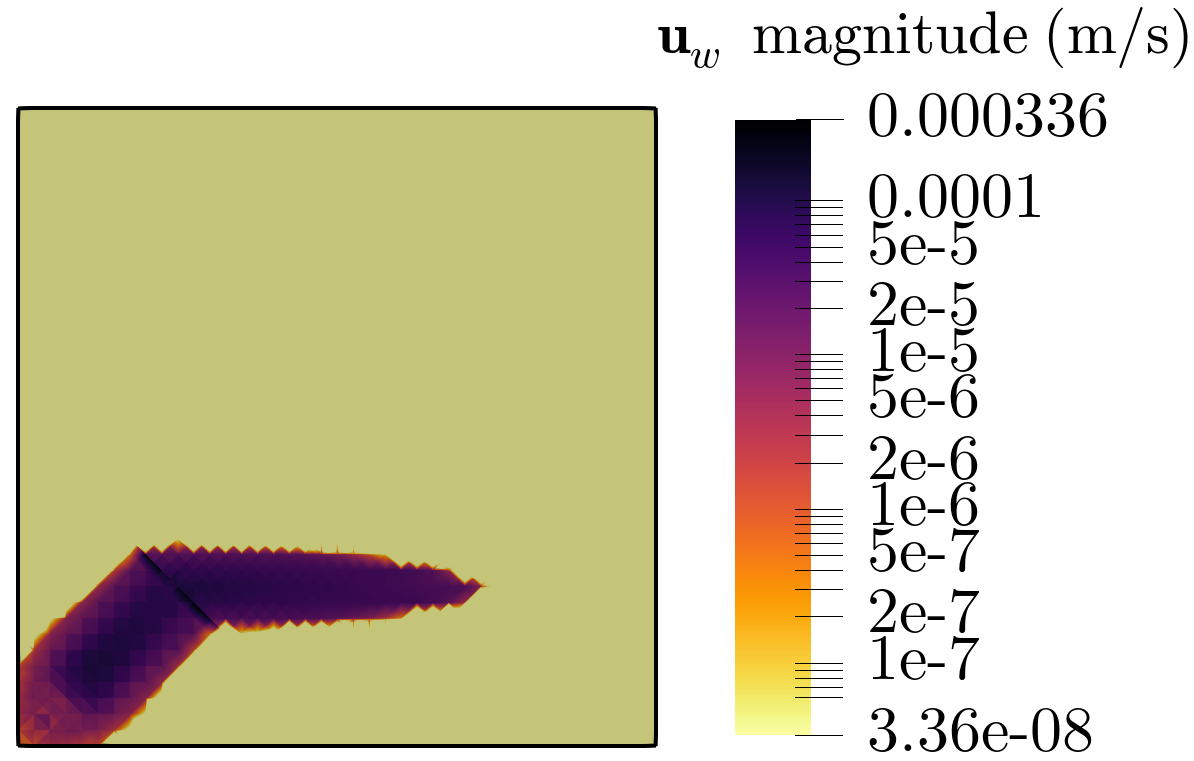}}
        \hspace{.05cm}
    \subfigure[$t=2.5$ days; velocity \label{Fig:Q5_anisotropy_vel_b}]{
        \includegraphics[clip,scale=0.12,trim=0 0cm 0cm 0]{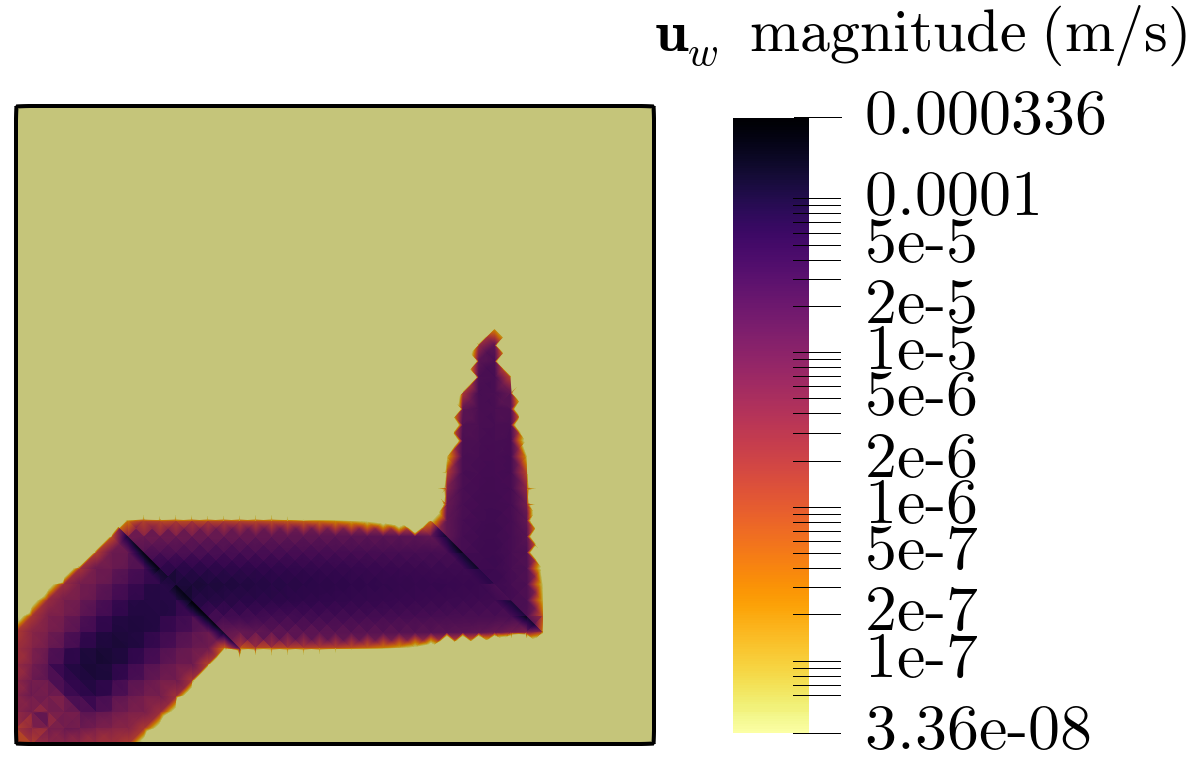}} 
        \hspace{.05cm}
    \subfigure[$t=5$ days; velocity \label{Fig:Q5_anisotropy_vel_c}]{
        \includegraphics[clip,scale=0.12,trim=0 0cm 0cm 0]{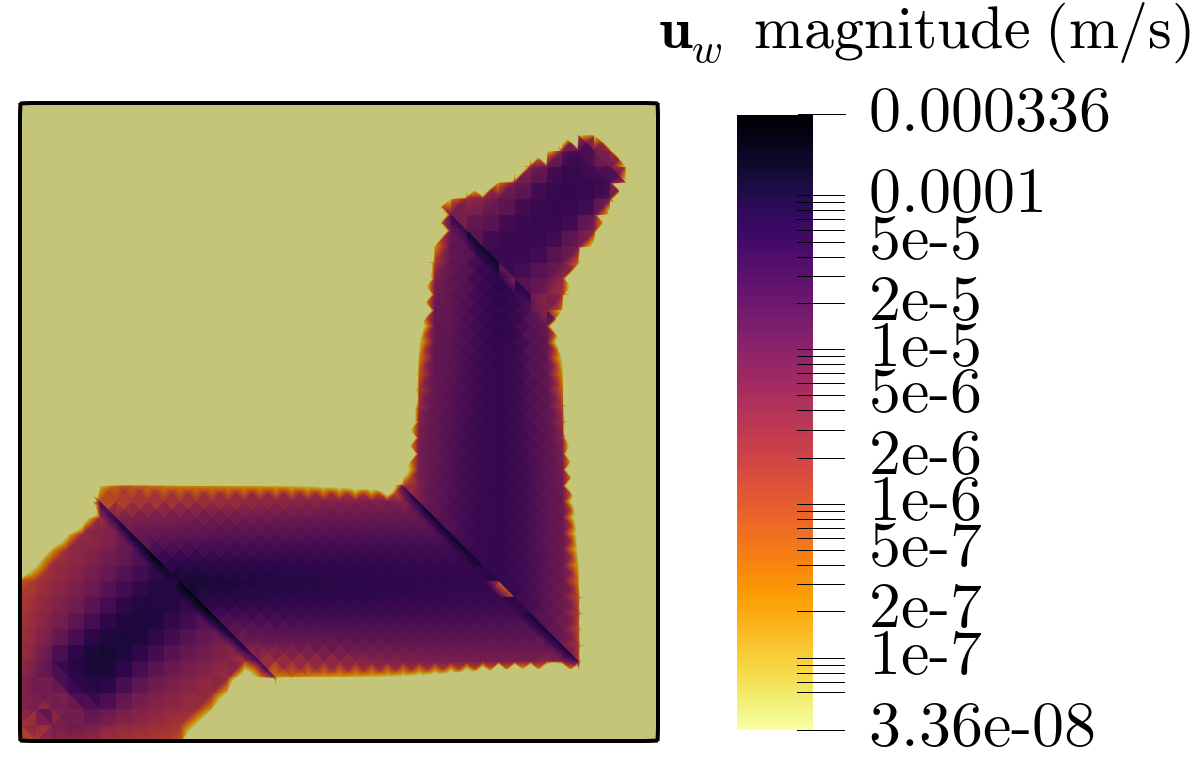}} 
        \caption{\textsf{Quarter five-spot problem in an anisotropic domain:}~this figure shows the computational results obtained from the limited DG scheme as time advances. The proposed scheme satisfies the maximum principle throughout the simulation.  The zigzag flow pattern confirms that solutions respect  the anisotropy of the domain.  
        \label{Fig:Q5_anisotropy_sol}}
\end{figure}

\label{sub:effect_of_anisotropy}
Finally, to examine the capability of the proposed limiting scheme to produce a correct and bound-preserving solution
for strongly anisotropic permeability fields, we present the results for the boundary-value problem illustrated in 
Figure \ref{fig:aniso_schematic}, which was adopted from \citep{galindez2020numerical,nikitin2014monotone}.
The permeability matrix  is defined as follows:
\begin{align}
    K=\begin{bmatrix}
        \cos \theta & -\sin \theta \\
        \sin \theta & \cos \theta
    \end{bmatrix}
    \begin{bmatrix}
        k_1 & 0 \\
        0 & k_2
    \end{bmatrix}
    \begin{bmatrix}
        \cos \theta & \sin \theta \\
        -\sin \theta & \cos \theta
    \end{bmatrix},
\end{align}
where principal permeabilities are set to $k_1=2.25\times10^{-12}$ and $k_2=2.25\times10^{-14}$. As shown in Figure \ref{fig:aniso_schematic} the permeability field is divided into four separate regions with distinct anisotropic $K$,
in which the angle $\theta$ is equal to 45 degree in the bottom left and upper right parts of the domain and alternate between 0 and 90 degree in the central region. The remaining parameters are the same as in Section \ref{sub:quarter_five_spot_problem}.
Figure \ref{Fig:Q5_anisotropy_sol} depicts the computational results computed with the proposed limited DG scheme at
three different time instances. It is clear that the limiting strategy honors the domain's heterogeneity and anisotropy and 
the channel flow with stair-case shape is captured. It should be also noted that no violation of maximum principle or spurious oscillations are obtained in the solutions.


\subsection{A note on the solver and scheme performance}%
\label{sub:solver}
We implement the proposed computational framework using the finite element capabilities in Firedrake Project \citep{rathgeber2016firedrake,McRae2016,Homolya2016,Homolya2017,homolya2017a} with GNU compilers. 
We resort to the MPI-based PETSc library \citep{petsc-user-ref,petsc-web-page,Dalcin2011} as the
linear algebra back-end to solve the nonlinear system. We use Newton’s method with
 step line search technique and set the relative convergence tolerance to $10^{-6}$.
At each time step, after the Newton solver convergence, we apply flux and slope limiters. 
Implementation of the flux limiter algorithm is discussed in Section \ref{Sub:Flux_limiter} 
and the global stopping criteria for all problem sets are taken as $\epsilon_1=\epsilon_2=10^{-6}$.
As for the slope limiter, we take advantage of the VertexBasedLimiter module embedded in the Firedrake project. 
All simulations are run on a single socket Intel i5-8257U node by utilizing a single MPI process.
Codes used to perform all experiments in this paper are publicly available at msarrafj/LimiterDG [2022] 
repository for reproducibility.

Table \ref{tab:performance} illustrates the Newton solver and flux limiter performance in terms of number of iterations. A few Newton iterations are needed at each time step for the convergence of either limited DG or unlimited DG approximations. It is evident that the limiters do not have a noticeable effect on the number of solver iterations. 
We also see that the number of flux limiter iterations is not significantly affected by the compressibility factors. However, 
the maximum number of iterations increases as heterogeneity, gravity, and anisotropy are added to the system. It should be noted that the reported number range for flux limiter is recorded throughout the simulation time and we observed that in fact for most time steps (over 85\% to 90\%), flux limiter iterations remain relatively small (less than 5 iterations).

\begin{table}
\caption{Number of nonlinear Newton iterations and flux limiter iterations per time step during the simulation}
\label{tab:performance}
\centering
\resizebox{1.0\textwidth}{!}{%
\begin{tabular}{ll@{\hskip 0.5in}l@{\hskip 0.5in}ll}
\Xhline{2\arrayrulewidth}\\[-0.9em]
\multicolumn{2}{c}{\multirow{2}{*}{Problem description}}  & Unlimited DG    & \multicolumn{2}{l}{Limited DG+FL+SL} \\
\cmidrule(l{-0.1in}r{0.42 in}){3-3} \cmidrule(l{-0.1in}r{0.05 in}){4-5}
\multicolumn{1}{c}{}                            &                       & Newton's iter. num. & Newton's iter. num. & FL iter. num. \\\\[-0.9em] \hline\\[-0.9em]
\textit{Sec 4.2 - Ex 1}  & homogeneous-bc-w/o gravity   & $3-5$     & $3-5$         & $2-8$           \\
\textit{Sec 4.2 - Ex 2}  & heterogeneous-bc-w/o gravity   & $3-5$     & $3-5$         & $2-11$          \\
\textit{Sec 4.2 - Ex 3}  & nonhomogeneous-bc-with gravity   & $3-7$     & $3-7$         & $3-30$          \\
\textit{Sec 4.3}    & homogeneous-wells-w/o gravity   & $4-5$     & $4-5$         & $4-28$          \\
\textit{Sec 4.5} & anisotropic-wells-w/o gravity   & $3-5$     & $3-5$         & $4-53$          \\\\[-0.9em]
\Xhline{2\arrayrulewidth}
\end{tabular}%
}
\end{table}

\section{Conclusions}

We have developed a numerical method that solves for primary unknowns the wetting phase saturation and pressure
of a compressible two-phase flows problem in a compressible rock matrix. A fully implicit discontinuous Galerkin scheme
is augmented with post-processing flux and slope limiters for the saturation.  The performance and accuracy of the method
is investigated for several benchmark problems including the quarter-five spot problem. Overshoot and undershoot
are completely eliminated throughout the whole simulation time. The impact of the limiters
on the local mass conservation is shown to be negligible.  The use of flux and slope limiters does not change
the number of Newton iterations compared to the case of unlimited DG. Numerical simulations show that the limited DG scheme
significantly improves the monotonicity of the saturation compared to the one obtained with the unlimited method. The 
limited DG method produces sharp saturation fronts with minimal numerical diffusion, and can handle anisotropic media.
The method is also shown to be robust for flows under gravitational forces.
%

\bibliographystyle{plainnat}

\end{document}